\input amstex
\catcode`\@=11

 \def\AMSTeXfeatures{\Plainheads 
   \let\current@vert=\AMS@vert}

 \def\Plainheads{\sh@ftdiam=0.05em
   \getlabeldims
   \let\vshaftfill=\plnvsolidfill
   \let\hshaftfill=\plnhsolidfill
   \let\th@rhead=\plnrhead
   \let\th@lhead=\plnlhead
   \let\th@dnhead=\plndnhead
   \let\th@uphead=\plnuphead}
 
 \def\glet{\global\let}

 \def\LaTeXfeatures{\catcode`\@=11
   \ifx\@clnwd\undefined \nol@g
      \input ltxcode.tex \dol@g \fi
   \ltxheads \let\current@vert=\new@vert
   \providelto \catcode`\@=\active}

 \def\nol@g{\def\wlog{\edef\garbage}}
 \def\dol@g{\let\wlog=\wl@g} \let\wl@g=\wlog
 \nol@g 

 \newbox\ltobox
 \def\providelto{{\setbox\z@=
   \hbox{$\to$}\minharrlen=\wd\z@
   \global\setbox\ltobox=\hbox{$\activeat>>>$}}
   \def\lto{\mathrel{\copy\ltobox}}}

 \def\ltxheads{\sh@ftdiam=\@wholewidth
   \getlabeldims
   \let\vshaftfill= \ltxvsolidfill
   \let\hshaftfill=\ltxhsolidfill
   \let\th@rhead=\ltxrhead
   \let\th@lhead=\ltxlhead
   \let\th@dnhead=\ltxdnhead
   \let\th@uphead=\ltxuphead}
 {\catcode`\@=\active
   \gdef@#1{\csname #1\string@at\endcsname}
   \glet\activeat=@}
 \def\def@#1{\expandafter\def\csname #1@at\endcsname}

 \def@>#1>#2>{\@rrow R{#1}{#2}}
 \def@<#1<#2<{\@rrow L{#1}{#2}}
 \def@ V#1V#2V{\@rrow V{#1}{#2}}
 \def@ A#1A#2A{\@rrow A{#1}{#2}}
 \def@/#1/#2/#3/{\@rrow{#1}{#2}{#3}}
 \def@.{\ifodd\row\ifmmode\noharrow
     \else\leavevmode.\spacefactor3000 \fi
   \else\novarrow\fi}
 \def@={\ifodd\row\harrow\hequalfill{}{}%
   \else\varrow\vequalfill{}{}\fi}
 \def@:#1{\ifx=#1\harrow\deffill{}{}%
   \else\leavevmode\null:#1\fi}
 \def@|{\current@vert}
  \def\AMS@vert{\varrow\vequalfill{}{}}
  \def\new@vert#1|#2|{\ifodd\row
   \let\nextarrow\vertexvarrow
   \else\let\nextarrow\varrow\fi
   \nextarrow\vshaftfill{#1}{#2}}
 \def@-{\ifmmode\let\next\hl@ne
   \else\let\next\AMSatdash \fi \next}
  \def\hl@ne#1-#2-{\harrow\hshaftfill{#1}{#2}}
  \def\AMSatdash{\let\next\relax\leavevmode
    \def\next@{\ifx\next-%
      \def\next-{\futurelet\next\nextii@}%
     \else\def\next{\hbox{-}}\fi\next}%
    \def\nextii@{\ifx\next-\def\next-{\hbox{---}}%
      \else\def\next{\hbox{--}}\fi\next}%
    \futurelet\next\next@}
 \def@(#1){\tweenarrows{#1}}
 \def@[#1]{\setsp@n#1\relax\activeat}
 \def\fiberbox{\hbox{$\vcenter{\hr@le\hbox{\vr@le
   \kern1ex\vbox{\kern1.2ex}\vr@le}\hr@le}$}}
  \def\hr@le{\hrule height \sh@ftdiam}
  \def\vr@le{\vrule width \sh@ftdiam}
 \def@+#1+#2+#3+{\ifodd\row \harrow{#1}{#2}{#3}%
   \else \varrow{#1}{#2}{#3}\fi}


 \def\Dnarrfill{\vequalfill\Dnhe@d}
 \def\Uparrfill{\Uphe@d\vequalfill}
 
 \def\ontofill{\rtarrfill\kern-0.3em 
   \th@rhead\kern 0.3em} 

 \def\rtarrfill{\hshaftfill\th@rhead}
 \def\ltarrfill{\th@lhead\hshaftfill}
 \def\dnarrfill{\vshaftfill\th@dnhead}
 \def\uparrfill{\th@uphead\vshaftfill}
 \def\hequalfill{\plnhfill=}
 \def\deffill{:\plnhfill=}
 \def\plnvextfill#1{\setbox\z@
   \hbox{\the\textfont3 #1}%
   \dimen@=\dp\z@\advance\dimen@\ht\z@
   \copy\z@ \kern-\dimen@ 
   \cleaders\copy\z@ \vfill
   \kern-\dimen@ 
   \box\z@}
 \def\plnhfill#1{$\m@th\mkern-1.5mu\mathord#1\mkern-6mu
    \cleaders\hbox{$\mkern-2mu\mathord#1\mkern-2mu$}\hfill
    \mkern-6mu\mathord#1\mkern-1.5mu$}
 \def\vequalfill{\plnvextfill{\char'167}}
 \def\plnvsolidfill{\plnvextfill{\char'077}}
 \def\plnhsolidfill{\plnhfill-}
 \def\ltxhsolidfill{\leaders\hrule height\topofshaft depth\botofshaft
   \hfill}
 \def\ltxvsolidfill{\leaders\vrule width\sh@ftdiam\vfill}
 \def\hdashfill{\hd@sh\wd@sh
   \xleaders \hbox{\wd@sh\hd@sh\wd@sh}\hfill
   \wd@sh\hd@sh}
 \def\vdashfill{\vd@sh\wd@sh
   \xleaders \vbox{\wd@sh\vd@sh\wd@sh}\vfill
   \wd@sh\vd@sh}
 \def\dashed{\ifinmeasureCD\else
    \ifodd\row\option{\let\hshaftfill=\hdashfill}%
   \else\option{\let\vshaftfill=\vdashfill}\fi\fi}


 \newdimen\CDstrutht  \newdimen\CDstrutdp
   \CDstrutht=0.875\baselineskip
   \CDstrutdp=0.375\baselineskip
 \newdimen\CDstrutlen \CDstrutlen=\CDstrutht
   \advance\CDstrutlen by \CDstrutdp

 \def\CDstrut{\vrule
   height \ifnum\row=1 \z@\else\CDstrutht \fi
   depth \ifnum\row=\numrows \z@ \else\CDstrutdp \fi
   width\z@}

 \newdimen\CDarrsurr \CDarrsurr=0.375em
 \newdimen\CDdashlen
    \CDdashlen= 0.1875\baselineskip
 \newdimen\CDvarrlen \CDvarrlen=1.5\baselineskip
 \newdimen\minharrlen 
  \setbox\z@\hbox{$\longrightarrow$} \minharrlen=\wd\z@
 \newdimen\minCDharrlen \minCDharrlen=2.5em 
\newdimen \minc@lwd
\def\findminc@lwd{\minc@lwd=2\CDarrsurr
  \advance\minc@lwd\minCDharrlen}

 \newdimen\sh@ftdiam


 \newdimen\labelsurr \labelsurr=1.25 em

\newcount\sp@ncnt \sp@ncnt=\@ne
\newcount\sp@ncnt@ \sp@ncnt@=\@ne
\newdimen\@rrwd \newdimen\@rrdp


 \def\adjustbot#1{\option{\advance\@rrdp#1\relax}}
 
\def\pushvertex#1{\global\p@shlen#1\relax
   \global\let\maybepush=\dopush}


 \newdimen\p@shlen \p@shlen=\z@

 
 \let\maybepush=\relax
 \def\dopush{\ifinmeasureCD 
   \advance\locdimen by -\p@shlen 
   \else\advance \@rrwd by -\p@shlen \fi 
   \global\let\maybepush=\relax \global\p@shlen=\z@\relax}


 \def\span@ne{\global\sp@ncnt=\@ne\relax}
 \def\setsp@n#1#2{\global\sp@ncnt=#1\relax
   \ifx\relax#2\relax\else\global\sp@ncnt@=#2\relax\fi}

 \def\plnrhead{\llap{$\rightarrow\mkern-1.5mu$}}
 \def\plnlhead{\rlap{$\mkern-1.5mu\leftarrow$}}

 \def\clap#1{\hbox to \z@{\hss #1\hss}}

 \def\plndnhead{\hbox{\the\textfont3 \char'171}}
 \def\plnuphead{\hbox{\the\textfont3 \char'170}}
 \def\Dnhe@d{\hbox{\the\textfont3 \char'177}}
 \def\Uphe@d{\hbox{\the\textfont3 \char'176}}

 \def\ltxrhead{\raise\@xisheight
   \llap{\smash{\@linefnt\@getrarrow(1,0)}}}
 \def\ltxlhead{\raise\@xisheight
   \rlap{\@linefnt\@getlarrow(-1,0)}}
 \def\ltxuphead{\setbox\z@=\rlap{%
   \kern\@halfwidth\@linefnt\char'66}%
   \copy\z@\kern-\ht\z@}
 \def\ltxdnhead{\setbox\z@=\rlap{%
   \kern\@halfwidth\@linefnt\char'77}%
   \ht\z@=\z@\box\z@}

 \def\wd@sh{\kern0.5\CDdashlen}
 \def\hd@sh{\vrule height\topofshaft depth\botofshaft
    width\CDdashlen}
 \def\vd@sh{\hrule height\CDdashlen
   depth\z@ width\sh@ftdiam}

\def\xylist{14{3434}13{2414}12{1723}%
  23{1413}34{1153}11{0867}43{0707}%
  32{0580}21{0414}31{0291}41{0}}
\newcount\tgtcnt@
\def\find@xyargs{\dimen@=\@rrdp
  \advance\dimen@ by \CDstrutlen
  \tgtcnt@=\dimen@ \dimen@=\@rrwd 
  \divide\dimen@ by \@m 
  \divide \tgtcnt@ by \dimen@ 
  \expandafter\testxy\xylist\relax
  \unitlength=\@xarg\@rrdp
  \divide\unitlength by\@yarg\relax}
\def\testxy#1#2#3{\ifnum\tgtcnt@>#3
    \@xarg=#1\relax \@yarg=#2\relax
    \let\next=\ignorerest
  \else\let\next\testxy\fi\next}
\def\ignorerest#1\relax{\relax}

\let\scalefactor=\@ne
\def\SWarrow{\find@xyargs\vector
  (-\@xarg,-\@yarg)\scalefactor\hskip-\wd\@linechar}
\def\NWarrow{\find@xyargs\vector
  (-\@xarg,\@yarg)\scalefactor\hskip-\wd\@linechar}
\def\NEarrow{\find@xyargs\vector
  (\@xarg,\@yarg)\scalefactor}
\def\SEarrow{\find@xyargs\vector
  (\@xarg,-\@yarg)\scalefactor}
\def\rightupline{\find@xyargs\@linelen=\scalefactor
     \unitlength\@sline}
\def\rightdownline{\find@xyargs\@yarg=-\@yarg\relax
     \@linelen=\scalefactor\unitlength\@sline}

\def\Sim{\ifodd\row\setbox\z@=\hbox{$\sim$}\dimen@=\ht\z@
 \advance\dimen@ by -\@xisheight
  \vbox{\box\z@\kern-\@xisheight\kern\dimen@}%
  \else\hbox{$\wr$}\fi}

%
\def\harrow#1#2#3{\inmeasureCDtrue\findminarrwd
  {#2}{#3}{\sp@ncnt\minharrlen}\inmeasureCDfalse\span@ne
  \mathrel{\hbox{\options\hplace{#1}\ulabel{#2}\dlabel{#3}}}}

\def\noharrow{\harrow\hfill{}{}}
\def\vertexvarrow#1#2#3{\findarrdp \@rrwd=\z@ \setsp@n\@ne\@ne
  \vbox to \z@{\kern-1.2\CDstrutht
  \rlap{\options\vplace{#1}\llabel{#2}\rlabel{#3}}\vss}}

\newif\ifinmeasureCD
\def\measurelabel#1{\setbox\z@
  \hbox{$\scriptstyle#1\kern\labelsurr$}%
  \ifdim\wd\z@>\@rrwd \@rrwd=\wd\z@\fi}
\def\findminarrwd#1#2#3{\@rrwd=#3\relax
   \measurelabel{#1}\measurelabel{#2}}
\def\findCDarrwd#1#2{\@rrwd=\minCDharrlen
   \measurelabel{#1}\measurelabel{#2}%
  }

\newcount\row \row=\@ne \newcount\col \col=\@ne 
 \newcount\numrows 
\numrows=\@ne
 \newcount\numcols
\newcount\arrspan \newdimen\vrtxhalfwd  \newbox\tempbox

\def\DANABUG{\advance\col by \@ne
 \@rrwd=\minCDharrlen
  \advance\@rrwd by \vrtxhalfwd
  \advance\@rrwd by \CDarrsurr
  \ifnum\col>\numcols \numcols=\col
     \newlocdimen{col\the\col}\locdimen=\@rrwd 
  \else \ifdim\@rrwd>\c@l \c@l=\@rrwd\fi\fi}

\def\drop#1\\{
  \findvrtxhalfsum\DANABUG\advance\row by 2 \measureinit}

\def\measureinit{\col=\@ne \vrtxhalfwd=-\CDarrsurr\arrspan=\@ne\@rrwd=\z@
   \setbox\tempbox=\hbox\bgroup$}
\def\measure{
  \let\harrow\measureCDarrow
  \let\CDCR=\measureCR 
   \findminc@lwd 
  \inmeasureCDtrue
  \row=\@ne \numcols=\z@ \measureinit}

\def\endmeasure{\findvrtxhalfsum\DANABUG
  \numrows=\row 
  \inmeasureCDfalse}




\def\newlocdimen#1{\advance\dimenc@unt by \@ne
  \ifnum\dimenc@unt<\insc@unt
     \else\errmessage{No room for the CD}\fi
  \dimendef\locdimen=\dimenc@unt
  \expandafter\dimendef\csname#1\endcsname=\dimenc@unt}

 \def\r@wc@l{\csname row\the\row col\the\col\endcsname}
 \def\c@l{\csname col\the\col\endcsname}

 \def\findvrtxhalfsum{$\egroup
  \newlocdimen{row\the\row col\the\col}
  \locdimen=\vrtxhalfwd 
  \vrtxhalfwd=0.5\wd\tempbox 
  \advance\vrtxhalfwd by \CDarrsurr
  \advance\locdimen by \vrtxhalfwd 
  \advance\@rrwd by \locdimen 
  \maybepush
  \divide\@rrwd by \arrspan\relax
  \ifdim\@rrwd<\minc@lwd
    \ifnum\col>\@ne \@rrwd=\minc@lwd\fi \fi
  \loop 
    \ifnum\col>\numcols \numcols=\col
       \newlocdimen{col\the\col}
       \locdimen=\@rrwd 
    \else \ifdim\@rrwd>\c@l \c@l=\@rrwd\fi \fi
   \ifnum\arrspan>\@ne
      \advance\arrspan by -1 \advance\col by \@ne
  \repeat }

 \def\measureCDarrow#1#2#3{\findvrtxhalfsum
   \arrspan=\sp@ncnt\relax\global\sp@ncnt=1\relax
   \advance\col by \@ne
   \findCDarrwd{#2}{#3}%
   \setbox\tempbox=\hbox\bgroup$}

 \newcount\dr@tn \dr@tn=\z@
 \def\locate#1:#2{\ifinmeasureCD\else
   \count@=-#1
   \multiply\count@ by 2
   \advance\count@ by #2
   \dimen@=\count@\@rrwd
   \ifnum\dr@tn=\@ne\relax \else\dimen@=-\dimen@ \fi
   \dimen@i=\@rrdp
   \ifnum\dr@tn>\z@\advance\dimen@i by \CDstrutlen \fi
   \dimen@i=\count@\dimen@i
   \count@=#2 \multiply\count@ by 2
   \divide\dimen@ by \count@
   \divide\dimen@i by \count@
   \lift\dimen@i\nudge\dimen@\fi}

\def\betweenCDrows{\advance\row by \@ne \col=\@ne
\options}


\def\hbegin{\hbox\bgroup\kern\c@l \kern-\r@wc@l$}
\def\hend{$\glet\maybepush\relax \CDstrut\egroup}
\def\vbegin{\setbox\tempbox=\hbox\bgroup$}
\def\vend{$\egroup\ht\tempbox=\z@\dp\tempbox\CDvarrlen
  \box\tempbox}
\def\setCD{\let\harrow=\setCDarrow
  \let\CDCR=\setCR 
  \row=\@ne \col=\@ne \hbegin}
\let\endsetCD=\hend 

\def\findarrwd{\@rrwd=\z@ \count@=\col \advance\count@ by\sp@ncnt
  \loop\ifnum\count@>\col \advance\count@ by -1
      \advance\@rrwd by\csname col\the\count@\endcsname\repeat}
\def\setCDarrow#1#2#3{\kern\CDarrsurr\advance\col by \@ne
  \findarrwd \advance\@rrwd by -\r@wc@l  
  \@rrdp=\z@ 
  \maybepush
  \advance\col by -\@ne \advance\col by \sp@ncnt \span@ne
  \hbox to \@rrwd{\options
   \@rrwd=\scalefactor\@rrwd\hss
   \hplace{#1}\ulabel{#2}\dlabel{#3}\hss}%
   \kern\CDarrsurr}

\newdimen\labspacei 
\newdimen\labspaceii 

\newdimen\@xisheight
  \@xisheight=\the\fontdimen22\textfont2
\newdimen\labelskip
  \labelskip=\the\fontdimen10\textfont3 
\newdimen\topofshaft
\newdimen\botofshaft
\newdimen\botofulabel
\newdimen\topofdlabel
\def\getlabeldims{
  \topofshaft=0.5\sh@ftdiam
  \botofshaft=\topofshaft
  \advance\topofshaft by \@xisheight  
  \advance\botofshaft by -\@xisheight  
  \botofulabel=\topofshaft
  \advance\botofulabel by \labelskip
  \topofdlabel=\botofshaft
  \advance\topofdlabel by \labelskip}

\def\ulabel{\ifnum\row=\@ne\let\next\ulabeli
   \else\let\next\ulabellap\fi\next}
\def\ulabeli#1{\vbox{
  \clap{\kern-\@rrwd$\scriptstyle#1$}%
  \kern\botofulabel}\maybeoffset}
\def\ulabellap#1{\vbox to \z@{\vss
  \clap{\kern-\@rrwd$\scriptstyle#1$}%
  \kern\botofulabel}\maybeoffset}
\def\dlabel{\ifnum\row=\numrows\let\next\dlabeli
   \else\let\next\dlabellap\fi\next}
\def\dlabeli#1{\vtop{\kern\topofdlabel
  \clap{\kern-\@rrwd$\scriptstyle#1$}%
  }\maybeoffset}
\def\dlabellap#1{\vbox to \z@{\kern\topofdlabel
  \clap{\kern-\@rrwd$\scriptstyle#1$}%
  \vss}\maybeoffset}
\def\rlabel#1{\vbox to \z@{\vss
  \rlap{\kern\labelskip$\scriptstyle#1$}%
  \vss\kern-\@rrdp}\maybeoffset}
\def\llabel#1{\vbox to \z@{\vss
  \llap{$\scriptstyle#1$\kern\labelskip}%
  \vss\kern-\@rrdp}\maybeoffset}
\def\swlabel#1{\vtop{\kern0.5\@rrdp
  \llap{$\scriptstyle#1$\kern\labelskip\kern-0.5\@rrwd}
  }\maybeoffset}
\def\nwlabel#1{\vbox{
  \llap{$\scriptstyle#1$\kern\labelskip\kern-0.5\@rrwd}%
  \kern-0.5\@rrdp}\maybeoffset}
\def\selabel#1{\vtop{\kern0.5\@rrdp
  \rlap{\kern0.5\@rrwd\kern\labelskip$\scriptstyle#1$}%
  }\maybeoffset}
\def\nelabel#1{\vbox{
  \rlap{\kern0.5\@rrwd\kern\labelskip$\scriptstyle#1$}%
  \kern-0.5\@rrdp}\maybeoffset}
\def\cplace#1{\vbox to \z@{\vss
  \clap{$#1$\kern-\@rrwd}%
  \kern-\@rrdp\vss}\maybeoffset}
\def\hplace#1{\hbox to \@rrwd{#1}\maybeoffset}
\def\vplace#1{\clap{\vbox to \z@{#1\kern-\@rrdp}}\maybeoffset}

\newdimen\nudgeamount \nudgeamount=\z@
\newdimen\liftamount \liftamount=\z@
\let\maybeoffset\relax
\newbox\offsetbox \newdimen\lastheight
\def\dooffset{
  \setbox\offsetbox=\lastbox \lastheight=\ht\offsetbox 
  \setbox\offsetbox=\vbox{\kern-\liftamount\box\offsetbox}%
  \ht\offsetbox=\lastheight
  \kern\nudgeamount\box\offsetbox\kern-\nudgeamount
  \global\nudgeamount=\z@ \global\liftamount=\z@
  \glet\maybeoffset=\relax}
\def\nudge#1{\ifinmeasureCD\else
  \global\advance\nudgeamount#1\relax
  \global\let\maybeoffset\dooffset\fi}
\def\lift#1{\ifinmeasureCD\else
  \global\advance\liftamount#1\relax
  \global\let\maybeoffset\dooffset\fi}

\def\findarrdp{\@rrdp=\CDvarrlen
  \ifnum\sp@ncnt@>1
    \advance\@rrdp by \CDstrutlen
    \multiply\@rrdp by \sp@ncnt@
    \advance\@rrdp by -\CDstrutlen \fi
 }

\def\varrow#1#2#3{\ifnum\sp@ncnt>\@ne 
     \sp@ncnt@=\sp@ncnt\relax\fi
  \findarrdp \@rrwd=\z@ 
  \kern\c@l
   \hbox to \z@{\options
   \@rrdp=\scalefactor\@rrdp
    \hss\vplace{#1}\llabel{#2}\rlabel{#3}\hss}%
  \global\advance\col by \@ne \setsp@n\@ne\@ne
  }

\def\novarrow{\varrow\vfill{}{}}

\def\tweenarrows#1{\findarrwd \findarrdp \setsp@n\@ne\@ne
  \rlap{\options\cplace{#1}}}

\def\usarrow #1#2#3{\dr@tn=\@ne
  \findarrwd \findarrdp \setsp@n\@ne\@ne 
  \rlap{\options\cplace{#1}\nwlabel{#2}\selabel{#3}}%
  \dr@tn=\z@}
\def\dsarrow #1#2#3{\dr@tn=\tw@
  \findarrwd \findarrdp \setsp@n\@ne\@ne 
  \rlap{\options\cplace{#1}\swlabel{#2}\nelabel{#3}}%
  \dr@tn=\z@}
 \def\@rrow#1{\csname #1@rrow\endcsname}
 \def\R@rrow{\harrow \rtarrfill}
 \def\L@rrow{\harrow \ltarrfill}
 \def\V@rrow{\varrow \dnarrfill}
 \def\A@rrow{\varrow \uparrfill}
 \def\SE@rrow{\dsarrow \SEarrow}
 \def\NW@rrow{\dsarrow \NWarrow}
 \def\SW@rrow{\usarrow \SWarrow}
 \def\NE@rrow{\usarrow \NEarrow}
 \def\DS@rrow{\dsarrow \dnslope}
 \def\US@rrow{\usarrow \upslope}
 \def\upslope{\find@xyargs
       \@linelen=\unitlength\@sline}
 \def\dnslope{\find@xyargs\@yarg=-\@yarg\relax
       \@linelen=\unitlength\@sline}

\newtoks\optionlist 
\optionlist={}
\let\options\relax
\def\dooptions{\the\optionlist\global\optionlist={}%
  \glet\options=\relax}
\def\option#1{\ifinmeasureCD\else
  \glet\options=\dooptions
  \global\optionlist=\expandafter{\the\optionlist\relax#1}\fi}
\def\wider#1{\ifinmeasureCD\else
   \option{\advance\@rrwd by #1}\fi}
\def\deeper#1{\ifinmeasureCD\else
   \option{\advance\@rrdp by #1}\fi}


{\def\\{\global\let\sptoken= }\\ }

\def\CR{\futurelet\nexttok\testCR}
\def\testCR{\ifx\nexttok\sptoken
   \let\next\eatspaceCR\else\let\next\CDCR\fi\next}
\def\eatspaceCR#1 {\CR}
\def\measureCR{\ifx\nexttok\endmeasure\let\nextCR\relax
    \else\let\nextCR\drop\fi\nextCR}
\def\setCR{\ifodd\row
  \ifx\nexttok\endsetCD\else\hend\betweenCDrows\vbegin\fi
  \else\vend\betweenCDrows\hbegin\fi}

\countdef\dimenc@unt=11
\def\CD#1\endCD{
   \begingroup\let\\=\CR
  \m@th\offinterlineskip
   \measure#1\endmeasure\null\,\vcenter{\setCD#1\endsetCD}\,
   \endgroup
    }

\ifx\@clnwd\undefined \nol@g\else\catcode`\ =14\relax\fi
 \font\@linefnt=line10 
 \newcount\@tempcnta
 \newcount\@tempcntb
 \newdimen\@tempdima
 \newdimen\@tempdimb
 \newdimen\@wholewidth
 \newdimen\@halfwidth
   \@wholewidth\fontdimen8\@linefnt \@halfwidth .5\@wholewidth
 \newdimen\unitlength
 \newcount\@xarg
 \newcount\@yarg
 \newcount\@yyarg
 \newbox\@linechar
 \newdimen\@linelen
 \newdimen\@clnwd
 \newdimen\@clnht
 \newif\if@negarg
 
 \def\@whilenoop#1{}

 \def\@whiledim#1\do #2{\ifdim #1\relax#2\@iwhiledim{#1\relax#2}\fi}

 \def\@iwhiledim#1{\ifdim #1\let\@nextwhile=\@iwhiledim 
         \else\let\@nextwhile=\@whilenoop\fi\@nextwhile{#1}}

 \def\@sline{\ifnum\@xarg< 0 \@negargtrue \@xarg -\@xarg \@yyarg -\@yarg
   \else \@negargfalse \@yyarg \@yarg \fi
 \ifnum \@yyarg >0 \@tempcnta\@yyarg \else \@tempcnta -\@yyarg \fi
 \ifnum\@tempcnta>6 \@badlinearg\@tempcnta0 \fi
 \ifnum\@xarg>6 \@badlinearg\@xarg 1 \fi
 \setbox\@linechar\hbox{\@linefnt\@getlinechar(\@xarg,\@yyarg)}%
 \ifnum \@yarg >0 \let\@upordown\raise \@clnht\z@
    \else\let\@upordown\lower \@clnht \ht\@linechar\fi
 \@clnwd=\wd\@linechar
 \if@negarg \hskip -\wd\@linechar \def\@tempa{\hskip -2\wd\@linechar}\else
      \let\@tempa\relax \fi
 \@whiledim \@clnwd <\@linelen \do
   {\@upordown\@clnht\copy\@linechar
    \@tempa
    \advance\@clnht \ht\@linechar
    \advance\@clnwd \wd\@linechar}%
 \advance\@clnht -\ht\@linechar
 \advance\@clnwd -\wd\@linechar
 \@tempdima\@linelen\advance\@tempdima -\@clnwd
 \@tempdimb\@tempdima\advance\@tempdimb -\wd\@linechar
 \if@negarg \hskip -\@tempdimb \else \hskip \@tempdimb \fi
 \multiply\@tempdima \@m
 \@tempcnta \@tempdima \@tempdima \wd\@linechar \divide\@tempcnta \@tempdima
 \@tempdima \ht\@linechar \multiply\@tempdima \@tempcnta
 \divide\@tempdima \@m
 \advance\@clnht \@tempdima
 \ifdim \@linelen <\wd\@linechar
    \hskip \wd\@linechar
   \else\@upordown\@clnht\copy\@linechar\fi}
 
 \def\@getlinechar(#1,#2){\@tempcnta#1\relax\multiply\@tempcnta 8
 \advance\@tempcnta -9 \ifnum #2>0 \advance\@tempcnta #2\relax\else
 \advance\@tempcnta -#2\relax\advance\@tempcnta 64 \fi
 \char\@tempcnta}
 
 \def\vector(#1,#2)#3{\@xarg #1\relax \@yarg #2\relax
 \@tempcnta \ifnum\@xarg<0 -\@xarg\else\@xarg\fi
 \ifnum\@tempcnta<5\relax
 \@linelen=#3\unitlength
 \ifnum\@xarg =0 \@vvector 
   \else \ifnum\@yarg =0 \@hvector \else \@svector\fi
 \fi
 \else\@badlinearg\fi}
 
 \def\@svector{\@sline
 \@tempcnta\@yarg \ifnum\@tempcnta <0 \@tempcnta=-\@tempcnta\fi
 \ifnum\@tempcnta <5
   \hskip -\wd\@linechar
   \@upordown\@clnht \hbox{\@linefnt  \if@negarg 
   \@getlarrow(\@xarg,\@yyarg) \else \@getrarrow(\@xarg,\@yyarg) \fi}%
 \else\@badlinearg\fi}
 
 \def\@getlarrow(#1,#2){\ifnum #2 =\z@ \@tempcnta='33\else
 \@tempcnta=#1\relax\multiply\@tempcnta \sixt@@n \advance\@tempcnta
 -9 \@tempcntb=#2\relax\multiply\@tempcntb \tw@
 \ifnum \@tempcntb >0 \advance\@tempcnta \@tempcntb\relax
 \else\advance\@tempcnta -\@tempcntb\advance\@tempcnta 64
 \fi\fi\char\@tempcnta}
 
 \def\@getrarrow(#1,#2){\@tempcntb=#2\relax
 \ifnum\@tempcntb < 0 \@tempcntb=-\@tempcntb\relax\fi
 \ifcase \@tempcntb\relax \@tempcnta='55 \or 
 \ifnum #1<3 \@tempcnta=#1\relax\multiply\@tempcnta
 24 \advance\@tempcnta -6 \else \ifnum #1=3 \@tempcnta=49
 \else\@tempcnta=58 \fi\fi\or 
 \ifnum #1<3 \@tempcnta=#1\relax\multiply\@tempcnta
 24 \advance\@tempcnta -3 \else \@tempcnta=51\fi\or 
 \@tempcnta=#1\relax\multiply\@tempcnta
 \sixt@@n \advance\@tempcnta -\tw@ \else
 \@tempcnta=#1\relax\multiply\@tempcnta
 \sixt@@n \advance\@tempcnta 7 \fi\ifnum #2<0 \advance\@tempcnta 64 \fi
 \char\@tempcnta}
\catcode`\ =10

\dol@g 
\catcode`\@=\active
\LaTeXfeatures

\documentstyle{amsppt}
\magnification=1200   \NoBlackBoxes
\define\Ker{\operatorname{Ker}}
\define\Bool{\operatorname{Bool}}
\define\Pol{\operatorname{Pol}}   \define \a{\alpha}
\define \be{\beta}  \define \dl{\delta} \define
\g{\gamma} \define \G{\Gamma}  
\define \Om{\Omega} \define \om{\omega}  
  \define \vp{\varphi}

    \define \iy{\infty}

\define\ol{\overline}
\define \do{\dots}

\define\pt{\Phi\Theta}
\define\vare{\varepsilon}

  \define \ov{\overline}

 \define \ep{\endproclaim}

  \define \bk{\bigskip}

\define \1{^{-1}} \define \2{^{-2}}

\def\capl{\operatornamewithlimits{\bigcap}\limits}

\define \Var{\operatorname{Var}}
\define \Grp{\operatorname{Grp}}
 \define \Alv{\operatorname{Alv}} \define
\Cl{\operatorname{Cl}}
\define \Hom{\operatorname{Hom}}  \define
\End{\operatorname{End}}     
      
\define\Th{\Theta}

\define\Val{\operatorname{Val}}

\define\Sub{\hbox{\rm Sub}}
\define\Log{\hbox{\rm Log}}
\define\Hal{\hbox{\rm Hal}}

\define \e{\varepsilon}
\def\snu{{\mathop\nu\limits^{\simeq}}}
\def\pr{^{\prime}}

\define\pmf{\par\medpagebreak\flushpar}
\define\knowf{Know$_{\Phi \Theta}(f)$}
\define\em{{_{\Phi \Theta}}}
\def\know{Know$_{\Phi \Theta}$}

\def\b{\beta}

 \define \Int{\operatorname{Int}}
\define \Aut{\operatorname{Aut}}
\vsize=7.50in

\baselineskip 12pt

\parindent 20pt
\NoBlackBoxes

\define\C{\Bbb C}

\def\pmf{\par\medpagebreak\flushpar}

.
\bigskip
\bigskip
\bigskip
\bigskip
\bigskip
\bigskip

\baselineskip 12pt
\topmatter \rightheadtext{} \leftheadtext{B.
Plotkin}
\title Seven Lectures on the Universal Algebraic Geometry
\endtitle
\endtopmatter

\bigskip
\bigskip
\centerline{{\smc B. PLOTKIN}\footnote{Partially supported by the Edmund Landau Center for Research in
Mathematical analysis and Related Areas, sponsored by the Minerva
Foundation (Germany)}}
\bigskip
\centerline{\smc Institute of Mathematics, Hebrew University}
\bigskip
\centerline{\smc Jerusalem, Israel}
\bigskip
\centerline{\smc borisov\@math.huji.ac.il}
\bigskip

\comment
\author
B. Plotkin,\\

Institute of Mathematics, Hebrew University,\\

Jerusalem, Israel,\\

 borisov\@math.huji.ac.il

\endauthor

\endcomment

\newpage

\centerline{\smc Contents}

\bigskip
\pmf

\pmf
{\smc Preface}

\pmf
{Lecture 1. \ \smc What is the Universal Algebraic Geometry}

\pmf
{Lecture 2. \ \smc Algebraic sets and algebraic varieties}

\pmf
{Lecture 3. \ \smc Geometric properties and geometric relations of algebras}

\pmf
{Lecture 4. \ \smc Algebraic geometry in group-based algebras}

\pmf
{Lecture 5. \ \smc Isomorphisms of categories of algebraic sets}

\pmf
{Lecture 6. \ \smc Algebraic geometry in the First Order Logic (FOL)}

\pmf {Lecture 7. \ \smc Databases and knowledge bases. Geometrical
aspect}

\pmf
{\smc Bibliography}


\newpage

\centerline{\smc Preface}
\bigskip
\bigskip

These lectures are devoted to the subject I've got interested
in during my work in Hebrew University. They are naturally connected with my
previous works on universal algebra and algebraic logic in database
theory. On the other hand, an essential influence on them is provided
by talks with E.Rips and Z.Sela and discussions on algebraic geometry
in a free group.

A great role for me play communications with A.Bokut, I.Dolgachev,
B.Ku\-nyav\-skii, R.Lypianskii,\  G.Mashevitzky, \ A.Miasnikov,
A.Mikhalev, A.Ol\-shan\-ski,  E.Plot\-kin,  V.Re\-meslen\-nikov,
A.Shmelkin, N.Vavilov, E.Winberg, E.Zelmanov and with many other
colleagues.

I would like to distinguish specially A.Berzins, who, in particular,
completely studied the situation in the classical variety
$Var-P$ (see 5.3), and whose key idea gives rise to the general
theorems 5 and $5'$ from the same Section.

Lecture 7 returns me to database and knowledge base theory and
relies on the joint work with T.Plotkin.

This text was prepared for lectures in USA and Canada in
September-October, 2000 and does not contain, as a rule, proofs of
results. On the other hand, it contains a list of problems of
different  levels of difficulty.

I am very grateful to Hebrew University for the support during
preparation of the lectures.

The bibliography {\it is not complete and consists  of the
works which have direct and indirect relation to the topic of
the lectures}. Most of the proofs can be found there.

The lectures are mostly prepared by E.Plotkin. I am also very
grateful to Mrs. M.Beller for typing.

I would like to add that I thank circumstances for bringing me to
 the problems I describe here.

\newpage

\baselineskip 12pt
\topmatter \rightheadtext{Lecture 1}
\leftheadtext{B. Plotkin}
\title Lecture 1\\
\quad\\
What is the Universal Algebraic Geometry\endtitle
\endtopmatter

\centerline{\smc Contents}
\bigskip
\roster\item "1." {\smc Generic view}
\newline
\item"2." {\smc
Classical algebraic geometry from the point of view of universal
algebraic geometry}
\newline
\item"3." {\smc Some information from universal algebra}
\newline
\item"4." {\smc Algebras with  a fixed algebra of constants}
\newline
\item"5." {\smc Special categories}
\endroster

\newpage

\subheading{1. Generic view}

Universal algebraic geometry is a system of algebraic
concepts and
problems going back to classical algebraic geometry.
It is mainly algebra, not quite geometry.
However, this algebra inherits some geometric intuition; preserves a
certain geometric
characteristic which is not necessarily explicit.

The key idea of universal algebraic geometry is that every variety of
algebras $\Th$ has its
own algebraic geometry.
It varies if the variety $\Th$ and an algebra $H\in\Th$ vary, but the
core notions are common.
We look at algebras $H\in\Th$ from the point of view of algebraic geometry
in $\Th$. On the other hand, supposing that an algebra $H$ is given,
we consider algebraic varieties and their invariants.
This is what is used to study in the classical algebraic geometry.

Universality of geometry means that the variety $\Theta$ can be
an arbitrary fixed variety of algebras. For every algebra $H\in\Theta$
we consider the category of algebraic sets $K_{\Theta}(H)$. This is
a geometric invariant of the algebra $H\in\Theta$. We are looking for
conditions which provide the categories $K_{\Theta}(H_1)$ and
$K_{\Theta}(H_2)$ to be isomorphic or equivalent for different $H_1$,
$H_2$ in $\Theta$. This is one of the main problems of the theory.
For the whole variety $\Theta$ its geometrical invariant is presented
by the special category $K_\Theta$. All $K_{\Theta}(H)$ are the subcategories
in $K_\Theta$. There is a problem of isomorphism or equivalence of
the categories $K_{\Theta_1}$ and $K_{\Theta_2}$ for the different
varieties $\Theta_1$ and $\Theta_2$.

In the above setting, the classical algebraic geometry occupies a
distinguished place.
It is based on the variety of commutative and associative rings with unity
and  is connected
with the corresponding varieties of algebras over fields.
Every such variety over a field $P$ is denoted by $\Var$-$P$ and is called
a classical variety
over $P.$

If we replace the variety of associative and commutative  algebras by a
variety of
associative but not necessarily commutative algebras over different $P$,
then we come
to non-commutative
algebraic geometry.

On the other hand, in group theory there is a strong interest to
investigate equations in groups and, especially, in free groups.
The corresponding results give rise to the sketch of the special
algebraic geometry. The most interesting results are obtained for
groups with a given group of constants $G$.
 The corresponding
variety is denoted by $\Grp$-$G.$

Something similar is known for semigroups and other algebraic systems.

Now, one can say that universal algebraic geometry stands in the same
position with respect to
all these special theories as does universal algebra with
respect to groups,
rings, semigroups, etc.

To every variety of algebras $\Theta$ its elementary logic (i.e.,
first order logic) is associated. We consider algebraic logic in
$\Theta$. The elementary properties of an algebra $H\in\Theta$
are considered from the point of view of logic
in $\Theta$. Geometrical
and logical properties of algebras in $\Theta$ are often well
interacted.  The same is valid for geometrical and logical
relations between different algebras in $\Theta$.

Besides, we will consider algebraic geometry in first order logic, which
generalizes usual
equational algebraic geometry.
This outlet to logic is motivated geometrically and
also by various applications, in particular,
in computer science.

\subheading{2. Classical algebraic geometry from the point of view of
universal algebraic
geometry}

Let us fix an infinite ground field $P.$
This field generates the variety of algebras which was denoted by
$\Var$-$P,$ i.e., the
variety of all commutative and associative algebras with unity over the
given $P.$
The polynomial algebra $P[X]=P[x_1,\dots,x_n]$ is the free algebra in
$\Var$-$P$ over the set
$X.$

A field $L$ is an extension of the field $P$ if there is an injection $P\to L.$
Every such field $L$ is also an algebra in $\Var$-$P$.

To every map $\mu: X\to L$ there is a $P$-homomorphism $\mu:
P[X]\to L$ in one-to-one correspondence. Namely, if
$f=f(x_1,\dots, x_n)\in P[X]$, then
$f^\mu=f(x_1^\mu,\dots,x_n^\mu)$ is an element in $L$. The affine
space $L^{(n)}=L^X$ can be viewed as the set of all homomorphisms
$\Hom(P[X],L)$. A  point $a=(a_1,\dots,a_n)\in L^{(n)}$ is
identified with the homomorphism $\mu=\mu_a:P[X]\to L $ by the
rule $\mu(x_i)=a_i$. Then $f^\mu=f(a)$ and the point $a$ is a root
of the polynomial $f,$ if $f\in \Ker \mu.$

We consider affine algebraic sets over the given $L.$
These sets are defined by the systems of equations with coefficients from
$P$ and solutions
in the field $L.$
Since the algebraic sets are supposed to be subsets in $L^{(n)}$, we can
consider them  also as
subsets in $\Hom(P[X],L).$

Let $A$ be an arbitrary set of points $\mu:P[X]\to L.$
It corresponds an ideal $U$ in $P[X]$ defined by
$$U=A'=\capl_{\mu\in A}\Ker \mu.$$
This is the set of all polynomials $f,$ such that every point $\mu\in A$ is
a root of $f.$
Every such ideal is called an $L$-{\it closed ideal.}

On the other hand, let $U$ be a subset in the algebra $P[X].$
One can consider $U$ as a system of equations.
It corresponds the set of points $A\subset \Hom(P[X],L),$ defined by the rule
$$A=U_L'=U'=\{\mu: P[X]\to L\bigm| U\subset\Ker\mu\}.$$
This means that a point $\mu$ belongs to $A$ if $\mu$ is a root of every
polynomial from $U.$
Such sets $A$ are called {\it algebraic sets} or {\it closed sets}.

Define two closures:
$$A''=(A')';\quad U_L''=U''=(U_L')'.$$

\definition{Definition 1}

Two extensions $L_1 $ and $L_2$ of the ground field $P$ are  called {\it
geometrically
equivalent}, if for every finite set $X=\{x_1,\dots,x_n\}$ and for every
$U\subset P[X],$ the
equality
$$U_{L_1}''=U_{L_2}''$$
takes place.
\enddefinition

\proclaim{Theorem 1}
\roster \item"1." If the field $P$ is algebraically closed, then all its
extensions are geometrically
equivalent.
\item"2." If every two extensions of $P$ are geometrically equivalent, then
$P$ is
algebraically closed.
\endroster
\endproclaim

The first part of the theorem follows from Hilbert Nullstellensatz, the
second one is deduced
easily from definitions.

We consider the category of algebraic sets over the given $L$ and fixed $P.$
Denote this category by $K_P(L).$
Its objects are algebraic sets $A$ with different $X.$

To the category $K_P(L)$ corresponds the dual category $C_P(L)$, whose
objects are
$P$-algebras of the form $P[X]/A'.$
We later define these categories in an arbitrary variety $\Th.$

If the fields $L_1$ and $L_2$ are geometrically equivalent, then the
categories $C_P(L_1)$ and
$C_P(L_2)$ coincide and the categories $K_P(L_1)$ and $K_P(L_2)$ are
isomorphic.

The category $K_P(L)$ can be considered as a natural geometric invariant of
the $P$-algebra
$L.$
It can also be viewed as a kind of characteristic which measures
$P$-algebraic closeness of
the field $L.$
In other words, it characterizes possibilities to solve in $L$ systems of
equations with
coefficients from $P.$
We say that a field $L$ is $P$-algebraically closed if, for any finite $X$
and a proper ideal $U\subset P[X],$
the set $U_L'$ is not empty.
If the field $P$ is algebraically closed, then each of its extensions $L$ is
$P$-algebraically closed.
If $L_1 $ and $L_2$ are $P$-algebraically closed, then they are
geometrically equivalent.
If $L_1$ is $P$-algebraically closed, and $L_1$ and $L_2$ are geometrically
equivalent, then
$L_2$ is $P$-algebraically closed.
The definition above can be viewed as a definition of strict
algebraic closeness of the field $L$ over $P$. The weak algebraic closeness
of the field $L$ over $P$ means that every polynomial with one variable and
coefficients from $P$ has roots in $L$. In fact, strict and weak
closeness coincide. Indeed, algebraic (absolute) closure $\bar P$
of the field $P$ coincides with the weak $P$-closure of the $P$.
Therefore, $\bar P$ is contained in every $P$-closed extension of
$L$. Then, Hilbert's Nullstellensatz works.

In the definitions above the classical Nullstellensatz can be formulated
as follows: if the field $L$ is $P$-algebraically closed then for every
ideal $U\in P[X]$ the equality
$$
U''_L=\sqrt U
$$
holds.

The pointed above equivalence of definitions plays a crucial
role in the algebraic geometry in different varieties $Var-P$.
This equivalence fails if we doing general varieties $\Theta$
and $H\in \Theta$.

Note that we distinguish algebraic sets in affine space and algebraic
varieties.
It is assumed that an algebraic variety is an algebraic set considered up
to isomorphisms in
the category $K_P(L).$

The following problem is of special interest.

Let $L_1 $ and $L_2$ be two extensions of the field $P.$
When are the categories $K_P(L_1)$ and $K_P(L_2)$  isomorphic?

Let $\mu_1:P\to L_1$ and $\mu_2:P\to L_2$ be injections defining extensions.

An isomorphism of extensions means that there is a commutative diagram

$$
\CD
P @>h_1>> L_1\\
@. @/SE/h_2// @VV\mu V\\
@. L_2\\
\endCD
$$

where $\mu$ is an isomorphism of rings.

Two extensions $h_1: P\to L_1$ and $h_2: P\to L_2$ are called {\it
semi-isomorphic} if there
is a diagram
$$\CD P @>h_1>> L_1\\ @V\sigma VV @VV\mu V\\
P@>h_1>> L_2\endCD$$
where $\mu$ is an isomorphism of rings and $\sigma$ is an automorphism of
the field $P.$

The semi-isomorphism does not imply the equivalence but if $L_1$
and $L_2$ are semi-isomorphic then the categories $K_P(L_1)$ and
$K_P(L_2)$ are isomorphic.

\definition{Definition 2}
Two extensions $L_1$ and $L_2$ are said to be geometrically
equivalent up to a semi-isomorphism, if there exists a field $L,$
such that $L$ is semi-isomorphic to $L_1$ and geometrically
equivalent to $L_2.$\enddefinition

\proclaim{Theorem 2} Two categories $K_P(L_1)$ and $K_P(L_2)$ are
isomorphic if and only if the extensions $L_1$ and $L_2$ are
geometrically equivalent up to a semi-isomorphism.\ep

The proof of this theorem is based on consideration of
automorphisms of the category of polynomial algebras $P[X]$ with
different $X$. It can be proved that every such an automorphism is
a semi-inner automorphism. Definitions of inner and semi-inner
automorphisms are natural and will be given in Lecture 3.

In order to get further information on geometric equivalence of
fields, we use the idea of a quasi-identity.

A quasi-identity $u$ in the variety $\Var$-$P$ is a formula of the
form
$$f_1\equiv 0\wedge f_2\equiv 0\wedge\dots\wedge f_n\equiv 0\Rightarrow
f\equiv 0$$
where $f_i\in P[X],$\ $i=1,\dots,n.$

\proclaim{Theorem 3} Two extensions $L_1 $ and $L_2$ of a field
$P$ are geometrically equivalent if and only if they satisfy the
same quasi-identities. \ep

Thus in this case geometric equivalence coincides with equivalence
in logik of quasi-identities in the variety $Var-P$.

This theorem implies that any extension $L$ of $P$ is geometrically
equivalent to each its
ultrapower. Every $L$ and every its ultrapower have the same
elementary theories. The following problem looks quite natural:

\proclaim{Problem 1}
When two geometrically equivalent extensions $L_1 $ and $L_2$ of a field
$P$ have different elementary theories in the logik of $Var-P$?
\ep

There are examples of such kind. Let, for example, $P$ be an algebraically
closed field and $L$ be its non-algebraically closed extension.
Then $L$ and $P$ are geometrically equivalent, while they are not
elementary equivalent. In particular, they are not elementary equivalent
in the logic of $Var-P$.

Indeed, let $f(x)=\alpha_0+\alpha_1 x\cdots +\alpha_n x^n $ be a
polynomial with the coefficients in $L$ and without roots in $L$.
Take a polynomial over $P$
$$
\varphi(x,y_0,\ldots,y_n)= y_0+y_1x+\cdots +y_nx^n
$$
and consider a $P$-formula
$$
\forall y_0\ldots y_n\exists x (\varphi(x,y_0,\ldots,y_n)=0).
$$
This formula holds in $P$ and does not hold in $L$.

In the Problem 1 one has to consider the general situation.

\subheading{3. Some information from universal algebra}

Almost all that which has been determined for the classical case admits
natural generalizations.

Let us recall some information from universal algebra.
Fix a signature $\Om,$ i.e., a set of symbols of operations of the
arbitrary arity.
Consider the class of all $\Om$-algebras.
Every set of identities, satisfied by the operations from $\Om,$ specify
some subclass in the
class of all $\Om$-algebras.
Every such subclass $\Th$ is called {\it a variety of algebras.}
Groups, rings, semigroups, associative and Lie rings, associative and Lie
algebras over
fields, $\Var$-$P$ are all examples of varieties.

Speaking of algebra, we think of universal algebra, i.e., algebra in an
arbitrary but fixed variety $\Theta$. One can assume that this is a group,
or an associative or Lie algebra. Algebra can be viewed as an universal or
a concrete one.

In each variety $\Th$ there are free algebras $W(X),$ where $X$ is a set.
This is a significant feature of varieties which plays the crucial role
in the theory under consideration.

The characteristic property of $W(X)$ can be described by the commutativity
of the following
diagram

$$
\CD
X @>id >> W(X)\\
@. @/SE/\mu// @VV{\bar \mu} V\\
@. H\\
\endCD
$$

Here, $H$ is an arbitrary algebra in $\Th$,\ id is the identity map, $\mu$
is an arbitrary
map, and $\overline \mu$ is the homomorphism in $\Th$ uniquely determined
by $\mu.$

In other words, this diagram means that the set $X$ freely generates $W(X)$
(or, $X$ is a
basis for $W$), the variables from $X$ take values from $H$ by the rule
$\mu,$ and the
homomorphism $\overline \mu$ computes in $H$ the values of arbitrary
elements $w\in W.$

Elements from $W$ are used to call words, terms, polynomials, or
$\Theta$-polynomials.
They are constructed from the variables from $X$ by the rules of the
variety $\Theta.$
Every $w\in W$ is uniquely presented as $w=w(x_1,\dots,x_n),$\ $x_i\in X,$
and uniqueness
means that $w(x_1,\dots,x_n)=w'(x_1,\dots,x_n)$ holds if and only if
it is an identity in $\Theta.$

Every variety $\Theta$ can be regarded as category whose morphisms are
homomorphisms in
$\Theta.$
To each homomorphism $\mu: G\to H$ corresponds its kernel $\Ker \mu.$
This is a special binary relation, defined by the rule:
$g_1(\Ker\mu)g_2 $
if and only if
$g_1^\mu=g_2^\mu.$
A relation $T$ in $G$ is called a congruence if $g_1Tg_1'\wedge\dots\wedge
g_n Tg_n' $ implies
$g_1\dots g_n\omega Tg_1'\dots g_n'\omega$ for every $n$-ary operation
$\om$ in $\Om.$ The kernel $Ker\mu$ is always a congruence.

For every congruence $T$ in $G$ one can take a factor algebra $G/T$ with
the natural
homomorphism $G\to G/T.$
The diagram

$$
\CD
G @>\mu >> H\\
@. @/SE/\mu_0// @AA\mu_1 A\\
@. G/Ker \mu\\
\endCD
$$

gives a canonical decomposition of a  homomorphism $\mu$.
The set of all homomorphisms $G\to H$ is denoted by $\Hom(G,H).$

Every variety of algebras is closed under taking Cartesian products,
subalgebras and
homomorphic images.
In particular, if $I$ is a set and $H\in\Th$, then the Cartesian power
$H^I$ also belongs to
$\Theta.$
If $|I|=n$, then $H^I$ is denoted by $H^{(n)}.$
The algebra $H^{(n)}$ is viewed as an affine space of points
$(a_1,\dots,a_n),$ \ $a_i\in H.$
Let $X=\{x_1,\dots,x_n\}.$
Then there is a bijection
$$\alpha_X:\Hom(W(X),H)\to H^{(n)},$$
where for every $\nu: W\to H,$
$$\a_X(\nu)=(\nu(x_1),\dots,\nu(x_n)).$$
This bijection induces a structure of algebras in $\Theta$ on $\Hom(W(X),H).$
By definition, if $\om\in\Om$ is an $m$-ary operation,
$\nu_1,\dots,\nu_m\in\Hom(W(X),H),$
then
$$(\nu_1\dots\nu_m\om)(x)=\nu_1(x)\dots \nu_m(x)\om,$$
for every $x\in X.$

We define the notion of a commutative algebra $H.$
For a commutative algebra $H$ the equality
$$(\nu_1,\dots\nu_m\om)(w)=\nu_1(w)\dots\nu_n(w)\om$$
always holds for every $w\in W.$

Let $w_1,w_2\in\Om$ be two operations of arity $m_1,\ m_2$, respectively.
Consider the matrix
$$\pmatrix &x_{11}&&\dots && x_{1m_1}\\
&\dots&&\dots&&\dots\\
&x_{m_21}&&\dots &&x_{m_2m_1}\endpmatrix$$
Two operations $w_1$ and $w_2$ are said to be commutative if
$$ ((x_{11}\ldots x_{1m_1}\om_1)\ldots(x_{m_21}\ldots x_{m_2m_1})\om_1)\om_2=$$
$$
=((x_{11}\ldots x_{m_21}\om_2)\ldots(x_{1m_1}\ldots x_{m_2m_1}\om_2)\om_1.$$
An algebra $H$ is called commutative if every two operations, not
necessarily different,
commute.
For nullary operations this means that all of them coincide.

If $H$ is a commutative algebra, then the algebra $\Hom(W(X),H)$ is also
commutative and it
is a subalgebra in the algebra of mappings $G^{W(X)}.$

In the general case the set $\Hom(W(X),H),$\ $|X|<\iy$ is considered as an affine space over an
algebra $H.$

Let now $\Theta^0$ denote the category of all free algebras $W(X)$ in $\Theta$
with finite $X.$
The category $\Theta^0$ is a full subcategory in the category $\Theta.$
Morphisms in $\Theta^0$ are presented by arbitrary homomorphisms $s:
W(X)\to W(Y).$

Let us make remarks on equations.
Let us fix a finite set $X=\{x_1,\dots,x_n\}$ and consider equations of the
form $w=w',$\
$w,w'\in W(X).$
Every such equation with the given $X$ is considered also as a
formula in the
logic in $\Theta.$ In the later case we write $w\equiv w'$.

A point $\nu: W(X)\to H$ is a root of the equation
$w(x_1,\dots,x_n)=w'(x_1,\dots,x_n),$ if
$w(x_1^\nu,\dots ,x_n^\nu)=w'(x_1^\nu,\dots, x_n^\nu).$
This also means that the pair $(w,w')$ belong to $\Ker \nu.$
We will identify the pair $(w,w')$ and the equation $w=w'.$

In order to get a reasonable geometry in $\Theta$ one has to consider the
equations with constants. The next subsection deals with constants.
\bk
\subheading{4. Algebras with  a fixed algebra of constants}

Let $\Theta$ be an arbitrary variety of algebras, $G$ be a fixed algebra
in $\Th,$\
$|G|>1.$
Consider the new variety, denoted by $\Th^G.$
First we define the category $\Th^G.$
The objects in $\Th^G$ have the form $h: G\to H,$ where $h$ is a
homomorphism in $\Th,$ not
necessarily injective.
Morphisms in $\Th^G$ are presented by the commutative diagrams

$$
\CD
G @>h >> H\\
@. @/SE/h'// @VV\mu V\\
@. H'\\
\endCD
$$

where $\mu,h,h'$ are homomorphisms in $\Th.$

An algebra $H,$ considered as a $G$-algebra, is denoted by $(H,h).$
In this way, to elements $g\in G$ correspond constants, i.e., nullary
operations in $H.$
We add them to the signature $\Om.$
This signature allows us to consider a variety of $G$-algebras.
Identities of a variety of $G$-algebras are presented by identities of
$\Th$ and by the
defining relations of the algebra $G$ (see \cite{Pl1} for details).

A free in  $\Th^G$ algebra $W=W(X)$ has
the form  $G\ast W_0(X),$ where $W_0(X)$ is the free in the $\Th$
algebra over $X,$ $\ast$ is the free product in $\Theta$ and
the embedding $i_G: G\to W(X)=G\ast W_0(X)$ follows from the definition of
free product.

A $G$-algebra $(H,h)$ is called a {\it faithful G}-algebra if $h:
G\to H$ is an injection. A free algebra $(W,i_G)$ is  faithful. A
$G$-algebra $G$ with the identical $G\to G$ is also faithful. All
other $G$-algebras $G$ are isomorphic to this one. All of them are
simple, i.e., they do not have faithful subalgebras and
congruences. Let $(H,h)$ be a $G$-algebra, and $\mu:H\to H'$ is a
homomorphism in $\Th.$ Then, by $h'=\mu h$, \ $H'$ becomes a
$G$-algebra, and $\mu$ is a homomorphism of $G$-algebras. Since
one can start from an arbitrary congruence $T$ in $H$ and from the
natural congruence $H/T,$ we say that $T$ is called faithful if
the $G$ algebra $H/T$ is faithful. A congruence $T$ is faithful if
and only if $g_1^h=g_2^h$ is equivalent to $g_1=g_2.$

Let a homomorphism

$$
\CD
G @>h >> H\\
@. @/SE/h'// @VV\mu V\\
@. H'\\
\endCD
$$

be given, and let $(H',h')$ be a faithful $G$-algebra.
Then $(H,h)$ is a faithful $G$-algebra.
If $T=\Ker \mu,$ then $T$ is a faithful congruence and $H/T$ is also faithful.

Note that the inclusion $G \in \Theta$ means that all the constants of
$G$-algebra $H$ are covered by the elements from $G$.

\example{Examples}
A variety $\Var$-$P$ is a variety of the type $\Th^G$, where $\Th$ is a
variety of
associative and commutative rings with $1$, and $G$ is a field $P.$
In this example, elements of the field $P$ are the constants in the
$P$-algebras.
They considered as nullary operations, and, simultaneously, using
multiplication, we
can look at them as unary operations.

$G$-groups is another example of $G$-algebras.
Here, elements from $G$ also can be viewed as unary operations.

We denote also a free $G$-algebra $W=W(X)$ by $G[X].$
Since $\Th^G$ is a variety of algebras, all constructions like Cartesian
and free products,
subalgebras and homomorphisms are naturally defined for $\Th^G.$

Consider, further, the special condition on $\Th^G$ denoted by $(\ast)$.
Namely, we assume that the algebra $G$ generates the whole variety $\Th^G,$
i.e., in $G$
there are no non-trivial identities with the coefficients from $G.$
This condition is fulfilled in $\Var$-$P$ if the field $P$ is infinite and
in $\Grp-F=(\Grp)^F.$

Every faithful $G$-algebra $H$ contains $G$ as a subalgebra.
Thus, $(\ast)$ implies that every faithful $G$-algebra $H$ generates the
whole variety
$\Th^G,$ i.e., in $\Th^G$ there are no proper subvarieties containing
faithful algebras.

In the category $\Th^G$ along with morphisms, one can consider also
semimorphisms.
They have the form
$$\CD G @>h >> H\\ @V\sigma VV @VV\mu V\\
G@>h'>> H'\endCD$$
where $\sigma\in\End G.$
Then, we can consider semi-isomorphic $G$-algebras.

Another possibility is to vary also the algebra of constants $G.$
This leads to the diagram of the form
$$\CD G @>h >> H\\ @V\sigma VV @VV\mu V\\
G'@>h'>> H'\endCD$$
with  component-wise multiplication.

Let us make a remark on equations.
The equations of the form $w=w'$ with $w,w'\in W(X)=G\ast W_0(X)$ are
equations with
constants.
Consider systems of such equations $T.$
If $T$ is a congruence, then $T$ has a solution in a faithful $G$-algebra
$H$ if and only if
$T$ is a faithful congruence in $W=W(X).$
Thus, a system $T$ has a solution, if $T$ is contained in a faithful
congruence in $W.$
Note that, by definition, all faithful congruences are proper.
\endexample
\bk

\subheading{5. Special categories}

Let $\Th$ be a variety of algebras and $H\in\Th.$
Consider the category of affine spaces $K_\Th^0(H).$
Its objects are represented by affine spaces $\Hom(W(X),H).$
Morphisms have the form
$$\tilde s: \Hom(W(X),H)\to\Hom(W(Y),H),$$
where $s:W(Y)\to W(X)$ is a morphism in $\Th^0.$
Here, $\tilde s$ is defined by $\tilde s(\nu)=\nu s$ for every point $\nu:
W(X)\to H.$

If an algebra $H$ generates the variety $\Th$ then $K_\Th^0(H)$ is dual to
the category
$\Th^0.$

Another category is denoted by $\Pol_\Th(H).$
Its objects have the form of affine spaces $H^{(n)},$ while morphisms are
the special
polynomial maps
$$s^\alpha:H^{(n)}\to H^{(m)}.$$
If, further, $s:W(Y)\to W(X)$ is a morphism in $\Th$ and $|Y|=m,$\ $|X|=n,$
then there is a
commutative diagram

$$\CD \Hom(W(X),H) @>\tilde s>>\Hom(W(Y),H)\\ @V\alpha_X VV @VV\alpha_Y V\\
H^{(n)}@>s^\alpha>> H^{(m)}\endCD$$

This yields a canonical duality between $K_\Th^0(H)$ and $\Pol_\Th$-$H$.
In this diagram, all maps are algebra homomorphisms, if in $\Th$ all
operations commute.

Denote $\Bool_\Th(H)$ to be a special category of Boolean algebras.
Its objects $Bool(W(X),H)$   are the Boolean algebras
of subsets in $\Hom(W(X),H).$

Morphisms
$$s_\ast:\Bool(W(X),H)\to\Bool(W(Y),H)$$
are defined for every $s:W(X)\to W(Y).$
Here, $\tilde s: \Hom(W(Y),H)\to\Hom(W(X),H)$ is given, and for every $A\subset
\Hom(W(X),H)$, \
$B=s_\ast A$ is the inverse image of the set $A$ under the map $\tilde s.$

\newpage

\baselineskip 12pt
\topmatter \rightheadtext{Lecture 2}
\leftheadtext{B. Plotkin}
\title Lecture 2\\
\quad\\
 Algebraic sets and algebraic varieties\endtitle
\endtopmatter

\centerline{\smc Contents}
\bigskip
\roster\item "1." {\smc Systems of equations and algebraic sets}
\newline
\item"2." {\smc Lattices of algebraic sets}
\newline
\item"3." {\smc Categories of algebraic sets}
\newline
\item"4." {\smc On the notion of algebraically closed algebras}
\newline
\item"5." {\smc Changing the basic variety of algebras}
\newline
\item"6." {\smc Zariski topology}
\endroster

\newpage

\subheading{1. Systems of equations and algebraic sets}

Let $\Th$ be a variety of algebras, $W=W(X)$ be the free algebra in $\Th$
with $|X|<\iy.$

Fix an algebra $H$ from $\Th.$
Consider equations of the form $w\equiv w',$\ $w,w'\in W.$
Denote the value of this equation in $H$ by $\Val_H(w\equiv w').$
It is defined by
$$\Val_H(w\equiv w')=\{\mu:W\to H\bigm |w^\mu=w'{^\mu}\}.$$
Systems of equations in $W$ are denoted by $T.$
They can be viewed as binary relations in $W.$
Some of these relations can be congruences in the algebra $W$.

Consider, also, sets of points $A\subset \Hom(W,H).$
The Galois correspondence between $T$ and $A$ is as follows:
$$\cases T'=A=\{\mu:W\to H\bigm| T\subset \Ker \mu\}=T'_H\\
A'=T=\capl\limits_{\mu\in A}\Ker\mu\endcases$$
The same can be written in the form,
$$T'=\capl_{(w,w')\in T}\Val_H\ (w\equiv w')$$
$$A'=\{w\equiv w'\bigm| A\subset \Val_H\ (w\equiv w')\}$$

\definition{Definition 1}
Every set $A$ such that $A=T'$ for some $T$ is called a closed or an $H$-closed
set.\enddefinition
These are algebraic sets over the algebra $H.$

\definition{Definition 2}
Algebraic sets considered up to isomorphisms in the category of algebraic
sets $K_\Th(H)$ are
called {\it algebraic varieties}.\enddefinition

The definition of the category $K_\Th(H)$ is given in the item 3
of this lecture.

If $A$ is an arbitrary set of points, then its closure is defined to be
$A''=(A')'.$

\definition{Definition 3} Congruences $T$ in $W$ of the form $T=A'$ are
called $H$-closed
congruences.\enddefinition

For an arbitrary system of equations $T$ its closure is defined by
$T_H''=(T_H')'.$

\proclaim{Proposition 1}
 A congruence $T$ in $W$ is $H$-closed if and only if for some set $I$ there
is an injection
$$W/T\to H^I.$$\endproclaim

Let, further, $T$ be a congruence in $W,$ and $\mu_0: W\to W/T$ be a
natural homomorphism.
Consider the sets $\Hom(W/T,H)$ and $\Hom(W/T,H)\mu_0,$ where
$$\Hom(W/T,H)\mu_0=\{\nu\mu_0: W\to H\},\ \nu\in\Hom(W/T,H).$$
Denote $\tilde\mu_0(\nu)=\nu\mu_0.$

\proclaim{Proposition 2}
The following formula
$$T_H'=A=\Hom(W/T,H)\mu_0$$
holds.
Moreover, $\tilde\mu_0:\Hom(W/T,H)\to A$ is a bijection.\endproclaim

Let us consider, separately, the case $H$ is a commutative algebra.
Then the sets $\Hom(W,H)$ and $\Hom(W/T,H)$ are also commutative algebras.
Every algebraic set $A$ is a subalgebra in $\Hom(W,H).$
The bijection $\tilde\mu_0:\Hom(W/T,H)\to A$ is an isomorphism of algebras.

Now, let us make some remarks on trivial cases.

A congruence $T$ is called a zero congruence if $wTw'$ means that
$w$ and $w'$
coincide in $W.$
Here $T$ is a relation of equality.
$T$ is called a unity or a non-proper congruence, if $wTw'$ is fulfilled
for every $w,w'\in
W.$
These congruences are denoted $T=0,$\ $T=1,$ respectively.

We have $0'=\Hom(W,H)$ and, therefore, $\Hom(W,H)$ is an algebraic set.
As for $1',$ this is either an empty set, or a zero point in $\Hom(W,H),$
sending $W$ to the
zero subalgebra in $H,$ if the latter exists.
A zero subalgebra consists of one element and is defined by a unique
nullary operation.

It is easy to see that if $A=\Hom(W,H)$, then $T=A'$ is the congruence of
all identities of
an algebra $H$ in the free algebra $W.$
This is the minimal closed congruence in $W.$
If an algebra $H$ generates the whole variety $\Th,$ then $A'=0.$
If $A$ is the empty set, then $A'=1.$
If $A$ consists of one zero point, then $A'=1$ as well.
If $T$ is the empty set, then assume that $T'=\Hom(W,H).$

For arbitrary algebras $H$ and $G$, denote
$$(H-\Ker)(G)=\capl_{\nu: G\to H}\Ker \nu.$$

Let $T$ be a congruence in $W$ and take $\mu_0: W\to W/T.$
Then
\proclaim{Proposition 3} \rom{(Hilbert's Nullstellensatz)}
$$T_H''=\mu_0^{-1}((H-\Ker\mu)(W/T)).$$
\endproclaim

This proposition can be considered as one of the forms of Hilbert's theorem.
There are, also, others.

Let $A$ be an algebraic set.
To each element $w\in W$ corresponds a map
$$\ov w: A\to H$$
defined by $\ov w(\nu)=w^\nu$ for every $\nu\in A.$
Such $\ov w$ are called {\it regular maps}.
One can define the algebra of regular maps.
This algebra is isomorphic to an algebra $W/A'$, which is called the {\it
coordinate algebra}
of the algebraic set $A.$

From the Proposition 1 follows that an algebra $G\in\Theta$ can
be presented as a coordinate algebra of an algebraic set $A$
over given algebra $H\in\Theta$ if and only if $G$ is finitely
generated algebra and there is an injection $G\to H^I$ for some set $I$.

\bk

\subheading{2. Lattices of algebraic sets}

Given a set $I$, the following formulas hold:

1. \ $(\cup A_\a)'=\cap A_\a'$

2.\ $(\cup T_\a)'=\cap T_\a'$

3.\ $\cup T_\a'\subset (\cap T_\a)'$

4.\ $\cup A_\a'\subset (\cap A_\a)'$

\noindent where $\a\in I.$

If all $A_\a$ and $T_\a$ are $H$-closed sets, then

5.\ $(\cup T_\a')''=(\cap T_\a)'$

6.\ $(\cup A_\a')''=(\cap A_\a)'.$

Thus, one can state that the intersection of algebraic sets is an algebraic
set, and the
intersection of closed congruences is again a closed congruence.

It is evident that the union of two closed congruences
is not a closed congruence,
and the union of
two algebraic sets is not necessarily an algebraic set (say, if $H$ is
a commutative group it is not
an algebraic set).

For given $H\in\Th$ and $W=W(X),$ \ $|X|<\iy,$ denote by

$\Alv_H(W)$\ the set of all algebraic sets in $\Hom(W,H),$

$\Cl_H(W)$ \ the set of all $H$-closed congruences in $W.$

This give rise to the functors
$$\Alv_H:\Th^0\to \text{Set},$$
$$\Cl_H:\Th^0\to \text{Set}.$$
The first functor is covariant, while the second one is contravariant.
If $s: W(Y)\to W(X)$ is a morphism in $\Th^0,$ then for every $B\in
\Alv_H(W(Y))$ the
corresponding $s_*B=A$ is contained in $\Alv_H(W(X)).$

There is a map
$$s_*:\Alv_H(W(Y))\to \Alv_H(W(X)).$$
If $T\in Cl_H(W(X))$ and $s_*T$ is the inverse image of $T$ in $W(Y),$ then
$s_*T\in\Cl_H(W(Y)).$
This gives the map
$$s_*:\Cl_H(W(X))\to\Cl_H(W(Y)).$$
The sets $\Alv_H(W)$ and $\Cl_H(W)$ constitute the lattices.
If $A$ and $B$ are the elements from $\Alv_H(W),$ then we set
$$A\ov\cup B=(A\cup B)''.$$
If $T_1$ and $T_2$ belong to $\Cl_H(W),$ then
$$T_1\ov\cup T_2=(T_1\cup T_2)''.$$

\proclaim{Proposition 4}
Lattices $\Alv_H(W)$ and $\Cl_H(W)$ are dual.\endproclaim

The duality is determined by $A\to A'$ and the properties 5,6 above.

In the general case, we cannot state that the functors $\Alv_H$ and $\Cl_H$
are coordinated
with lattice operations.

\definition{Definition 4}
An algebra $H$ is called geometrically stable, if for every $W(X)=W$ and any
$A,B\in\Alv_H(W),$ the equality
$$A\ov\cup B=A\cup B$$
takes place.\enddefinition

This is equivalent to the fact that for  closed systems of equations $T_1$
and $T_2$ in $W,$
the equality
$$T_1'\cup T_2'=(T_1\cap T_2)'$$
holds.

If $H$ is a stable algebra, then every lattice $\Alv_H(W)$ is a
distributive lattice in
$\Bool_\Th(W,H).$
A dual lattice $\Cl_H(W)$ is also distributive.

In this case, the functors $\Cl_H$ and $\Alv_H$ are the functors to the
category of
distributive lattices.

\proclaim{Proposition 5}
If $s: W(Y)\to W(X)$ is surjective, then the maps
$$s_*:\Cl_H(W(X))\to \Cl_H(W(Y))\quad\text{and}$$
$$s^*:\Alv_H(W(Y))\to\Alv_H(W(X))$$
are homomorphisms of lattices.\ep

Here the map $s^*$ is conjugated to the map $s_*.$

For every algebraic set $A\in \Alv_H(W)$, denote by $L(A)$ the lattice of
all algebraic sets
in $A.$
This lattice is considered as a lattice invariant of $A.$

Another invariant is based on the semigroup $\End(H).$
For every $\dl: H\to H$ we have $\tilde\dl:\Hom(W,H)\to \Hom(W,H)$ by the rule
$\tilde\dl(\nu)=\dl\nu.$
Thus, the semigroup $\End(H)$ acts in the affine space $\Hom(W,H).$
It is clear that every algebraic set $A$ in this space is invariant with
respect to the
action of $\End(H).$
Thus, we can look for the structure of $A$ with respect to this action.

\bk

\subheading{3. Categories of algebraic sets}

Denote the category of algebraic sets in $\Th$ over the given $H$ by
$K_\Th(H).$
Its objects are denoted by
$$(X,A)$$
where $A$ is an algebraic sets in the affine space $\Hom(W(X),H).$

In order to define morphisms
$$(X,A)\to(Y,B)$$
we start from $s: W(Y)\to W(X).$
We say that $s$ is admissible with respect to $A$ and $B,$ if $\tilde
s(\nu)=\nu s\in B$ for
every $\nu\in A.$
The equivalent condition is $A\subset s_*B.$
If $s$ is admissible with respect to $A$ and $B$, it corresponds the map
$[s]: A\to B.$
We write
$[s]:(X,A)\to (Y,B).$
The multiplication of such morphisms is defined in a natural way and the
category $K_\Th(H)$
is defined.

Let us define now the category $C_\Th(H).$
Its objects are algebras from $\Th$ of the form $W(X)/T,$ where $|X|<\iy$,
and $T$ is an $H$-closed congruence.
Morphisms are presented by homomorphisms
$$\sigma: W(Y)/T_2\to W(X)/T_1.$$
Here, $C_\Th(H)$ is a full subcategory in the category $\Th.$

Every such $\sigma$ can be represented as $\sigma=\ov s,$ where $s: W(Y)\to
W(X)$ is an
admissible homomorphism of free algebras with respect to congruences $T_2 $
and $T_1.$
This means that $wT_2w'$ implies $w^sT_1w'{}^s.$

\proclaim{Proposition 6}
A homomorphism $s:W(Y)\to W(X)$ is admissible with respect to $A$ and $B$
if and only if $s$
is admissible with respect to $B'=T_2$ and $A'=T_1.$\ep

Besides, if $s_1:W(Y)\to W(X)$ is admissible with respect to $A$ and $B,$
and, hence, with
respect to $T_2$ and $T_1,$ then the equality $[s_1]=[s]$ holds if and only
if $\ov s_1=\ov
s.$

This leads to the following theorem
\proclaim{Theorem 1}
The transition
$$(X,A)\to W(X)/A'$$
defines the duality of the categories $K_\Th(H)$ and $C_\Th(H).$\ep

This fact confirms the correctness of the definitions of $K_\Th(H)$ and
$C_\Th(H).$

Define now the categories $K_\Th$ and $C_\Th.$
In these categories an algebra $H$ is not fixed.
The objects of $K_\Th$ have the form
$$(X,A,H),$$
where $A$ is an algebraic set in $\Hom(W(X),H).$
The objects in $C_\Th$ have the form
$$(W(X)/T,H),$$
where $T$ is an $H$-closed congruence in $W(X).$

Consider morphisms
$$(X,A,H_1)\to(Y,B,H_2).$$
Take $s:W(Y)\to W(X)$ and a homomorphism $\dl: H_1\to H_2.$

Consider a commutative diagram
$$\CD W(Y) @>s>> W(X)\\ @V\nu' VV @VV\nu V\\
H_2@<\dl<< H_1\endCD$$
For every point $\nu: W(X)\to H_1$ denote $(s,\dl)(\nu)=\nu'=\dl\nu s.$

A pair $(s,\dl)$ is said to be admissible with respect to $A$ and $B,$
if $\nu'=(s,\dl)(\nu)\in B$ for every $\nu\in A.$
Let $(s,\dl)$ be an admissible pair with respect to $A$ and $B.$
Fix $\dl$ and consider the map $[s]_\dl: A\to B.$
The pair $([s]_\dl,\dl)$ is considered to be a morphism
$$(X,A,H_1)\to(Y,B,H_2).$$
Let now
$$([s]_\dl,\dl):(X,A,H_1)\to(Y,B,H_2)\quad\text{and}$$
$$([s']_{\dl'},\dl'):(Y,B,H_2)\to (Z,C,H_3)$$
be two morphisms.
Then, set
$$([s']_{\dl'},\dl')([s]_\dl,\dl)=([ss']_{\dl'\dl},\dl'\dl):(X,A,H_1)\to
(Z,C,H_3).$$

These definitions are correct and the  category $K_\Th$ is defined.

The categories $K_\Th$ and $C_\Th$ are dual.

The categories $K_\Th(H)$ and $C_\Th(H)$ are subcategories in $K_\Th$ and
$C_\Th$, and
duality between $K_\Th$ and $C_\Th$ induces duality between $K_\Th(H)$ and
$C_\Th(H).$

Let $T$ be a system of equations in $W(X).$

Define a full subcategory $K_\Th(T)$ in $K_\Th.$

The objects have the form
$$(X,A,H)$$
where $X$ is fixed, $A=T_H',$ and $H,A$ are changed.
The category $K_\Th(T)$ characterizes possibilities to solve the system $T$ in
different algebras
$H\in\Th.$

Consider separately the case when there is an injection $\dl: G\to H$.
It can be proven that this injection induces in $K_\Th$ an isomorphism
between the category
$K_\Th(G)$ and a subcategory in $K_\Th(H).$

\bk

\subheading{4. On the notion of algebraically closed algebras}

As previously in the classical situation, we consider $K_\Th(H)$ as a
geometrical invariant
of the algebra $H.$
This category is viewed also as a measure of algebraic closeness of an
algebra $H.$
An algebra $H$ can be not an algebraically closed algebra, but
$H$ is algebraically
closed with some
measure.

Consider the notion of algebraically closed algebras.

\definition{Definition 5}
An algebra $H\in\Th$ is called {\it algebraically closed} in $\Th,$ if for
every finite $X$
and every non-unity congruence $T$ in $W(X)$ this $T$ has a non-zero root
$\mu: W(X)\to H.$
\enddefinition

Here, $\Ker \mu$ is a non-unity congruence and $T\subset \Ker \mu.$
In particular, if $T$ is a maximal congruence in $W=W(X)$, then $\Ker
\mu=T$ and $T$ is an
$H$-closed congruence.
Every maximal congruence $T$ is $H$-closed.

It is clear that if an algebra $H$ is contained in some other $H_1\in\Th,$
and $H$ is
algebraically closed, then $H_1$ is also algebraically closed.

Consider this definition with respect to varieties of the kind $\Th^G.$
An algebra $H\in\Th$ is considered as a $G$-algebra for different $G$ with
an injection $h:
G\to H$.
For such $(H,h)$ and different $G,$ the definition of algebraic closeness
with respect to
different varieties $\Th^G$ differ.

In $\Var$-$P$, an algebra $L$ over $P$ can be $P$-closed but not
absolutely algebraically
closed.
Here and in the general situation, absolutely algebraically closed means
algebraic closeness
in the variety $\Th^L.$

In general case all equations and systems of equations
have coefficients in $H.$
So, fix $G\in\Th$ and consider a homomorphism in $\Th^G$

$$
\CD
G @>i_G>> G\ast W_0=W\\
@. @/SE/\varepsilon// @VV\mu V\\
@. G\\
\endCD
$$

For every such $\mu$ the congruence $T=\Ker \mu$ is faithful and
maximal. The condition of algebraic closeness of $G$ in $\Th^G$
means that if $T$ is a maximal faithful congruence in $W$, then
$T=\Ker \mu$ for some $\mu.$

This is well-coordinated with the situation of the field $P$ in
$\Var$-$P$. Note, that there are other approaches to the notion of
algebraic closeness.

Consider systems of equalities $T$ of the form $w\equiv w'$ and inequalities
$w\not\equiv w'.$
Such a system has a solution if there exists a proper congruence $T_1$ such
that all $w\equiv
w'$ from $T$ belong to $T_1$ and all $w\not\equiv w'$ from $T$ do not
belong to $T_1.$
Now, one can say that an algebra $H$ is an algebraically closed, if for any
compatible system  $T$ of
equalities and inequalities  for every finite $X$ there is a point $\mu:W(X)\to H,$ which
is a root of such $T.$

Such a point of view on algebraic closeness of an algebra $H$ is closer to
algebraic geometry
in logic (see lecture 6) than to equational algebraic geometry.

\bk

\subheading{5. Changing of the basic variety of algebras}

Let $\Th$ be a variety of algebras, $\Th_0\subset \Th$ be a subvariety, and
$X$ be a finite
set.
Then $W=W(X)$ and $W_0=W_0(X)$ are the free algebras in $\Th$ and $\Th_0,$
respectively.
The homomorphism $\mu_0:W(X)\to W_0(X)$ with the kernel $T=\Ker \mu$
consists of identities
in $W$ of the variety $\Th_0.$
For every such $H$ there is a decomposition

$$
\CD
W @>\mu_0>> W_0\\
@. @/SE/\mu// @VV\nu V\\
@. H\\
\endCD
$$

To every point $\mu: W\to H$ uniquely corresponds $\nu:
W_0\to H$ with
$\mu=\nu\mu_0=\tilde \mu_0(\nu).$
We have the bijection $\mu_0:\Hom(W_0,H)\to\Hom(W,H).$

\proclaim{Proposition 7}
The bijection $\tilde\mu_0$ induces the isomorphism of lattices of the
algebraic sets
$\Alv_H(W)$ and $\Alv_H(W_0);$ and the lattices
$\Cl_H(W)$ and $\Cl_H(W_0).$\ep

The same bijections $\mu_0$ for different finite $X_0$ is used to
prove that

\proclaim{Theorem 2}
The categories $K_\Th(H)$ and $K_{\Th_0}(H)$ are canonically isomorphic.\ep

For $\Th_0\subset \Th$ one can consider the categories $K_{\Th_0}$ and $K_\Th.$
It can be proved that the first one is isomorphic to a subcategory of the
second.

Let, again, $H\in\Th_0\subset\Th.$
Then,

\proclaim{Proposition 8}
If $H$ is algebraically closed in $\Th,$ then $H$ is algebraically closed
in $\Th_0.$\ep

The opposite statement is not true in general.

One can imagine the situation when a system of equations $T$ in the given
$W=W(X)$ in $\Th$
contradicts to identities of the algebra $H.$
Then $T$ does not have a common root in $H.$

One can consider algebraic closeness of the given $H$ with respect to a
different finite $X.$
Then the question is whether the absolute algebraic closeness follows from
the algebraic
closeness with respect to one-element set $X.$
The answer depends on $\Th$ and on the structure of algebraic closure of an
arbitrary $H$ in
$\Th.$

\bk

\subheading{6. Zariski topology}

This is the minimal topology in the space $\Hom(W,H)$ in which all
algebraic sets are closed.

If the algebra $H$ is stable, then algebraic sets are, precisely, all
closed sets.
In the general case, the closed sets in Zariski topology are represented by
finite unions of
algebraic sets and their arbitrary intersections.

There is another approach which we now discuss.

Consider formulas of a more general form called {\it pseudoequalities}.
They have the form
$$w_1=w_1'\vee\dots\vee w_n\equiv w_n',\quad w_i,w_i'\in W=W(X).$$

One can consider also pseudoequations.
If $u$ is the pseudoequality above,
then its value $\Val_H(u)$ in the algebra $H$
is defined by
$$\Val_H(u)=\Val_H(w_1\equiv w_1')\vee\dots\vee \Val_H(w_n\equiv w_n').$$
Similar to the usual equalities, for every system $T$ of pseudoequalities,
one can define
$$T_H'=A=\capl_{u\in T}\Val_H(u).$$

The sets $A$ of such form are called {\it pseudoalgebraic sets}.

Every algebraic set is pseudoalgebraic.
Conversely, if an algebra $H$ is stable, then every pseudoalgebraic set is
algebraic.

\proclaim{Proposition 9}
Closed sets in Zariski topology coincide with pseudoalgebraic sets.\ep

One can consider the category of affine spaces $K_\Th^0(H)$ and every
category $K_\Th(H)$ as
  categories with topological objects.
The morphisms defined above turn out to be continuous maps.
\newpage

\baselineskip 12pt
 \topmatter \rightheadtext{Lecture 3}
\leftheadtext{B. Plotkin}
\title Lecture 3\\
\quad\\
Geometric properties and geometric relations of algebras\endtitle
\endtopmatter

\bigskip
\centerline{\smc Contents}
\bigskip
\roster\item "1." {\smc Geometric equivalence of algebras}
\newline
\item"2."{\smc  Geometric equivalence and quasi-identities}
\newline
\item"3." {\smc Geometric similarity of algebras}
\newline
\item"4." {\smc  Similarity and equivalence}
\newline
\item"5." {\smc Similarity and semi-isomorphisms}
\newline
\item"6." {\smc Geometrically Noetherian algebras}
\newline
\item"7." {\smc Geometric stability}
\endroster

\newpage

\subheading{1. Geometric equivalence of algebras}

Given a variety $\Th.$

\definition{Definition 1} Algebras $H_1$ and $H_2$ in $\Th$ are geometrically
equivalent if
for every finite $X$ and every system of equations $T$ in the free algebra
$W=W(X),$ the
equality
$$T_{H_1}''=T_{H_2}''$$
takes place.
\enddefinition

This is equivalent to the fact that always $\Cl_{H_1}(W)=\Cl_{H_2}(W)$,
i.e., the functors
$\Cl_{H_1}$, and $Cl_{H_2}$ coincide.

Also the lattices    $\Cl_{H_1}(W)$ and $\Cl_{H_2}(W)$ coincide and the
lattices
$\Alv_{H_1}(W)$ and $\Alv_{H_2}(W)$ are isomorphic.

\proclaim{Proposition 1}
Every algebra $H$ is geometrically equivalent to any of its cartesian
powers $H^I.$\ep

This proposition follows from the fact that the algebras $H_1$ and $H_2$
are geometrically
equivalent if and only if for any $W$ and a congruence $T$ in $W$ the equality
$$(H_1-\Ker)(W/T)=(H_2-\Ker)(W/T)$$
takes place.

\proclaim{Proposition 2}
Let $\Th$ be a variety and $\Th_0$ a subvariety in $\Th.$
The algebras $H_1$ and $H_2$ from $\Th_0$ are geometrically equivalent in
$\Th_0$ if and only
if they are equivalent in $\Th.$\ep

Consider one more point of view on the notion of geometric equivalence.

Let $H_1$ and $H_2$ from $\Th$ be given.
There exists a canonical functor
$$F=F(H_1,H_2): K_\Th(H_1)\to K_\Th(H_2).$$
If $(X,A)$ is an object in $K_\Th(H_1)$, then
$$F((X,A))=(X,(A')_{H_2}').$$

The definition of the action of $F$ on morphisms is more complicated and is
omitted.

\proclaim{Proposition 3}
The algebras $H_1$ and $H_2$ are geometrically equivalent if and only if
the functors
$F(H_1,H_2)$ and $F(H_2,H_1)$ are mutually inverse.\ep

In particular, if $H_1$ and $H_2$ are geometrically equivalent, then the
categories
$K_\Th(H_1)$ and $K_\Th(H_2)$ are isomorphic.

Isomorphism of categories can be established also from the following
commutative diagram
$$\CD C_\Th(H_1) @>>> C_\Th(H_2)\\ @AAA @AAA\\
K_\Th(H_1)@>>> K_\Th(H_2)\endCD$$
Here, the upper arrow is an identity map, and the vertical arrows are dual.

As an example, let us consider the following result which is due to A. Berzins.

\proclaim{Proposition 4}
Two abelian groups $H_1$ and $H_2$ are geometrically equivalent if and only if
\roster\item"1." They have the same exponents,
\item"2." For every prime $p$ the exponents of Sylow subgroups $H_{1p}$ and
$H_{2p}$
coincide.\endroster\ep

\bk

\subheading{2. Geometric equivalence and quasi-identities}

A quasi-identity in $\Th$ has the form
$$w_1\equiv w_1'\wedge\dots\wedge w_n\equiv w_n'\to w\equiv w'.$$
We can assume that all $w_i,w_i',w,w'$ belong to the same $W(X)$ with
finite $X.$

\proclaim {Theorem 1} If   algebras $H_1$ and $H_2$ are
geometrically equivalent, then they have the same
quasi-identities.\ep

In particular, if  algebras $H_1$ and $H_2$ are geometrically equivalent, then
$$\Var(H_1)=\Var(H_2).$$

\proclaim{Corollary}
If two groups $G_1$ and $G_2$ are geometrically equivalent and one of them
is torsion-free,
then the second one is also torsion-free.\ep

\proclaim{Problem 2} Is it true that if $H_1$ and $H_2$ have the
same quasi-identities then they are geometrically equivalent? Is
it true that if $H_1$ and $H_2$ are elementary equivalent, then
they are geometrically equivalent?\ep

There is a negative answer to both of these questions. Namely,
there exists a group $H$ and its ultrapower $\tilde H$ such that
$H$ and $\tilde H$ are not geometrically equivalent. Similar fact
is true also for $G$-groups. This very important expected result
is obtained in \cite{MR2} and uses ideas of the dissertation of
V.A. Gorbunov (1996) \cite{Gor}. For groups there are also
beautiful solutions in \cite{GoSh} and \cite{BlG}.

For associative and Lie algebras this result is also valid (see
\cite{Pl11} for details. For associative algebras it uses the
result from \cite{Li}).

 On
the other hand, recall that in classical situation every extension
$L$ of the field $P$ is geometrically equivalent to every its
ultrapower.

We will consider also generalized (infinitary) quasi-identities.
They have the form
$$\wedge_{\a\in I}(w_\a\equiv w_\a')\to w\equiv w'.$$
Here $I$ is an arbitrary set, and all $w_\a,w_\a',w,w'$ belong to one and
the same $W(X)$
with finite $X.$

\proclaim{Theorem 2} Algebras $H_1$ and $H_2$ are geometrically
equivalent if and only if they have the same generalized
quasi-identities.\ep

This theorem is based on the following form of the theorem about zeroes.

Let $T$  be an arbitrary system of equations in $W(X)$ with finite $X.$
Denote by $I$ the set of all indices $\a$ of equations
$w_\a\equiv w_\a'$ from $T.$ Then
\proclaim{Proposition 5}
Inclusion $w\equiv w'\in T_H''$ is satisfied if and only if
$$\wedge_{\a\in I}(w_\a\equiv w_\a')\to (w\equiv w')$$
is fulfilled in  the algebra $H$.\ep

If $T$ is a finite set, then with $T_H''$ the usual quasi-identity is
associated.
Thus,
\proclaim{Theorem 2$'$}
Algebras $H_1$ and $H_2$ are geometrically equivalent with respect to
finite $T$ if and only
if $H_1$ and $H_2$ generate the same quasivariety.\ep

Note also
\proclaim{Theorem 2$''$}
Algebras $H_1$ and $H_2$ are geometrically equivalent if and only if
$$LSC(H_1)=LSC(H_2).$$
\ep

Here, $L,$ $S,$ and $C$ are the standard closure operators on classes
of algebras, used in the
characterization of prevarieties.

For any class $\frak X$ the class $LSC(\frak X)$ is a locally closed
prevariety over
$\frak X$ which is contained in the quasivariety, generated by $\frak X$
{\cite{PPT}.

For every algebra $H\in\Theta,$ finitely generated algebras in the
prevariety $SC(H)$ are the algebras presented as coordinate algebras
of algebraic sets over $H$.

It follows from \cite{MR2} that the class $LSC(\frak X)$ is not a
quasivariety and, moreover, not an axiomatized class. In this
sense, the relation of geometric equivalence of algebras is not an
axiomatizable relation. This relation is axiomatizable in terms of
generalized quasi-identities \cite{Pl10}.

\bk

\subheading{3. Geometric similarity of algebras}

This notion generalizes the notion of geometric equivalence of algebras,
and, like  the
notion of geometric equivalence, is associated with the problem of isomorphism
of categories of algebraic sets. It is the main notion of the lecture 5.
Geometric equivalence means that the functors $\Cl_{H_1}$ and $\Cl_{H_2}$
coincide.
Geometric similarity assumes the more complicated connection between
$\Cl_{H_1}$ and
$\Cl_{H_2}.$

If $H_1$ and $H_2$ are geometrically equivalent, then $\Var(H_1)=\Var(H_2)$
and the categories $\Var(H_1)^0$ and $\Var(H_2)^0$ are also coincide.

Let us consider the functors
$$\Cl_{H_1}:\Var(H_1)^0\to \text{Set}\quad\text{and}\quad \Cl_{H_2}:
\Var(H_2)^0\to\text{Set}.$$
Similarity of algebras means that there is an isomorphism $\vp:
\Var(H_1)^0\to\Var(H_2)^0$ with
the commutative diagram
$$\alignat 3
&\Var(H_1)^0\quad &&\overset \vp\to\longrightarrow \quad &&\Var(H_2)^0\\
&\ \Cl_{H_1}\searrow && &&\ \swarrow \Cl_{H_2}\\
& && \text{Set}\endalignat$$

Commutativity of the diagram indicates an isomorphism (not necessarily
equality) of the functors
$\Cl_{H_1}$ and
$\Cl_{H_2}\vp. $
This isomorphism $\a=\a(\vp)$ depends on the isomorphism of categories
$\vp$ and is
constructed in a special way.

The notion of geometric similarity appears to be sophisticated, but in
``good'' cases it
reduces to geometric equivalence.

We will use, further, some remarks on congruences in free algebras.

Let the algebra $W=W(X)$ be given.
For a congruence $T$ in $W$, consider a relation $\rho=\rho(T)$ in the
semigroup $\End(W).$
Let us set $\nu\rho\nu',\nu,\nu'\in\End(W),$ if $\nu(w)T\nu'(w)$ holds for
every $w\in W.$

Here, $\rho$-equivalence on $\End(W),$ and $\nu\rho\nu'$ implies
$(\nu_1\nu)\rho(\nu_1\nu')$
for every $\nu_1\in \End(W).$

Let now $\rho$ be an arbitrary congruence on $\End(W).$
Define a relation $T=T(\rho)$ on the algebra $W:$
$$w_1Tw_2\Leftrightarrow\exists w\in W,\nu,\nu'\in\End(W)$$
such that $w^\nu=w_1,$\ $w^{\nu'}=w_2$, and $\nu\rho\nu'.$

\proclaim{Proposition 6}
If $T$ is a congruence on the algebra $W$ and $\rho=\rho(T),$ then $T(\rho)=T.$
\ep

If $T$ is a fully characteristic congruence on $W,$ then $\rho=\rho(T)$ is
a congruence in
the semigroup $\End(W),$ and the semigroup $\End(W)/\rho$ is isomorphic to
the semigroup
$\End(W/T).$

Now let an isomorphism of categories
$$\vp:\Var(H_1)^0\to\Var(H_2)^0$$
be given.

We want to define a function $\a=\a(\vp): \Cl_{H_1}\to\Cl_{H_2}\vp$ which
yields an isomorphism
of functors.
A function $\a$ for every $W$ from $\Var(H_1)^0$ is given by the map
$\a_W:\Cl_{H_1}(W)\to\Cl_{H_2}\vp(W),$ and for every $s:W\to W'$, the
commutative diagram

$$\CD \Cl_{H_1}(W') @>\Cl_{H_1}(s)>> \Cl_{H_1}(W)\\ @V \a_{W'}VV @V
\a_W VV \\
\Cl_{H_2}\varphi(W') @>\Cl_{H_2}\varphi(s)>> \Cl_{H_2}\varphi(W)\endCD$$

should be fulfilled.
This means that $\a=\a(\vp)$ is an isomorphism of functors $\Cl_{H_1}$ and
$\Cl_{H_2}\vp.$
Let us construct such an $\a.$

First for any $\Th$ define functions $\be$ and $\g$ which for an object $W$
from $\Th^0$ give
the maps $\be_W$ and $\g_W.$

The function $\be_W$ is defined by $\be_W(T)=\rho=\rho(T),$ where $T$ is a
congruence in $W$
and $\rho(T)$ is defined above.

The function $\g_W$ is given by
$\g_W(\rho)=T=T(\rho),$ where $\rho$ is an equivalence in $\End(W).$
We have $$\g_W(\be_W(T))=T$$
if $T$ is a congruence in $W.$

For every $T\in\Cl_{H_1}(W)$ define
$$\a(\vp)_W(T)=\g_{\vp(W)}(\vp(\be_W(T)).$$

Here, if $\rho$ is a relation on $\End(W)$ then $\vp(\rho)=\rho^*$ is a
relation on $\vp(W)$
defined by $\mu\rho^*\mu'$ if and only if $\mu=\vp(\nu),$\
$\mu'=\vp(\nu'),$\ $\nu\rho\nu'.$
Here $\a(\vp)_W(T)$ is a relation on the algebra $\vp(W)$ which is not
necessarily a
congruence.

The meaning of the function $\a=\a(\vp)$ is that it represents the action
of the isomorphism
$\vp$ on the congruences of free algebras.

Define now a function $\tau.$
Let $W,W'$ be two objects in $\Th^0.$
Then for any congruence $T$ in $W'$ define $\tau_{W,W'}(T)=\rho$ to be a
relation over the
set $\Hom(W,W')$ defined by
$$s\rho s'\Leftrightarrow w^sTw^{s'},\ \forall w\in W,$$
where $s$ and $s'$ are homomorphisms $W\to W'.$
We have $\tau_{W,W}=\be_W.$
We say that $\a$ is compatible with $\tau$ if for every $W$ and $W',$ and
a congruence $T\triangleleft W',$
$$\vp(\tau_{W,W'}(T))=\tau_{\vp(W),\vp(W')}(\a(\vp)_{W'}(T))$$
holds.

Now we state the main definition.

\definition{Definition 2}
Two algebras $H_1$ and $H_2$ are called {\it geometrically similar} if
there exists an
isomorphism
$$\vp:\Var(H_1)^0\to\Var(H_2)^0,$$
such that the function $\a=\a(\vp):\Cl_{H_1}\to\Cl_{H_2}\vp$
is compatible with the function $\tau$ and for every object $W$ in
$\Var(H_1)^0$ there is a
bijection
$$\a(\vp)_W:\Cl_{H_1}(W)\to\Cl_{H_2}(\vp(W)).$$
\enddefinition
{}From this definition it follows that a function $\a=\a(\vp)$ defines an
isomorphism of
functors $\Cl_{H_1}$ and $\Cl_{H_2}\vp.$

In addition, the same $\a(\vp)$ defines an isomorphism of lattices
$$\a(\vp)_W: \Cl_{H_1}(W)\to\Cl_{H_2}(\vp(W)).$$

\proclaim{Theorem 3}
The algebras $H_1$ and $H_2$ are geometrically similar if and only if there
exists an
isomorphism $\vp:\Var(H_1)^0\to\Var(H_2)^0,$ and an isomorphism of
functors
$\a=\a(\vp):\Cl_{H_1}\to\Cl_{H_2}\vp$ which is
compatible with the
function $\tau.$\ep

Note that
\roster\item"1." If $H_1$ and $H_2$ are geometrically equivalent, then they
are similar.
Here $\Var(H_1)=\Var(H_2)$ and for $\vp$ take the identity map.
\item"2." A relation of geometric similarity is symmetric, reflexive and
transitive.
\item"3." If $\Th_0$ is a subvariety in $\Th,$ and $H_1,H_2\in\Th_0, $ then
$H_1$ and $H_2$
are similar in $\Th$ if and only if they are similar in $\Th_0.$\endroster

\bk

\subheading{4. Similarity and equivalence}

Consider the case when $\Var(H_1)=\Var(H_2)=\Th.$
In this case $\vp$ is an automorphism of the category $\Th^0.$

\definition{Definition 3}
An automorphism $\vp: C\to C$ of an arbitrary category $C$ is called {\it
inner} if it
corresponds to a function $s$ which attaches an isomorphism
$s_A:A\to\vp(A)$ to every object
$A$ from $C$ such that for any $\nu: A\to B$ there is a commutative diagram

$$\CD A @>\nu>> B\\ @V s_AVV @Vs_B VV\\
\vp(A) @>\vp(\nu)>> \vp(B)\endCD$$

i.e.,
$$\vp(\nu)=s_B\nu s_A^{-1}.$$
\enddefinition

From this definition follows that an automorphism $\varphi$ is inner
if and only if the functor $\varphi$ is isomorphic to identity functor
of the category $C$.

\proclaim{Theorem 4}
If the algebras $H_1$ and $H_2$ are geometrically similar with respect to
the inner
automorphism $\vp$, then they are equivalent.\ep

Let us make some remarks on inner automorphisms of categories.

If an inner automorphism $\vp$ is realized by a function $s,$ then we write
$\vp=\hat s.$
Let now $\vp$ be an arbitrary automorphism and $\hat s$ be an inner
automorphism.
Take an automorphism $\vp\hat s\vp^{-1}$ and consider a function $\vp(s)$
defined by the rule
$$\vp(s)_B=\vp(s_{\vp^{-1}B})$$
for any object $B$ in the given category $C.$

Let us check that
$$\vp\hat s\vp^{-1}=\widehat{\vp(s)}.$$
Thus if $\Aut(C)$ is a group of all automorphisms of the category $C,$ then
there is a normal
subgroup $\Int(C)$ of all inner automorphisms.

Let now $\vp\in\Aut(C)$ and $A$ be an arbitrary object in $C.$
Then $\vp$ induces an isomorphism of semigroups $\End(A)$ and
$\End(\vp(A))$ and the groups
$\Aut(A)$ and $\Aut(\vp(A)).$
If $\vp(A)=A$ then $\vp$ induces an automorphism of the semigroup $\End(A)$
and of the group
$\Aut(A).$
If $\vp$ is an inner automorphism, then the corresponding induced
automorphism of the
semigroup $\End(A)$ and the group $\Aut(A)$ is also inner.

\definition{Definition 4}
A variety $\Th$ is called {\it perfect} if every automorphism of the
category $\Th^0$ is an
inner automorphism.\enddefinition

We will consider the examples of perfect $\Th.$
We can state that if $\Th$ is a perfect variety and
$\Var(H_1)=\Var(H_2)=\Th,$ then $H_1$ and
$H_2$ are similar if and only if they are equivalent.

Let us define the weakly inner automorphisms $\vp: C\to C.$
An automorphism $\vp$ is called {\it weakly inner} if there exists a
function $s$ with
isomorphisms $s_A:A\to\vp(A)$, such that for every $\nu: A\to A$ holds
$$\vp(\nu)=s_A\nu s_A^{-1}:\vp(A)\to\vp(A).$$

It can be proven that in Theorem 4 one can proceed from a weakly inner
automorphism.
Using this fact, one can prove

\proclaim{Proposition 7}
Two abelian groups $H_1$ and $H_2$, such that each of them generates the
whole variety of
abelian groups, are similar if and only if they are equivalent.\ep

A reasonable conjecture states that the same is true for $H_1$ and $H_2$,
generating the
variety of all groups Grp.

The solution of this problem is closely connected with the following problem.

\proclaim{Problem 3}
Let $F=F(X)$ be the non-commutative free group with finite $X.$
It is known that the group $\Aut(F)$ is a perfect group.
Is the same true for the semigroup $\End(F)?$
\ep

Let $\vp$ be an automorphism of the semigroup $\End(F).$
This $\vp$ induces an automorphism of the group $\Aut (F).$
Thus, there exists an element $\sigma\in\Aut(F),$ such that $\vp(x)=\sigma
x\sigma^{-1}$ for
every $x\in\Aut(F).$
Is it true that the same holds for every $x\in\End(F)?${\footnote{
This question is solved positively in \cite{For}. From this follows that
the conjecture is also true. From the result of Formanek follows that
the variety $Grp$ is weekly perfect. Thus, arises the following problem.
Is it true that the variety $Grp$ is perfect?  The positive answer can be
obtained using the technique of Formanek. The same question stands
for the variety of all abelian groups.
}

\bk

\subheading{5. Similarity and semiisomorphisms}

Let $\Th$ be an arbitrary variety, $G\in\Th,$ and consider the variety $\Th^G.$
Assume that for $G\in\Th^G$ the condition (*) (see Lecture 2) is
fulfilled.
This means that if $(H_1,h_1)$ and $(H_2,h_2)$ are two faithful
$G$-algebras, then each of
them generates the variety $\Th^G.$

In the category $\Th^G$ one can consider semiisomorphisms and semiisomorphic
algebras.

\proclaim{Theorem 5}
If algebras $(H_1,h_1)$ and $(H_2,h_2)$ are semiisomorphic, then they are
similar.\ep

Thus, if there is a sequence of $G$-algebras
$$(H_1,h_1)=(H_1^0,h_1^0),(H_1^1,h_1^1),\dots (H_1^n,h_1^n)=(H_2,h_2),$$
such that the neighbors are geometrically equivalent or semiisomorphic,
then $(H_1,h_1)$ and
$(H_2,h_2)$ are similar.

\definition{Definition 5} Algebras $(H_1,h_1)$ and $(H_2,h_2)$ are called
geometrically
equivalent up to a semiisomorphism, if there exists $(H,h)$ such that
$(H_1,h_1)$ and $(H,h)$
are semiisomorphic, and $(H,h)$ is equivalent to $(H_2,h_2).$
\enddefinition

We will consider semi-inner automorphisms of the category
$(\Th^G)^0.$

Let $\vp:(\Th^{G})^0\to(\Th^G)^0$ be an automorphism of the category of
free algebras in
$\Th^G.$
These algebras have the form $$W(X)=G\ast W_0(X)=W=G\ast W_0.$$
Consider a pair $(\sigma,s)$, where $\sigma$ is an automorphism of the
algebra $G$ and $s_W:
W\to\vp(W)$ is an algebra isomorphism in $\Th$ for every $W$.
The commutative diagram

$$\CD G @>i_G>> W\\ @V \sigma VV @V s_W VV \\
G @>i_G'>> \vp(W)\endCD$$

defines a semiisomorphism
$$(\sigma,s_W):(W,i_G)\to(\vp(W),i_G').$$

It is supposed that with the automorphism $\vp$ one can associate such a
pair $(\sigma,s)$
and, besides, for every morphism $\nu: W\to W'$ there is a commutative diagram

$$\CD W @>\nu>> W'\\ @V (\sigma,s_W) VV @V (\sigma,s_{W'})VV\\
\vp(W) @>\vp(\nu)>> \vp(W')\endCD$$

Then,
$$\vp(\nu)=(\sigma,s_{W'})(1,\nu)(\sigma,s_W)^{-1}=(1,s_{W'}\nu
s_W^{-1})=s_{W'}\nu s_W^{-1},$$
and this is a morphism in $(\Th^G)^0.$

\definition {Definition 6}
An automorphism $\vp$ is called {\it semi-inner}, if there exists
a pair $(\sigma,s)$ with the property above.\enddefinition

We denote $\vp=(\widehat{\sigma,s}).$

It can be proved that semi-inner automorphisms constitute a
subgroup in $\Aut((\Th^G)^0)$ which contains the normal subgroup
of inner automorphisms.

From the definition 6 follows that an automorphism $\varphi$ is
semi-inner if and only if the functor $\varphi$ is semiisomorphic
to the identity functor of the category $(\Theta^G)^0$.

The main theorem here is as follows:

\proclaim{Theorem 6} If the similarity of algebras $(H_1,h_1)$ and
$(H_2,h_2)$ is given by a semi-inner automorphism $\vp,$ then they
are equivalent up to a semiisomorphism.\ep

\bk

\subheading{6. Geometrically Noetherian algebras}

We define a property of an algebra $H\in\Th,$ which always takes place in
the classical
situation and plays a crucial part.

\definition{Definition 7}
An algebra $H$ in the variety $\Th$ is called {\it geometrically Noetherian}
if for every
finite set $X$ and every system of equations $T$ in $W=W(X),$ there exists
a finite
subsystem $T_0\subset T,$ such that $T_H'=T_{0H}',$ or, in other words,
$T_H''=T_{0H}''.$

We say that $T$ and $T_0$ are an equivalent systems of equations.\enddefinition

\proclaim{Proposition 8}
An algebra $H$ is geometrically Noetherian if and only if for every finite
$X$ the lattices
$\Cl_H(W(X))$ and $\Alv_H(W(X))$ satisfy the maximal and minimal condition,
respectively.\ep

\proclaim{Proposition 9}
Let $\Th_0$ be a subvariety in $\Th$ and $H\in\Th^0.$
Then $H$ is geometrically Noetherian in $\Th$ if and only if it is
geometrically Noetherian
in $\Th_0.$
\ep

\definition{Definition 8}
A variety $\Th$ is called {\it Noetherian} if every finitely generated
algebra from $\Th$ is
Noetherian.\enddefinition

\example{Examples}
1. The variety $Var-P$ is noetherian.

2. Arbitrary variety of nilpotent groups  is noetherian.

3. In the variety of arbitrary associative algebras over the field $P$
   every noetherian subvariety can be distinguished by special known
   identities.
\endexample

\proclaim{Proposition 10} If the variety $\Th$ is Noetherian, then
every algebra $H$ in $\Th$ is geometrically Noetherian. In
particular, if the variety $\Var(H)$ is Noetherian, then the
algebra $H$ is geometrically Noetherian.\ep

\proclaim{Theorem 7} Let $H_1$ and $H_2$ be two geometrically
Noetherian algebras. They are geometrically equivalent if they
have the same quasi-identities.\ep

\proclaim{Theorem 8} Let $H_1$ and $H_2$ have the same
quasi-identities. Then, if $H_1$ is geometrically Noetherian, then
$H_2$ is  also geometrically Noetherian.\ep

\proclaim{Corollary}
If an algebra $H$ is geometrically Noetherian, then each of its cartesian
power and
ultrapower is geometrically Noetherian.\ep

The same is true for arbitrary filtered power.

\proclaim{Proposition 11}
If algebras $H_1$ and $H_2$ are geometrically similar and one of them is
geometrically
Noetherian, then the second one is also geometrically Noetherian.\ep

It can be easily seen that every subalgebra of a geometrically Noetherian
algebra is also
geometrically Noetherian, and a finite cartesian product of geometrically
Noetherian algebras
is also a Noetherian algebra.
Every finite algebra is geometrically Noetherian.

Geometrically Noetherian algebras admit a Noetherian Zariski topology in
the corresponding
affine spaces and the Lasker-Noether Theorem on the decomposition of
algebraic varieties.

\definition{Definition 9}
A variety of the type $\Th^G$ is called {\it faithfully Noetherian} if it
is Noetherian with
respect to faithful congruences.\enddefinition

If the variety $\Th^G$ is faithfully Noetherian, then each of its faithful
algebras $(H,h)$ is
geometrically Noetherian.

There are many important results on geometrically Noetherian
groups (see \cite{Br}, \cite{Guba}, \cite{BMR1}, \cite{BMRo} and
others. There are also many interesting problems, especially for
the case of $G$ in $\Th^G.$

\bk

\subheading{7. Geometric stability}

This notion was already defined and will be used and discussed in the next
lecture.
Let us make now some remarks on $G$-algebras.

\proclaim{Proposition 12}
Let the algebras $(H_1,h_1)$ and $(H_2,h_2)$ be semiisomorphic.
Then,
\roster\item"1." If one of them is stable, then the second one is also stable.
\item"2." If one of them is geometrically Noetherian, then the second one
is also geometrically Noetherian.
\endroster\ep

Besides that, if $(H_1,h_1)$ and $(H_2,h_2)$ are two arbitrary
$G$-algebras, then their
product $(H_1\times H_2,h_1\times h_2)$ cannot be a stable algebra \cite{Be}.
{}From this follows that if $(H_1,h_1)$ and $(H_2,h_2)$ are equivalent and
one of them is
stable, then the second one is not necessarily stable.

\definition{Definition 10} An algebra $H$ in the variety $\Th$ is called {\it
geometrically
distributive} if every lattice $\C_H(W)$ is distributive.
\enddefinition

If $H$ is stable, then $H$ is geometrically distributive.

\proclaim{Problem 4}  Is it true that an algebra $H$ is geometrically
distributive if and only
if it is similar to a stable algebra?\ep

One can consider also geometrically modular algebras.

It can be proven that if an algebra $G\in\Th$ admits faithful finite
dimensional
linearization over a field, then it is geometrically Noetherian in $\Th$
and in $\Th^G.$
This fact, in particular, relates to finite dimensional, associative and
Lie algebras.
(For groups, see \cite{BMR}, the general case is done by A.Belov
\cite{BelP}).
The proof uses the natural Zariski topology over the same field.

A.Miasnikov and V.Remeslennikov have noticed the following
interesting question. Is it true that a free Lie algebra $F(X)$
with finite $X$ is geometrically noetherian?

\newpage

\baselineskip 12pt
 \topmatter \rightheadtext{Lecture 4}
\leftheadtext{B. Plotkin}
\title Lecture 4\\
\quad\\
Algebraic geometry in group-based algebras.\endtitle
\endtopmatter

\bigskip
\centerline{\smc Contents}
\bigskip
\roster\item "1." {\smc Group-based algebras }
\newline
\item"2."{\smc  Zero divisors and nilpotent elements}
\newline
\item"3." {\smc Group-based algebras with an algebra of constants }
\newline
\item"4." {\smc  Domains}
\newline
\item"5." {\smc Algebraic geometry in group-based algebras}
\newline
\item"6." {\smc Stability}
\newline
\item"7." {\smc Theorems on zeroes}
\endroster

\newpage

\par\newpage\par\flushpar

\bigskip
\subheading{1. Group-based algebras}

\definition{Definition 1}
A group-based algebra $H$ is an additive group (not necessarily
 commutative) with some additional signature $\Omega$.
For every $\omega \in \Omega$ of the arity $n(\omega) = n>0$ the condition
$0 \cdots 0 \omega = 0$ should be fulfilled.
\enddefinition

Here $0$ is the zero element in the additive group $H$ and
on the left side we have
 $n$ times for $0$.

Such group-based algebras are called also $\Omega$-groups.
They were introduced by P. Higgins \cite{Hi} in 1956, see also \cite{Ku},
\cite{Pl7}.

A group is an $\Omega$-group with empty set $\Omega$,
in rings the set $\Omega$ consists
  of a single multiplication, in modules over a ring $R$ all elements of $R$
belong to the set $\Omega$. Groups over rings, considered by R.Lyndon
\cite{L2} are also
$\Omega$-groups.

In arbitrary $\Omega$-group $G$ we have the usual commutator $[a, b] = - a-b+a+b
= -a+a^b$,
and $\omega$-commutators for $\omega$ with $n(\omega) > 0$.
By definition we have:
$$
[a_1, \cdots, a_n; b_1, \cdots, b_n; \omega]= - a_1\cdots a_n\omega-b_1\cdots b_n\omega+(a_1+
b_1) \cdots (a_n + b_n) \omega.
$$
From this definition follows that if all $a_1, \cdots, a_n$ are zeros or all
$b_1, \cdots, b_n$ are zeros then the $\omega$-commutator is a zero.

In rings we have:
$$
[a_1, a_2; b_1, b_2; \cdot] = a_1b_2+b_1a_2.
$$
Now about ideals in $\Omega$-groups.
Let $\mu\colon H \to H_1$ be a homomorphism of two $\Omega$-groups, ${U}$
 is the coimage of the zero element from $H_1$.
Then:

\noindent 1) \   $U$ is closed relative to all $\omega\in\Omega$ with $n(\omega)>0$.

\noindent 2) \  $U$ is a normal subgroup in additive group $H$.

\noindent 3) \  The $\omega$-commutator
$[a_1, \cdots, a_n; b_1, \cdots, b_n; \omega]$  belongs to $U$
if $a_1, \cdots, a_n
 \in U$, and  $ b_1, \cdots, b_n$ are arbitrary elements of $H$.

\definition{Definition 2}  \ An ideal $U$ in $H$ is a set $U$ with the
three conditions
above.
\enddefinition

From the definition follows also that the commutator  $[a_1, \cdots, a_n; b_1, \cdots,
 b_n; \omega]$ belongs to $U$ if all $b_1, \cdots, b_n$ in $U$ and $a_1, \cdots, a_n$-
are arbitrary.
The ideal in a group is a normal subgroup, the ideal in a ring is a usual
ideal.
It can be proved that if $T$ is a congruence in an $\Omega$-group $H$,
then the class
 $[0]=U$ for $T$ is an ideal and we have:
$[a] = a+U$ for arbitrary $a\in H$.
So we have the one-to-one correspondence  between congruences and ideals.
We write now  $\Ker \mu = U$  instead of $\Ker \mu = T$ and $H/U$ instead of $H/T$.

Now let $A$ and $B$ be two sets in the $\Omega$-group $H$, $\{ A, B\}$
be an $\Omega$-subgroup
in $H$, generated by $A$ and $B$.

\definition{Definition 3} The mutual commutant
$[A, B]$ is the ideal in $\{A, B\}$,
generated by all commutators of the kind
$[a, b], a \in A, b \in B$ and
all $[a_1,\cdots, a_n;b_1,\cdots, b_n;\omega]$, were
 $a_1, \cdots, a_n \in A, b_1, \cdots, b_n \in B$.
\enddefinition

We have $[A, B]= [B, A]$ and for every two ideals $U_1$ and $U_2$ we have:
$[U_1, U_2]\subset U_1\cap U_2$.
Let now $U_1$ and $U_2$ be ideals in the associative ring $R$.
Then:
$$
[U_1, U_2]=U_1U_2+U_2U_1.
$$
In  the Lie case we have: $[U_1, U_2]= U_1U_2$, in groups it is usual
mutual commutant.
\definition{Definition 4}  The group-based algebra $H$ is abelian if $[H,H]=0$.
\enddefinition

This means that the group $H$ is abelian and for every
$\omega \in \Omega$, $n=n(\omega)>0$
 we have:
$$
\eqalign{
&(a_1+b_1)\cdots (a_n+b_n)\omega=\cr
&=a_1\cdots a_n\omega+b_1\cdots b_n\omega\cr}
$$
for every $a_1, \cdots, a_n$; $b_1, \cdots, b_n$.
In rings this means that $ab=0$ for every $a$ and $b$.
\definition{Definition 5}  A group-based algebra $H$ is antiabelian if:

\noindent
1. \ If $U$ is an abelian  ideal in $H$ then $U=0$.

\noindent
2. \ If  $U_1$ and $U_2$ are two nontrivial ideals then
the ideal $U_1 \cap U_2$ is also
non-trivial.
\enddefinition

\subheading{2. Zero divisors and nilpotent elements}

From now on all algebras are group-based algebras.
Let an algebra $H$ be given, $a \in G$
and $(a)$ be the ideal in $H$, generated by the element $a$.

\definition{Definition 6} A non-zero element $a\in G$ is called a zero
divisor if for some non-zero
 $b\in G$ we have
$$
[(a), (b)]= 0.
$$
\enddefinition

Consider this notion for associative and Lie rings and groups.
We are studying the case without zero divisors.

\proclaim{Proposition 1} A ring $R$ is without zero divisors if
and only if for every non-zero
 $a$ and $b$ there exist $c_1, c_2, c'_1, c'_2$ such that
$$
c_1ac_2\cdot c'_1bc'_2\neq 0 \; \; \hbox{\rm or} \; \; c'_1bc'_2\cdot
c_1 a c_2 \neq 0.
$$
\endproclaim
\proclaim{Corollary} The class of all $R$ without zero divisors is an axiomatizable
class.
\endproclaim

\proclaim{Proposition 2} The Lie ring $L$ is without zero divisors
if and only if  for every non-zero $a$ and $b$ in $L$ there exist
$m \ge 0$, $ n \ge 0$, and elements $c_1, \cdots, c_m,
 c'_1, \cdots, c'_n\in L$ such that:
$$
[[a, c_1, \cdots, c_m]; \; [b, c'_1, \cdots, c'_n]] \neq 0.
$$
\endproclaim
\proclaim{Corollary}  The class of all $L$ without zero divisors, may not be an axiomatizable
class.
\endproclaim

\proclaim{Proposition 3}  The group $G$ is a group without zero divisors iff for every non-unite
 $a$ and $b$ in $G$ there exists $c \in G$ such that:
$$
[a, b^c] \neq 1,
$$
where  $b^c = c^{-1} bc.$
\endproclaim
\proclaim{Corollary} The class of all groups without zero divisors is an axiomatizable
class.
\endproclaim
Groups without zero divisors were considered in [BMR] and
in [BPP] under the name
antiabelian. In [BMR] it was proved (in additional conditions)
that the free product of two groups without
zero divisors is a group without zero divisors as well.

From the corollary above we can conclude that the ultraproduct of
groups without
 zero divisors is also a group without zero divisors.

Now the following general result.

\proclaim{Proposition 4}  An $\Omega$-group $H$ does not have
zero divisors iff $H$
 is antiabelian.
\endproclaim

{\bf Examples.} An algebra $H$ is simple if $H$ has not  non-trivial
ideals.
Every non-abelian and simple algebra is antiabelian and so without zero-divisors.
Every free group and every free Lie algebra is antiabelian.
Every free associative or free associative and commutative algebras
are without zero divisors,
 and thus such an algebras are antiabelian.
It is easy to characterize all antiabelian
finite dimensional associative and Lie algebras.

Note that using commutants in $\Omega$-groups the nilpotent and solvable
$\Omega$-groups are naturally defined.

\definition{Definition 7} An element $a$ in an $\Omega$-group $H$ is called
strictly nilpotent if the ideal $(a)$ is nilpotent.
 An element $a$ in an $\Omega$-group $H$ is called weakly nilpotent
if the ideal $(a)$ is solvable.
\enddefinition

It is easy to see that if $a$ is a weakly nilpotent element then
in $(a)$ there is a zero divisor. If  $a$ is a strictly nilpotent
element then $a$ is a zero divisor.

\subheading{3. Group-based algebras with an algebra of constants}

For an arbitrary variety $\Theta$ and an arbitrary $G\in \Theta$
we consider a new variety of algebras with constants from $G$, which
is denoted by $\Theta^G$.
We take for $\Theta$ a variety of $\Omega$-groups and let $G$ be
a fixed $\Omega$-group from $\Theta$. As before, the key point for
considering such a variety is the fact that groups and rings in their  pure form are not
so convenient for constructing an algebraic geometry.
For algebraic geometry we need some
large set of constants in the signature $\Omega$, and
we need equations with constants.

For a given $h\colon  G \to H$ the ideal $U$ in $H$ is the
same as the ideal in $G$-algebra $H$.
For such $U $ we have the embedding $\bar h \colon G \to H/U$, induced
 by $h$.
If $h$ is faithful and $Im (h) \cap U=0$, then $\bar h$ is also faithful.

Let us define now relative ideals and relative zero divisors in respect
to the group $G$.

Now we consider the notion of relative ideal or $G$-ideal.
We connect the notion of ideal with the algebra of constants
$G$.
\definition{Definition 8} The set $U\subset H$ in the algebra $h\colon G \to H$ is a
$G$-ideal if
\pmf
1) \ $U$ is closed in respect to all $\omega \in \Omega$, $n(\omega) >0$, and
$U$ is a subgroup in the additive group $H$,
\pmf
2) \ For every $a \in {U} $ and $g \in G$ we have $[a, h(g)] \in U$,
\pmf
3) \ $[a_1, \cdots, a_n; h(g_1), \cdots, h(g_n); \omega] \in U$ \
if $\omega \in \Omega, n(\omega) >0$; $a_1, \cdots, a_n \in U,$ $ g_1, \cdots, g_n \in G$.
\enddefinition

If the Conditions 1,2,3, are fulfilled the set $U$ is called
invariant in respect to constants from $G$.

From the conditions 1), 2), 3) we have also:
$$
[h(g_1), \cdots, h(g_n); a_1, \cdots, a_n; \omega] \in U.
$$

Every ideal at the same time is an $G$-ideal.

The algebra $G$ can be considered as a $G$-algebra with the identical
$G\to G$.
In this case the $G$-ideal in $G$ is the same as ideal.

For every $a \in H$ by $(a)^{G}$ we denote the $G$-ideal, generated
 by $a$.
\definition{Definition 9}  A non-zero $a\in H$ we call a $G$-zero divisor
if for some
$ b \neq 0$  we have
$$
[(a)^{G}, (b)^{G}]=0.
$$
\enddefinition

Every zero divisor in $H$ is also a $G$-zero divisor.
So if $H$ is without $G$-zero divisors then $H$ is without (absolute)
zero divisors.
Such $H$ is an antiabelian.

\proclaim{Proposition 5}  Let $h\colon G\to H$ be a $G$-group.
Then $H$ is without $G$-zero divisors if and only if for every
non-unite $a$ and $b$
 there exists $c \in G$ such that
$$
[a, b^{h(c)}]\neq 1.
$$
\endproclaim

We can consider now a set $T$ of equations
$$
[x, y^{h(c)}] = 1,\qquad c \in G.
$$
Let $T' = A$ be the set of all solutions in the given $H$.
We can say that $H$ is without $G$-zero divisors iff from $(a, b) \in A$
follows $a= 1 $ or $b=1$.
In other words the infinitary formula
$$
\mathop{\wedge}\limits_{c\in G}( [x, y^{h(c)}]=1) \to (x=1) \vee (y=1).
$$
holds in $H$.
The similar can be considered in Lie-$L$ and in $Assoc-R$.

An algebra $G \to G$ is without $G$-zero divisors iff $G $ is
 antiabelian.

{\bf Examples.}  As before, if $\Theta = Grp$,
 then we write
$\Theta^G = Grp-G$.
In the similar Lie case we have Lie-$L$, and in associative one
 $\Theta^G =Assoc-R$.

\subheading{4. Domains}

 A $G$-algebra $H$ is called a domain if
$H$ is  without
 $G$-zero divisors.
It is trivial and very important that every $G$-subalgebra
in the domain $H$ is also
 a domain.

Now let  $h_1\colon G \to H_1$ and $h_2 \colon G \to H_2$ be two
faithful $G$-algebras.
We have in this case an isomorphism $h_0 \colon Im (h_1) \to Im (h_2)$.
Using this $h_0$ we can consider the amalgamated free product
$$
h\colon G \to {H_1}*_{G} H_2.
$$
That is a free product in $\Theta(G)$.
In [BMR] it is proven  that if $H_1$ and $H_2$ from
$Grp-G$ are domains then their
 free product is also a domain.
\definition{Definition 10}  An ideal $U$ in $G$-algebra $H$ is called prime if
$H/U$ is a domain.
\enddefinition

For every $G$-algebra $H$ denote by $Spec (H)$ the set of all prime ideals
 in $H$.

If $H/U$ is a domain then this algebra is antiabelian.
In this case the ideal $U$ is irreducible, $U$ cannot be
non-trivially represented as $U=U_1 \cap U_2$.

Let now $\mu \colon H \to H_1$ be a homomorphism of $G$-algebras and $H_1$
 be a domain.
Then $U = \Ker \mu$ is prime.

\proclaim{Proposition 6}  Let $H$ and $H_1$ be two semiisomorphic
$G$-algebras. Then $H_1$ is a domain if and only if $H$ is a
domain.
\endproclaim
\proclaim{Proposition 7}  Let
$\mu \colon H \to H_1$ be a homomorphism of $G$-algebras, $a\in H$.
Then we have
$$
((a)^G)^\mu=(a^\mu)^G.
$$
\endproclaim
This Proposition plays an important role in the proofs of Theorems
1 and 2 below.


\subheading{5. Algebraic geometry in group-based algebras}

Let us repeat the basic notions for the case when the variety
$\Theta$ is a variety of $\Omega$-groups.

Let $X$ be a finite set and $W=W(X)$ a free in $\Theta$ algebra over $X$.
For the given $H \in \Theta$ the set of
homomorphisms $Hom (W, H)$ we consider as an
 affine space.

For every point $\mu \colon W \to H$ we have the
kernel $U=Ker \mu$ and $\mu$ is a solution of some equation $w \equiv 0$ iff
 $w \in Ker \mu$.

Now let $U$ be a set of ``polynomials" in $W$ and $A$ be a set of points in the space
 $Hom (W, H)$.
We establish the following Galois correspondence:
$$
\cases
U'=A=\{ \mu | U \subset Ker \mu \}\cr
A'=U=\bigcap\limits_{\mu \in A} Ker \mu \endcases
$$
The set $A$ such that $A=U'$ for some $U$  we call an algebraic set
over the algebra $H$. The  ideal $U$ of the form $U=A'$ for some $A$
we call an  $H$-closed ideal.

As usual we call algebraic sets also algebraic varieties.

Now let two algebraic varieties $A$ and $B$ in $Hom (W, G)$ be
given, $A'=U_1$ and $B'=U_2$. Then: $(A\cap B) = (U_1\cup U_2)'$
and so the intersection of two algebraic varieties  is also an
algebraic variety. But we can not say the same for the union
$A\cup B$. In general, it is not true that $A\cup B = (U_1 \cap
U_2)'$. It depends on the choice of the algebra $H$.

\subheading{6. Stability}

\definition{Definition 11}  The algebra $H$ is said to be {\sl stable}
if for every $W = W(X)$
 and every two algebraic sets $A$ and $B$ in the space $Hom (W, H)$ the union
 $A\cup B$ is also an algebraic set.
\enddefinition

If $H$ is stable then in the category $K_\Theta (H)$  every object $A$
can be considered
 as a topological space.
The morphisms in $K_\Theta(H)$ are well coordinated with the topology,
 they are continuous maps.

In the classical case all fields and all domains are stable.
But if $\Theta = Grp$ is the variety of  all groups and $H$ is abelian
then $H$ is not stable.

Now the main question: when the algebra $H$ is stable?
And the main idea:  the notions to be stable and without zero divisors
are close to each other.

Consider now the varieties of the type $\Theta^G$.
\proclaim{Theorem 1} If a $G$-algebra $H$ is
a domain then $H$ is stable.
\endproclaim

Now we consider the following special identities,
the so-called $CD$-identities (commutator
 distributivity).
They are defined for all
$\omega \in \Omega$, $n=n(\omega)>0$ and have the form:
$$
\eqalign{
&[x_1+y_1, \cdots, x_n +y_n; z_1, \cdots, z_n; \omega]=\cr
&=[x_1, \cdots, x_n; z_1, \cdots, z_n; \omega]+\cr
&+[y_1, \cdots, y_n; z_1, \cdots, z_n; \omega].\cr}
$$
Here all $x_i, y_i, z_i $ in $X$.
\definition{Definition 12} An algebra $H \in \Theta$ is called a $CD$-algebra
if
 all $CD$-identities are hold in $H$.
\enddefinition
\proclaim{Proposition 8}  If the algebra $G$ is a $CD$-algebra and
$G$ has zero-divisors,
 then $G$ is not stable in $\Theta^G$.
\endproclaim

\proclaim{Theorem 2} Let $G$ be a $CD$-algebra in $\Theta^G$.
Then $G$ is stable if and only if $G$ has no zero divisors.
\endproclaim

$CD$-conditions are fulfilled in groups, associative and Lie algebras.
In these cases the theorem was proved by A.Berzins. The result is true also
in arbitrary rings.

\subheading{7. Theorems on zeroes}

We consider the case when the algebra $G$ is algebraically closed in
the variety $\Theta^G$.
Let first the variety $\Theta$ be an arbitrary variety, i.e.,
not necessarily the variety of group-based algebras.

If $T$ is a maximal faithful congruence in the
algebra $W=G\ast W_0$, then $T=Ker\mu$ for some point
$\mu$.

$G$-algebra $G$ is a simple algebra since it has
no faithful congruences.

\definition{Definition 13} $G$-algebra $H$ is called semisimple if
it is approximated by simple $G$-algebras.
$G$-algebra $H$ is called locally semisimple,
if every finitely generated subalgebra
in $H$ is simple.
\enddefinition

\proclaim{Theorem 3} Let algebra $G$ be algebraically closed
in the variety $\Theta^G$ and let $H$ be a faithful
locally semisimple $G$-algebra. Then for every
faithful congruence $T$ in $W=W(X)$ with finite
$X$, the congruence $T''_H$ is the intersection
of all maximal faithful congruences in $W$
which contain $T$.
\endproclaim

{Proof} By definition, $T''_H$ is the intersection
of all $Ker\mu$ where $T\subset Ker\mu$. For every such
$\mu$ take in $H$ the image $Im \mu=H_1$.
The algebras $H_1$ and $W/Ker \mu$ are
faithful and semisimple. Therefore, the congruence
$Ker \mu $ is the intersection of some
maximal faithful congruences. Then $T''_H$
is an intersection of maximal faithful congruences.

Let now $T_1$ be a maximal faithful congruence in
$W$ containing $T$. Show that $T_1$
contains $T''_H$. We have $T_1=Ker \mu$ for
some point $\mu : W\to G$. This point $\mu$
is a homomorphism $W\to H$. This means
that
$T_1=T''_{1H}$. Since $T\subset T_1$ then
$T''_H\subset T''_{1H}$. This implies the theorem.

\proclaim{Corollary}
If $H_1$ and $H_2$ are faithful locally simple
$G$-algebras, and $G$ is algebraically closed, then
$H_1$ and $H_2$ are geometrically equivalent.
\endproclaim

Let us come back to the situation of group-based algebras.
Define relative nilpotents, i.e., $G$-nilpotents.

Let $(H,h)$ be a $G-\Omega$-group and $a\in H$.

\definition{Definition 14} An element $a$ is called
$G$-nilpotent if in $\Omega$-group $(a)^G$ there is
a series of $G$-invariant ideals
$$
0=U_0\subset U_1\cdots \subset U_n=(a)^G
$$
with the abelian factors $U_{i+1}/U_i$. If ,
moreover, $[(a)^G,U_{i+1}]\subset U_i$
then $a$ is said to be strictly $G$-nilpotent
element.
\enddefinition

In groups ($\Omega$ is empty) an element $a$ is
a nilpotent element if the group $(a)^G$
is solvable. If $(a)^G$ is nilpotent then
$a$ is strictly nilpotent.

In associative rings both these notions coincide.
Namely, an element $a\in H$ is a nilpotent if there
exists $n$ such that
$$
b_0ab_1ab_2a\cdots b_{n-1}ab_n=0
$$
for any $b_i=h(g_i)$, $g_i\in G$.

If the ring is commutative then this condition
is equivalent to $a^n=0$.

If the algebra $H$ is a domain then it has no non-zero
nilpotent elements.

Given $G$-algebra $H$ denote $N(H)$ the ideal
in $H$ generated by all its $G$-nilpotent elements.
This is the $N$-radical of $H$. The function $N$ can be iterated
and this gives rise to the upper $N$-radical.
Denote it by $\tilde{N}(H)$. There are no
non-trivial nilpotent elements in $H/\tilde{N}(H)$.

If $U$ is an ideal in $H$ then denote by $\sqrt U$
the inverse image in $H$ of the radical
$\tilde{N}(H/U)$.

It is clear that if a $G$-homomorphism
$\mu: H\to H' $ is given and $U\subset Ker \mu$
then $\sqrt U$ is also contained in $Ker \mu$ in case
$H'$ is a domain.
This implies

\proclaim{Proposition 9}  Let $U$ be an ideal
in the free $G$-algebra $W=G\ast W_0$, and let
$H$ be a $G$-domain. Then
$$
\sqrt U\subset U''_H.
$$
\endproclaim

We are looking for conditions which
imply the equality in the formula above.

\definition{Definition 15} A variety
is called special, if for every finitely generated
$G$-algebra $H$ its radical $\tilde{N}(H)$
is an intersection of maximal faithful
congruences in $H$.
\enddefinition

\proclaim{Theorem 3$'$ } If an algebra $G$ is
algebraically closed and the variety
$\Theta^G$ is special then for every faithful $G$-domain $H$
the equality
$$
\sqrt U = U''_H.
$$
holds for every ideal $U$ and every $W$.
\endproclaim

In fact, this theorem, as well as the theorem 3,
immediately follows from definitions. The problem is to
study the situations when all this is
applicable.

Here the ideas of radical, semisimplicity and algebraic closure are
naturally intersected in the theorems 3 and 3$'$.

\newpage

\newpage

\baselineskip 12pt
 \topmatter \rightheadtext{Lecture 5}
\leftheadtext{B. Plotkin}
\title Lecture 5\\
\quad\\
Isomorphisms of categories of algebraic sets
\endtitle
\endtopmatter

\centerline{\smc Contents}
\bigskip
\roster\item "1." {\smc  Correct isomorphism}
\newline
\item"2." {\smc Isomorphism, similarity, equivalence}
\newline
\item"3." {\smc  Perfect and semiperfect varieties of algebras}
\newline
\item"4." {\smc  Other results}
\newline
\item"5." {\smc  Problems}
\newline
\endroster


\subheading{1. Correct isomorphism}

 Let $H_1$ and $H_2$ be two algebras in the variety
$\Theta$. We are looking for an isomorphism of categories
$$
F: K_\Theta(H_1)\to K_\Theta (H_2).
$$

Denote $\Theta_1=Var(H_1)$ and $\Theta_2=Var(H_2)$.

The isomorphism $F$ one-to-one corresponds
to the isomorphism
$$
 F: K_{\Theta_1}(H_1)\to K_{\Theta_2} (H_2).
$$

which is also denoted by $F$.

Let $W^1=W^1(X)$ and $W^2=W^2(X)$ be the free algebras in
$\Theta_1$ and $\Theta_2$ respectively, $X$ is a finite
set.

\definition{Definition 1}
 An isomorphism $F$ is called correct isomorphism, if

1. $ F(Hom(W^1(X), H_1)=Hom(W^2(Y),H_2)$ for some
$Y$ such that $|Y|=|X|$.

2.If $[s]: A\to Hom(W^1(X),H_1)$ is an identity
embedding, then
$$
F([s]); F(A) \to Hom (W^2(Y),H_2)
$$
is also an identity embedding.
\enddefinition

Categories $K_\Theta(H_1)$ and $K_\Theta(H_2)$
are correctly isomorphic if there exists a correct isomorphism
of categories $K_{\Theta_1}(H_1)$ and $K_{\Theta_2}(H_2)$.

Let us translate this definition to the language of dual
categories. Consider the diagram
$$\CD K_{\Theta_1}(H_1) @>F>>K_{\Theta_2}(H_2) \\ @V F_1 VV @VV F_2 V\\
 C_{\Theta_1}(H_1) @>\Phi>> C_{\Theta_2}(H_2)\endCD$$

Here, $F_1$ and $F_2$ are dual isomorphisms and
$\Phi$ is an isomorphism of categories induced by
$F$.

We consider the correctness condition in terms
of the functor $\Phi$. First, let us connect directly
the functors $\Phi$ and $F$.

Let us take objects $(X,A)$ in $K_{\Theta_1}(H_1)$
and $(Y,B)$ in $K_{\Theta_2}(H_2)$.

Let $T=A'_{H_1}$ and $T^*={B'}_{H_2}$, $A=T'_{H_1}$,
and $B={T^*}'_{H_2}$. Then $F(X,A)=(Y,B)$ if
and only if $\Phi(W^1(X)/T)=W^2(Y)/T^*$.

Consider morphisms. Let $[s]:(X_1,A_1)\to (X_2,A_2)$
be a morphism in $K_{\Theta_1}(H_1)$ and
$$
\bar s: W^1(X_2)/T_2\to W^1(X_1)/T_1
$$
be the corresponding morphism in $C_{\Theta_1}(H_1)$.

In the categories  $K_{\Theta_2}(H_2)$ and
  $C_{\Theta_2}(H_2)$ we have:
$$
[s_1]:(Y_1,B_1)\to(Y_2,B_2)
$$
and
$$
\bar{s_1}: W^2(Y_2)/T^*_2\to W^2(Y_1)/T^*_1.
$$

Now, note that
$$
F[s]:F((X_1,A_1))\to F((X_2,A_2)
$$
is
$$
[s_1]:(Y_1,B_1)\to(Y_2,B_2)
$$
if and only if
$$
\Phi(\bar s): \Phi (W^1(X_2)/T_2)\to
\Phi(W^1(X_1)/T_1)
$$
equals to
$$
\bar{s_1}: W^2(Y_2)/T^*_2\to W^2(Y_1)/T^*_1.
$$

\definition{Definition 2} An isomorphism
$\Phi:C_{\Theta_1}(H_1)\to C_{\Theta_2}(H_2)$
is a correct isomorphism if

1. $\Phi(W^1(X))=W^2(Y)$ for some $Y$, $|Y|=|X|.$

2. If $\mu_X: W^1(X)\to W^1(X)/T_1$ is the
natural homomorphism in $C_{\Theta_1}(H_1)$, then
$$
\Phi(\mu_X): \Phi(W^1(X))\to \Phi(W^1(X)/T_1)
$$
is also a natural homomorphism,
$$
\Phi(\mu_X)=\mu_Y: W^2(Y)\to W^2(Y)/T_2 .
$$
\enddefinition

The main result  here is as follows.
\proclaim{Theorem 1} An isomorphism
$ F: K_{\Theta_1}(H_1)\to K_{\Theta_2}(H_2)$
is a correct isomorphism if and only if
$\Phi: C_{\Theta_1}(H_1)\to C_{\Theta_2}(H_2)$
is a correct isomorphism.
\endproclaim

\subheading{2. Isomorphism, similarity, equivalence}

In the paper \cite{Pl5} the following theorems are proved

\proclaim{Theorem 2} The categories
$K_{\Theta}(H_1)$ and $ K_{\Theta}(H_2)$
is correctly isomorphic if and only if the algebras
$H_1$ and $H_2$ are geometrically similar.
\endproclaim

Let now $H_1$ and $H_2$ be abelian groups, each of
them generates the variety of all abelian groups,
and let $\Theta$ be the variety of all groups.

\proclaim{Theorem 3} Categories $K_{\Theta}(H_1)$
and $K_{\Theta}(H_2)$ are correctly isomorphic if and only
if $H_1$ and $H_2$ are geometrically equivalent.
\endproclaim

Let $H_1$ and $H_2$ be groups such that
$Var(H_1)=Var(H_2)=\Theta=Grp$.

\proclaim{Conjecture} Categories $K_{\Theta}(H_1)$
and $K_{\Theta}(H_2)$ are correctly isomorphic if and only
if the groups $H_1$ and $H_2$ are geometrically equivalent.
\endproclaim

The crucial point here is the problem 3 on $Aut(End(F))$, where
$F$ is free, is stated earlier.{\footnote {Since the problem 3 is
solved positively (E.Formanek \cite{For}), the conjecture  is
true.}}
\bigskip
{\bf Problems}
\bigskip

{\bf Problem 5}. What is the situation in semigroups, i.e.,
$H_1$ and $H_2$ are commutative semigroups.

{\bf Problem 6}. What is the situation in modules over commutative
rings.

{\bf Problem 7}. What is the relation between arbitrary isomorphisms of the
categories $K_\Theta(H_1)$ and $K_\Theta(H_2)$
and correct isomorphisms?. This problem depends on the choice
of $\Theta$ and
$H\in\Theta$.

The following two problems are of the very general nature.

{\bf Problem 8}. Consider the question of equivalence of the
categories $K_\Theta(H_1)$ and $K_\Theta(H_2)$.{\footnote{ This
problem is considered in the paper \cite{Pl11}}

{\bf Problem 9}. When the categories $K_{\Theta_1}$ and
$K_{\Theta_2}$ are isomorphic or equivalent. Consider, in
particular, the case when $\Theta_1$ and $\Theta_2$ are equivalent
categories. Consider also the cases $\Theta_1=\Theta^{G_1}$,
$\Theta_2=\Theta^{G_2}$ for the different algebras $G_1$ and $G_2$
in the given $\Theta$.

Note that in \cite{Pl6} and \cite{Pl8}
 there is an invariant approach (without equations) to the
category of algebraic varieties and
it is quite natural to proceed here from the idea
of equivalence of categories. There are some results in
this direction.

\subheading{3. Perfect and semiperfect varieties of algebras}

3.1. In the first lecture the theorem on isomorphism of the categories
of type $K_P(L)$ was formulated. The key role in this
theorem played the notion of semiisomorphism of $P$-algebras.
All this was done for the variety $Var-P$.

Now we consider the general situation.
We take varieties of the type $\Theta^G=\Theta-G$ and consider the
question
on isomorphisms of the categories of the type $K_{\Theta^G}(H)$
where $H$ is a faithful $G$-algebra. We are going to apply
Theorem 2 in the $G$-algebras case.

Fix a variety $\Theta$ and $G\in\Theta$. Consider the variety
$\Theta^G$ and the category $(\Theta^G)^0$ of free
in $\Theta^G$ algebras of the form $W(X)=G\ast W_0(X)$ with different
finite $X$.

We are interested in automorphisms of this category, and want to
find out when all of them are inner or semi-inner.

We assume that the condition $(\ast)$ is fulfilled. This means that
the algebra $G$ generates the variety $\Theta^G$.
Results on automorphisms of the category $(\Theta^G)^0$ we
apply to the problem of isomorphism of categories
of algebraic varieties with faithful $G$-algebras $H$.

In the given conditions the category $(\Theta^G)^0$
is dual to the category of affine spaces $K^0_{\Theta^G}(G)$,
while the last one is connected with the category of
polynomial maps $Pol-G$.

Let $\varphi:(\Theta^G)^0\to (\Theta^G)^0$
be an automorphism of the category $(\Theta^G)^0$.
It corresponds the automorphism $\tau$ of the
category $K^0_{\Theta^G}(G)$. The automorphisms $\varphi$
and $\tau$ are connected by the following rules.

First of all
$$
\tau(Hom(W(X),G)=Hom(\varphi(W(X)),G).
$$

Let further, $s: W(Y)\to W(X)$
be a morphism in $(\Theta^G)^0$. In the category
$K^0_{\Theta^G}(G)$ it corresponds
$$
\tilde s: (Hom(W(X),G)\to Hom(W(Y),G)).
$$
For every $\nu: W(X)\to G$ we have $\tilde s (\nu)=\nu s $.
Then
$$
\tau(\tilde s)(\nu)=\nu\varphi(s)= \tilde{\varphi(s)}(\nu)
$$
for every point $\nu:\varphi(W(X))\to G$,
$\tau(\tilde s)=\tilde{\varphi(s)}$. This gives one-to-one
correspondence between $\varphi$ and $\tau$.

3.2. We would like to explain  that every automorphism $\tau$ of
the category of affine spaces is, in a sense, a quasi-inner
automorphism. Denote the category of affine spaces by
$K^0_{\Theta^G} (G) = K^0_{\Theta^ G}$.

Recall that every point $\nu\colon W \to G$ satisfies the commutative diagram
$$
\CD
G @>i_G>> W=G*W_0\\
@.\!\!\!\!@ /SE/id_G// \!@VV\nu V\\
@.\!\!\!\!\!\!\!\!\!
G\\
\endCD
$$
\newline

Let $\nu\colon W \to G$ be a point.
Consider the homomorphism defined by
$$
\CD
W @>\nu >>G @>i_G>> W
\endCD
$$
Check that $i_G \nu \colon W \to W$ is an endomorphism of the algebra $W$ in the variety
 $\Theta^G$.
We have to check that there is the commutative diagram
$$
\CD
G @>i_G>> W\\
@. @/SE/i_G// @VVi_G\nu V\\
@. W\\
\endCD
$$
We have
$$
(i_G \nu) i_G = i_G (\nu i_G) = i_G id_G = i_G.
$$
Denote $i_G \cdot \nu = \bar \nu$. For every $w \in W$ the
element $\bar \nu(w) = i_G (\nu(w))$ is a constant in $W$
 and, therefore, $\bar \nu$ is called a constant endomorphism.
Every endomorphism $s\colon W \to W$ leaves constants and,
therefore, $(s\bar \nu) (w) = \bar \nu (w)$, $s \bar \nu = \bar \nu$.
Endomorphism $\bar \nu$ defines the map
$$
\mathop{\nu}\limits^{\simeq} \colon Hom (W, G) \to Hom (W, G).
$$

Note that for every $\nu_0 \colon W \to G$ we have
$\mathop \nu\limits^{\simeq}
(\nu_0) = \nu$.
Indeed, $\mathop \nu\limits^{\simeq} (\nu_0) = \nu_0 \cdot \bar \nu =
 \nu_0 (i_G \nu) = (\nu_0 i_G) \nu = id_G \cdot \nu = \nu$.

Thus, the map $\snu$ takes an arbitrary $\nu_0$ to one and
 the same element $\nu$,
 and, therefore $\snu$ is a constant map.

Consider an arbitrary $\sigma \colon W \to W$, such that $s \sigma = \sigma$
for every $s \colon W \to W$.
Since one can take for $s$ an endomorphism taking $w$ to a constant, $\sigma$
 takes any $w$ to a constant.

It can be shown that to every such $\sigma$ one-to-one corresponds
$\nu: W\to G$ such that
$\sigma = \bar\nu$.

Let now $\vp \colon (\Theta^G)^0  \to (\Theta^G)^0$ be an automorphism of the
category $(\Theta^G)^0$ and let $\tau:K^0_{\Theta^G} \to K^0_{\Theta^G}$
be the corresponding automorphism in the category of affine spaces.

In particular, the constant
$$
\snu\colon Hom (W, G) \to Hom (W, G)
$$
is characterized by the condition $s\bar \nu = \bar \nu$ for every $s$.
This condition can be rewritten as $\snu \tilde s = \snu$.
Apply $\tau$ to the equality $\snu \tilde s = \snu$.
We have
$$
(\snu \tilde s)^\tau = {\snu}^\tau \cdot \tilde s^\tau = {\snu}^\tau \colon
 Hom (W^1, G) \to Hom (W^1, G)
$$
Since $\tilde s^\tau$ is an arbitrary, ${\snu}^\tau$ is a constant, i.e.,
 ${\snu}^\tau =
{\snu}_1$, where $\nu_1\colon \varphi(W)=W^1 \to G$ is a point.
Denote $\nu_1 = \mu(\nu)$.
The map
$$
\mu=\mu_W\colon Hom (W, G) \to Hom(W^1, G)=Hom(\varphi(W),G)
$$
is a bijection.
For every point  $\nu\colon W \to G$ there is
${\snu}^\tau = \widetilde{\overline{\mu_W(\nu)}}$.

\definition{Definition 3}  An automorphism $\tau$  of the category
$K^0_{\Theta(G)}$ is called quasi-inner,  if for an arbitrary
$\tilde s\colon Hom (W^2, G)\to Hom (W^1, G)$, the formula
$$
\tilde s^\tau = \mu_{W^1}\tilde s \mu^{-1}_{W^2}.
$$
takes place.
\enddefinition

\proclaim{Theorem 4} \ \hbox{\rm (see [Be2] for $Var-P$)}
\ Every automorphism $\tau$
 of the category $K^0_{\Theta^G}$ is quasi-inner, i.e., for every
$s: W(Y)\to W(X)$ we have
$$
\tilde s^\tau = \mu_{W^1}\tilde s \mu^{-1}_{W^2}: Hom( \varphi(W(X), G)
\to Hom (\varphi (W(Y), G).
$$
\endproclaim

Note that the definitions of inner and semi-inner automorphism in
the category $(\Theta^G))^0$ is well correlated with the
definition of quasi-inner automorphism in the category
$K^0_{\Theta^G}$. Theorem 4 plays a key role in the proofs of the
theorems on conditions on an automorphism $\varphi$ to be inner or
semi-inner.

\bigskip
3.3. Let us define now a substitutional automorphism of the category
$(\Theta^G))^0$. Let, first, $\varphi$ be a substitution on the objects
$W=W(X)$ of the category $(\Theta^G))^0$. Suppose that if $\varphi(W(X))=W(Y)$
then $|X|=|Y|$.
Let $X = \{ x_1, \cdots, x_n\}$, $Y = \{ y_1, \cdots, y_n\}$.
Define isomorphism $s_W\colon W(X) \to \vp (W(Y))$ by
 $s_W(x_i) = y_i, i = 1, \cdots, n$.
If, further, $\nu\colon W^1 \to W^2$ is a morphism, then set
$$
\vp (\nu) = s_{W^2} \nu s^{-1}_{W^1} \colon \vp (W^1) \to \vp (W^2).
$$
From the substitution $\vp$ we come to the automorphism $\bar \vp$.
We call $\bar \vp$ substitutional automorphism.
\proclaim{Proposition 1} Every automorphism $\varphi$ in the category
$(\Theta^G))^0$ can be presented in the form
$$
\varphi=\varphi^2\varphi^1
$$
where $\varphi^2$ is a substitutional automorphism depending
on $\varphi$ and $\varphi^1=(\varphi^2)^{-1}\varphi$ does not change objects.
\endproclaim

Similar decomposition can be considered for any variety, not necessarily of
the form $\Theta^G$. Such a decomposition of an automorphism of the category
$\Theta^0$ implies, under some conditions, the decomposition of the
corresponding relation of similarity for algebras in $\Theta$.

3.4. Let us pass to the category of polynomial maps $Pol-G $ and
let $s: W(Y)\to W(X)$, where $X = \{ x_1, \cdots, x_n\}$, $Y = \{ y_1, \cdots, y_n\}$,  be given.
Consider the commutative diagram (Lecture 1)
$$
\CD
Hom(W(X), G) @>\tilde s >> Hom(W(Y),G)\\
@V\a_X VV @VV\a_Y V\\
G^{(n)} @>s^\a >> G^{(m)}\\
\endCD
$$
We have $s^\a=\a_Y \tilde s \a^{-1}_X, \; \; \tilde s = \a^{-1}_Y s^\a \a_X$
(Lecture 1).
Let $\vp$ be an automorphism in ${\Theta^G}^0$ and $\vp$ does not
change objects.
The corresponding automorphism $\tau $ also acts identically on objects.
It corresponds an automorphism  $\tau^\a$ of the category $Pol-G$
  which also preserves the objects.
Let $s^\a\colon G^{(n)} \to G^{(m)}$ be  given.
 We set  $\tau^\a (s^\a) = \a_Y \tau(\tilde s) \a^{-1}_X$, where $\tilde s =
\a_Y^{-1} s^\a \a_X$. Here $\tau^\a$ is an automorphism of the
category $Pol-G$. To automorphism $\tau$ corresponds a function
$\mu$, defining $\tau$ as a quasi-inner automorphism. Here
$$
\mu_X = \mu_{W(X)} = \mu_W \colon Hom (W, G) \to Hom (W, G).
$$
This is a bijection and $W=\vp(W)$.
 It corresponds a  bijection
$$
\mu_n\colon G^{(n)} \to G^{(n)}.
$$
Take $a=(a_1, \dots, a_n) \in G^{(n)}$.
Let
$$
\mu_n (a)= \a_X (\mu_X (\a_X^{-1} (a)) = \a_X \mu_X \a^{-1}_X (a).
$$
Then,
$$
\mu_n = \a_X \mu_X\a^{-1}_X, \; \mu_X = \a^{-1}_X \mu_n \a_X.
$$
For the given homomorphism $s\colon W(Y) \to W(X)$,
 $|X|=n, |Y|=m$, we have
$$\tilde s\colon Hom(W(X),G)\to Hom(W(Y),G)$$
 and
$$
\tilde s^\tau = \mu_Y \tilde s \mu^{-1}_X = \tilde s_1\colon Hom (W(X), G) \to Hom
(W(Y), G).
$$
Here
$\mu_Y \tilde s = \tilde s_1 \mu_X$.
In the category $Pol-G$ this equality has the form
$$
\mu_m s^\a = s^\a_1 \mu_n.
$$
3.5. Let us make some general remarks on the transition $W\to \varphi(W)$
in arbitrary $\Theta$. Here, $\varphi$ can change objects.
Consider $s\colon W(X) \to W(X)$, $X = \{ x_1, \dots, x_n\}$
and present it as $s=(s_1, \dots, s_n)$,
where all $s_i$, $i = 1, \dots, n$ are morphisms
$W(x) \to W(X)$.
Here, $s_i$ are defined by the condition
$$
s_i(x) = s(x_i) = w_i (x_1, \dots, x_n) = w_i.
$$
The presentation $s = (s_1, \dots, s_n)$ depends on the basis $X$.
We write $s = (w_1, \dots, w_n)$.
Consider an automorphism $\vp\colon \Theta^0 \to \Theta^0$.
What can be said about the equality
$$
\vp (s) = (\vp(s_1), \dots, \vp(s_n))?
$$
We will see that application of $\vp $ preserves the corresponding
presentation,
 but this is a presentation in some special base, connected with $\vp$.

Consider a system of injections $(\e_1, \dots, \e_n)$,
$$
\e_i\colon W(x) \to W(X).
$$
\definition{Definition 4}  We say that $(\e_1, \dots, \e_n)$ freely defines
an algebra $W$,
if for arbitrary morphisms $f_1, \dots, f_n$, $f_i\colon W(x) \to W(X)$,
there exist a unique $s
\colon W(X) \to W(X)$, such that $f_i = s\e_i$, $i = 1, 2, \dots, n$.
\enddefinition

\proclaim{Proposition 2} A collection $(\e_1, \dots, \e_n)$ freely defines an algebra
$W$ if and only if the elements $\e_1(x), \dots, \e_n(x)$ freely generate $W$.
\endproclaim
Consider, further, automorphism $\vp$ of the category $\Theta^0$ with the
condition $\vp (W(x)) = W(y)$.
\proclaim{Proposition 3} Let the set of morphisms
 $(\e_1, \dots, \e_n), \; \, \e_i
\colon W(x) \to W(X)$ freely define $W=W(X), X = \{ x_1, \dots, x_n\}$.
Then the set $(\vp(\e_1), \dots, \vp(\e_n))$, $$\vp (\e_i) \colon
\vp(W(x)) = W(y) \to \vp(W(X)) = W(Y)$$
freely defines $W(Y)$.
\endproclaim
\definition{Definition 5} A variety $\Theta$ is called a regular variety, if for any free algebra\break
 $W = W(X), |X| = n$, every other system of free generators of $W$, also consists of $n$
 elements.
\enddefinition

Now we can state that if $W(Y) = \vp (W(X))$, then $|Y| = |X|$ if $\Theta $ is regular.
Fix $(\e_1, \dots, \e_n), \e_i (x) = x_i$, take $(\vp(\e_1), \dots, \vp(\e_n))$
 and $y_1\pr = \vp (\e_1)(y), \dots, y\pr_n = \vp (\e_n)(y), Y\pr = \{ y\pr_1, \dots,
y\pr_n\}$.
\proclaim{Proposition 4} If an endomorphism $s: W\to W$ in
the basis $X = \{ x_1, \dots,
 x_n\}$ has presentation $s = (s_1, \dots, s_n)$,
then in the base $Y\pr$ we have the presentation
 $\vp (s) = (\vp(s_1), \dots, \vp(s_n))$.
\endproclaim
A transition from the base $X$ to another base
determines an automorphism used in the proof of the Theorem 5.

3.6. Now we will formulate the Theorem 5.
\definition {Definition 6} A variety $\Theta^G$
is called perfect if every automorphism $\varphi: (\Theta^G)^0 \to
(\Theta^G)^0 $ is inner. If every $\varphi$ is semi-inner, then
$\Theta^G$ is called semiperfect.
\enddefinition

We consider conditions which provide $\Theta^G$ to be perfect or
semiperfect.

To every decomposition $\varphi=\varphi^2\varphi^1$
from the proposition 1 corresponds a decomposition
$\tau=\tau^2\tau^1$. The function $\mu$ corresponds
to $\tau$ and we have a pair $(\mu,\varphi)$. There are
the corresponding pairs for $\tau^2$ and $\tau^1$. Then
$(\mu,\varphi)=(\mu^2\mu^1, \varphi^2\varphi^1)$ and
$(\mu^2\mu^1)_W=\mu^2_{\varphi^1(W)}\mu^1_W$.
This decomposition rule is used in the proofs of the theorems.

Let an automorphism $\varphi$
does not change the objects. For every $W=W(X)$ the automorphism $\varphi$
induces an
automorphism of the semigroup $End W$.

Take $W=W(x)$, $X=\{x\}$ and let $\varphi_0$
be the corresponding induced automorphism of the semigroup $End(W)$.
\bigskip
\definition{Definition 7} An automorphism $\psi$ of the semigroup
$End W$ is called semi-inner if $\psi$ is connected with the
diagram
$$
\CD
G @>i_G>> W(x) = G*W_0(x) \\
@V\sigma VV @VVsV\\
G @>i_G>> W(x) = G*W_0(x)\\
\endCD
$$
in such a way that for every endomorphism
$\nu: W(x)\to W(x)$, the equality $\varphi(\nu)=s\nu s^{-1}$ holds.
\enddefinition

Here $(\sigma, s)$ is a semiautomorphism of the
algebra $W=G\ast W_0$. The same $\psi$ is called inner if
take $\sigma = 1$.

\definition{Definition 8}  A semigroup $End W$ is called
perfect if every its automorphism is inner. A semigroup $End W$ is
called semiperfect if every its automorphism is semi-inner.
\enddefinition

Let now $\psi$ be an automorphism of the semigroup  $End W$ where
$W=W(x)$. It corresponds an automorphism $\bar \psi:
(\Theta^G)^0\to (\Theta^G)^0$. This $\bar \psi$ is constructed in
such a way that if   $\psi$ inner or semi-inner, then
 $\bar \psi$ is also inner or semi-inner.

Return to the initial $\varphi: (\Theta^G)^0 \to (\Theta^G)^0 $
and to the corresponding $\varphi_0$. Then $\bar \varphi_0$
as well as $\varphi$ preserves the objects. This $\bar \varphi_0$
also induces $\varphi_0$ in $End(W)$.

Let $\vp = \vp_1 \bar \vp_0$, $\vp_1 = \vp \cdot \bar\vp_0^{-1}$,
$\varphi_1$ preserves the objects and, besides, does not change constant
endomorphisms of the algebra $W(x)$.
Decomposition of $\vp$ gives rise to decomposition of $\tau$, $\tau = \tau_1 \tau_2$,
 where $\tau_2$ corresponds to the automorphism $\bar \vp_0$.
If now $\mu$ is a function for $\tau$, then $\mu = \mu^1 \mu^2$,
$\mu_W =
\mu^1_{W} \cdot \mu^2_{W}$. Let now $W_0=W(x)$. Then $\mu_{W_0} =
\mu^1_{W_0} \cdot \mu^2_{W_0}$. But $\mu_{W_0}= \mu^2_{W_0}$ since $\vp$ and $\bar \vp_0$, $\tau_1$ and $\tau_2$
 coincide on $W_0$.
Therefore $\mu^1_{W_0} = 1$. This precisely means that $\varphi_1$
 does not change constant
endomorphisms of the algebra $W(x)$.
Applying this fact and Proposition 4 it can be proved that
$\varphi_1$ is an inner automorphism.

The main theorem here is the following

\proclaim{Theorem 5 } If the semigroup $End W(x)$
is perfect in $\Theta^G$,
then the variety $\Theta^G$ is perfect too.
If the semigroup $End W(x)$ is semiperfect then the variety $\Theta^G$
is semiperfect.
\endproclaim

The proof of this theorem uses various reductions, based on transitions from
$\varphi$ to $\tau$, from $\tau$ to $\tau^\alpha$, and to the
corresponding functions $\mu$ and bijections of the form $\mu_X$ and
$\mu_n$.

Using Theorem 2 and the lecture 3, we have also
\proclaim{Theorem 5$'$} Let $H_1$ and $H_2$ be faithful
$G$-algebras. Consider the categories $K_{\Theta^G}(H_1)$ and
$K_{\Theta^G}(H_2)$. Then,

1. If the semigroup $End W(x)$ is semiperfect then the
categories are isomorphic if and only if $H_1$  and $H_2$ are geometrically
equivalent up to a semiisomorphism.

2. If the semigroup $End W(x)$ is perfect then the isomorphism of categories
coincides with geometrical equivalence of  $H_1$  and $H_2$.
\endproclaim
Automorphisms of categories of free algebras of varieties are
studied in the paper \cite{MPP}.


\subheading{ 4. Other results}

A.Berzins proved \cite{Be2} that \pmf 1. If $P$ is an infinite
field, and $P[x]$ is algebra of polynomials with one variable $x$,
then the semigroup $End P[x]$ is semiperfect, i.e. every its
automorphism is semi-inner.

\pmf
2. Let $F$ be a free non-commutative group.
The semigroup $End (F*\{ x\})$ is semiperfect.

So we have
\proclaim{Theorem 6} The variety $Var-P$ is semiperfect.
If the field $P$ does not have automorphisms, then $Var-P$ is perfect.
\endproclaim
\proclaim{Theorem 7}
The variety $Grp-F$ is semiperfect.
\endproclaim

And, further

1. If $L_1$ and $L_2$ are two extensions of the field $P$ then the
categories $K_P(L_1)$ and $K_P(L_2)$ are correctly isomorphic if and only
if $L_1$ and $L_2$ are geometrically equivalent up to a semiisomorphism.

2. If $H_1$ and $H_2$ are two faithful $F$-groups then the categories
$K_{Grp-F}(H_1)$ and  $K_{Grp-F}(H_2)$ are isomorphic if and only if
$H_1$ and $H_2$ are geometrically equivalent up to a semiisomorphism.
It can be shown that in this case $H_1$ and $H_2$ are, in fact, equivalent.

\subheading{5.Problems}

{\bf Problem 10}. What is the situation for Lie $F$-algebras,
where $F$ is a free Lie algebra.

{\bf Problem 11}. What is the situation for associative
$F$-algebras, where $F$ is a free associative algebra or an
infinite dimensional over its center skew field \cite{Pl11}.

It should be noted that Theorem 5 cannot be applied to associative
algebras over a field, since a field $P$ does not
generate the whole variety
of associative algebras over $P$ and the condition
$(\ast)$ does not fulfill.

This lecture concludes the part devoted to equational algebraic geometry.
We present now some general view on the situation in this part.

First of all note that there are problems
which relate to the universal theory. Some of them have been mentioned.
However, the principal thing is to consider situations in the various
special $\Theta$ and special $H\in\Theta$.

The algebraic geometry in groups is on rise now, see \cite{BMR},
\cite{KhM} \cite{MR1}, \cite{MR2}, \cite{BMRo}, \cite{Se} and
others. It is quite reasonable to expect the similar breakthrough
in Lie algebras and semigroups.

For the case of associative algebras over a field or over on some other algebra
of constants it is necessary to clarify how all this is connected
with the theory which is used to call non-commutative algebraic geometry.
In particular, it would be quite reasonable to compare the notions
of noetherian variety of algebras and geometrically noetherian algebras with
the notion of noetherian scheme in the non-commutative algebraic geometry.

One has to distinguish also the cases of noetherian and non-noetherian
non-commutative geometry. Algebraic set $A\subset Hom(W,H)$ is called
{\it correct} if for every system of equations $T$ in $W$ such that
$T'=A$ there exists a finite $T_0\subset T$, such that $T'_0=A$. In the
opposite case the set $A$ is called non-correct. A set $A$ is called
almost correct if $A$ is non-correct but all its proper algebraic subsets
are correct. If an algebra $H$ is not geometrically noetherian, then
there exist non-correct algebraic sets over $H$. Each non-correct
algebraic set contains an almost correct subset. This follows from the
following observation. Take in the lattice $Alv_H(W)$ an infinite
descending system $A_\alpha$, $\alpha\in I$, consisting of non-correct
algebraic sets. Denote $A=\cap_{\alpha\in I} A_\alpha$. Then
$A$ is also non-correct algebraic set. Indeed, take $T_\alpha=A'_\alpha$,
 $T'_\alpha=A_\alpha$, and let $T$ be $\cup T_\alpha$. Then $T'=A$ and
$A'=T''$. Suppose that the set $A$ is correct and $A=T'_0$ where
$T_0$ is a finite subset in $T$. Since $T_0$ is a finite set, $T_0$
is contained in some $T_\alpha$. Therefore, $T''_0=T'' \subset T_\alpha$.
This is impossible. Therefore, the set $A$ is non-correct.

There is the following general problem. Which almost correct sets
$A$ arise for the given non-geometrically noetherian algebra $H$?
This question , first of all, relates to associative
non-commutative algebras over a field.

The next natural object is modules over rings.

Note that there are some results
for algebraic geometry for group representations \cite{KP}.

There are many problems associated with the solution of equations
in specific groups, say in $GL_2(p)$, and with the investigation of the
corresponding algebraic sets.

Note one more problem.

Let $K$ be a ring with the unity and $Mod-K$ be the category of $K$-modules
which is considered as a variety.

The notion of semiisomorphism makes sense
also for this category.

{\bf Problem 12}. For which $K$ can be stated that all
automorphisms of the category of free $K$-modules are semi-inner?
{\footnote{This is true if a ring $K$ is left noetherian \cite{MPP}}}.

Is the theory similar to constructed one, possible in this case?
 The solution
of this problem is connected with the solution of the problem of similarity
for $K$-modules.

\newpage

\baselineskip 12pt
 \topmatter \rightheadtext{Lecture 6}
\leftheadtext{B. Plotkin}
\title Lecture 6\\
\quad\\
Algebraic Geometry in First Order Logic\endtitle
\endtopmatter

\bigskip
\centerline'{\smc Contents}
\bigskip
\roster\item "1." {\smc Introduction}
\newline
\item"2." {\smc Algebraic logic}
\newline
\item"3." {\smc Elementary (algebraic) sets}
\newline
\item"4." {\smc Galois theory in logic}
\newline
\item"5." {\smc Lattices of elementary sets}
\newline
\item"6." {\smc Zariski topology}
\newline
\item"7." {\smc Geometrical properties}
\newline
\endroster

\newpage

\subheading{1. Introduction}

\subheading\nofrills{1.1. The initial idea}

Equational algebraic geometry is a geometry whose algebraic sets are
determined by the
systems of equations of a special type: $w\equiv w'.$
These are equalities in logic.

Here we proceed from arbitrary formulas of elementary $\Th$-logic.
A formula $u$ is considered as an equation, and systems $T$ of such
generalized equations
determine generalized algebraic sets.
We call them {\it elementary sets}.
The point of view on the Zariski topology, which is a main topology in such
a geometry, is
correspondingly changed.

For this new algebraic geometry we need a special category of formulas
which takes the role that the
category $\Th^0$ plays in the equational theory.

This category of formulas assumes the transition from logic to
algebraic logic.
The logic is built accordingly to some variety of algebras $\Th.$

\subheading\nofrills{1.2. Algebra}

Fix a variety $\Th.$
Now, keeping in mind applications from Lecture 7, we proceed from the
situation when algebras
are multy-sorted (not necessarily one-sorted).
Fix for $\Th$ a set of sorts $\G$ which is now finite, but in general it
may be infinite.
Every algebra $G\in \Th$ is recorded as $G=(G_i,i\in \G).$
Operation in the signature $\Om$ is a $\G$-sorted one.
For every $\om\in\Om$ we have its type
$\tau=\tau(\om)=(i_1,\dots,i_n;j),$ \ $i,j\in\G.$
An operation $\om$ of the type $\tau$ is a mapping
$\om: G_{i_1}\times\dots\times G_{i_n}\to G_j.$
All operations of the signature $\Omega $ satisfy some set of
identities.  This fixes the variety $\Theta$ of $\G$-sorted
$\Omega$-algebras.  Let us switch to homomorphisms in $\Theta$ and
to free algebras.  A homomorphism of algebras in $\Theta$ has the
form
$$
\mu=(\mu_i, i \in \G) \colon G = (G_i, i \in \G) \to G'= (G'_i, i
\in \G).
$$
Here $\mu_i \colon G_i \to G'_i$ are mappings of sets, coordinated
with operations in $\Omega$. A congruence $Ker \mu = (\Ker {\mu_
i}, i \in \G)$ is the kernel of a homomorphism $\mu$.

We consider multisorted sets $X = (X_i, i \in \G)$ and the
corresponding free in $\Theta$ algebras
$$
W = W(X) = (W_i, i \in \G).
$$

A set $X$ and a free algebra $W$ can be presented as a free union
of all $X_i$ and all $W_i$, respectively.

Every (multisorted) mapping $\mu: X \to G$ is extended up to a
homomorphism $\mu: W \to G$.  Denote the set of all such $\mu$ by
$\Hom (W, G)$.  If all $X_i$ are finite, we treat this set as an
affine space. Homomorphisms $\mu\colon W \to G$ are points of this
space.

For the given $G = (G_i, i \in \G)$ and $X=(X_i, i \in \G)$ we can
consider the set
$$
G^X = (G^{X_i}_i, i \in \G).
$$
It is the set of mappings
$$
\mu = (\mu_i, i \in \G)\colon X \to G.
$$
Then we have the natural bijection $\Hom(W,G) \to G^X$.
More information about multisorted algebras can be found in [Pl1].

Now let us pass to the models.  Fix some set of symbols of
relations $\Phi$.  Every $\vp \in \Phi$ has its type $\tau = \tau
(\vp) = (i_1, \ldots, i_n)$.  A relation, corresponding to $\vp$,
is a subset in the Cartesian product $G_{i_1} \times\ldots\times
G_{i_n}$.
Now, $\Phi\Theta$ denotes the class of all
 models $(G, \Phi, f)$,
where $G \in \Theta$, and $f$ is a realization of the set $\Phi$
in $G$.  As for homomorphisms of models, they are
homomorphisms of the corresponding algebras which are
coordinated with relations.

\subheading{1.3 \ Logic}

We consider logic in the given variety $\Theta$.  For every finite
$X$ it is determined a logical signature
$$
L = L_X=\{\vee,\wedge,\lnot,\exists x, \; \; x \in X\},
$$
where $X$ is $\mathop\bigcup\limits_{i\in\G} X_i$ for a finite
$\G$.
We consider a set (more precisely, an $L$-algebra) of formulas
$L\Phi W$ over the free algebra $W=W(X)$.  This algebra is
an $L$-algebra of formulas of FOL over the given $\Theta$ and
$\Phi$ and for the given $X$.

First we define the atomic formulas.
They are equalities of the form $w\equiv w'$ with $w, w'\in W$ of
the same sort and the formulas $\vp (w_1, \ldots, w_n)$, where
$w_i \in W$ and all $w_i$ are positioned according to the type
$\tau=\tau(\vp)$ of the relations $\vp$ and to the sorts.  The set
of all atomic formulas we denote by $M = M_X$.  Define $L\Phi W$ as
the absolutely free $L_X$-algebra over $M_X$.

Let us consider another example of an $L_X$-algebra.

Given $W=W(X)$ and $G \in \Theta$, as before, we denote
 by $\Bool (W,G)$ the
Boolean algebra $\Sub (\Hom(W,G))$ of all subsets in $\Hom(W,G)$.
We define also the action of quantifiers in $\Bool (W,G)$.  Let $A$
be a subset in $\Hom(W,G)$ and $x \in X_i$ be a variable of the
sort $i$.  Then $\mu\colon W\to G$ belongs to
$\exists x A$ if there exists $\nu \colon W \to G$ in $A$
such that $\mu (y) = \nu (y)$ for every $y \in X$ of the sort $j,
j \neq i$, and for every $y \in X_i$, $y \neq x$.  Thus we get an
$L$-algebra $\Bool (W, G)$.

Now let us define a mapping
$$
\Val^X_f\colon  M_X \to \Bool (W,G),
$$
where $f$ is a model which realizes
the set $\Phi$ in the given $G$. If $w\equiv w'$ is an equality of
the sort $i$, then we set:
$$
\mu: W\to G \in \Val^X_f (w\equiv w') = \Val^X (w\equiv w')
$$
if $\mu_i(w)=\mu_i (w') $ in $G$.  Here the point $\mu$ is a
solution of the equation $w \equiv w'$.  If the formula is of the
form $\vp (w_1, \ldots, w_n)$, then
$$\mu \in \Val^X_f(\vp(w_1,\ldots, w_n))
$$
if $\vp(\mu(w_1),\ldots, \mu(w_n))$ is valid in the model $(G, \Phi,
f)$.  Here $\mu(w_j) = \mu_{i_j} (w_j)$, $i_j$ is the sort of
$w_j$.  The mapping $\Val^X_f$ is uniquely extended up to the
$L$-homomorphism
$$
\Val^X_f \colon L\Phi W\to \Bool (W,G).
$$
Thus, for every formula $u \in L\Phi W$ we have its value $\Val_f (u)$
in the model $(G,\Phi, f)$, which is an element in $\Bool (W,G)$.

Every formula $u \in L\Phi W$ can be viewed as an equation in the
given model.  The point $\mu\colon W\to G$ is the solution of the
``equation" $u$ if $\mu \in \Val_f (u)$.

\subheading{1.4 Geometrical Aspect}

In the $L$-algebra of formulas $L\Phi W$, $W = W(X)$, we consider
its various subsets $T$, i.e., sets of formulas.
We regard $T$ also as a system
of equations.
On the other hand, we consider subsets $A$  in the
affine space $\Hom (W, G)$, i.e., elements of the $L$-algebra
$\Bool (W, G)$.
For each given model $(G, \Phi, f)$  and for these $T$ and $A$ we
establish the following Galois correspondence:
$$
\eqalign{
T^f &= A = \mathop{\bigcap}\limits_{u\in T} \Val_f (u)\cr
A^f&=T=\{u|A\subset \Val_f (u) \}.\cr}
$$
Here $A=T^f$ is a locus of all points satisfying
the system of  equations  $u \in T$.  Every set $A$ of such kind
is said to be an
{\it algebraic set}
 (or closed set, or  elementary set), determined for the given model.

 The set
$A$ can be treated also as a relation between elements of $G$,
derived from equalities and relations of the basic set $\Phi$.
The relation $A=T^f$ belongs to the multisorted set
$$
G^X = \{ G^{X_i}_i, \; \; i \in \Gamma\}.
$$
The set $T$ of the form $T=A^f$ for some $A$ is an $f$-closed set.
 For an arbitrary $T$ we have its closure
 $
 T^{ff}=(T^f)^f$ and for every $A \subset \Hom(W, G)$ we have the
 closure $A^{ff} = (A^f)^f$.

 It is easy to understand that the following rule takes place:

\medskip
 {\it The formula $v$ belongs to the set $T^{ff}$ if and only if the
 formula
 $$(\mathop{\wedge}\limits_{u \in T} u) \to v$$ holds in the model
 $(G, \Phi, f)$.}
\medskip

 If the set $T$ is infinite then the corresponding formula is infinitary.

  Free in $\Theta$ algebras $W(X)$ with finite $X$ are the objects
 of the category, denoted by $\Theta^0$.   Morphisms of this
 category $s\colon W(X) \to W(Y)$ are arbitrary homomorphisms of
 algebras.  The category $\Theta^0$ is a full subcategory in the
 category $\Theta$.

 Basing on  the first order logic in the given $\Theta$,
we intend to build a category which is similar, in a sense,
 to the category  $\Theta^0$ in the equational logic.
  Thus we pass from pure logic to the algebraic logic.
 The sets of the type $T=A^f$ look here more attractive.

 \heading
 2. Algebraic logic
 \endheading

\subheading{2.1 The main idea}

 Algebraic logic deals with algebraic structures, related to
 various logical structures, i.e., with logical calculi.
Boolean algebras
 relate to classical propositional logic, Heyting algebras relate
 to non-classical propositional logic, Tarski cylindric algebras
 and Halmos polyadic algebras relate to FOL.

 Every logical calculus assumes a set of formulas of the calculus,
 axioms of logic and rules of inference.  On this basis the
syntactical  equivalence of formulas, well correlated with their
semantical equivalence, is defined.  The transition from pure logic
to the algebraic logic is based on treating logical formulas up to
a certain equivalence,  i.e., squeezed formulas.  Such transition
leads to various special algebraic structures, in particular to
the structures mentioned above.

 As for logical calculi, usually they are associated with some
 infinite sets of variables.  Denote such a set by $X^0$.  In our
 situation it is a multisorted set $X^0 = (X^0_i, i \in \Gamma)$.
 Keeping in mind algebraic geometry in logic,
 knowledge theory and its geometrical aspect we
 will use a system of all finite subsets $X = (X_i, \; i \in \G)$
 of $X^0$ instead of this infinite universum.  This leads to
 multisorted logic and multisorted algebraic logic.  Every formula
 has a definite type (sort) $X$.  Denote the new set of sorts by
 $\G^0$.  It is a set of all finite subsets of the initial set
 $X^0$.

\subheading{2.2 Halmos Categories}

 Fix some variety of algebras $\Theta$.   This means that a
 finite set of sorts $\G$, a signature $\Omega = \Omega (\Theta)$
 connected with $\G$, and a system of identities $Id(\Theta)$ are
 given.

 Define Halmos categories for the given $\Theta$.

 First, for the given Boolean algebra $B$ we define its
 existential quantifiers [HMT].
 Existential quantifiers are the mappings $\exists\colon B \to B$
 with the conditions:

 1) \ $\exists 0 = 0$,

 2) \ $a < \exists a$,

 3) \ $\exists (a \wedge \exists b) = \exists a \wedge \exists b$,
 $0, a, b \in B$.

 The universal quantifier $\forall \colon B\to B$ is defined
 dually:

 1) \ $\forall 1  = 1$,

 2) \ $a > \forall a$,

 3) \ $\forall (a\vee \forall b) = \forall a \vee \forall b.$



 Let $B$ be a Boolean algebra and $X$ a set.  We say that $B$ is a
 {\it quantifier $X$-algebra} if a quantifier $\exists
 x \colon B \to B$ is defined  for every $x \in X$
 and for every two elements $x, y \in
 X$ the equality $\exists x \exists y = \exists y \exists x$ holds
 true.

 One may consider also {\it quantifier $X$-algebras $B$ with equalities}
over $W(X)$. In such algebras to
 each pair of elements $w, w' \in W(X)$ of the same sort it
 corresponds an element $w \equiv w' \in B$ satisfying the
 conditions

1)  $w\equiv w$ is the unit in $B$

2)  $(w_1\equiv w'_1 \wedge\ldots \wedge w_n \equiv w'_n) < (w_1
\ldots w_n \omega \equiv w'_1\ldots w'_n \omega) $
where $\omega$ is an operation in $\Omega$ and everything is
coordinated with the type of operation.

Now we will give a general definition of the Halmos category for the
given $\Theta$,
which will be followed by examples.

{\it A Halmos category} $H$ for an arbitrary finite $X=(X_i, i \in
\Gamma)$ fixes some quantifier $X$-algebra $H(X)$ with $X$-equalities.
$H(X)$ is an object in $H$.

The morphisms in $H$ correspond to morphisms in the category
$\Theta^0$.
Every morphism $s_*$ in $H$ has the form
$$
s_*\colon H(X)\to H(Y),
$$
where  $s \colon W(X) \to W(Y) $ is a morphism in $\Theta^0$.
We identify $s_*$ and $s$.

We assume that

1)  The transitions $W(X) \to H(X)$ and $s \to s_*$ constitute a
(covariant) functor $\Theta^0 \to H$.

2) Every $s_* \colon H(X) \to H(Y)$ is a Boolean homomorphism.

3)  The coordination with the quantifiers is as follows:

$\qquad $ 3.1)  $s_1 \exists x a = s_2 \exists x a,
\quad a \in H(X)$, if $s_1 y = s_2y$ for every $y
 \in X, \; y \neq
x$.

$\qquad $ 3.2) $s\exists x a = \exists (sx) (sa) $ if $ sx =
y \in Y$ and $y = sx$ not in the support of $sx'$, $x' \in X, \;
x' \neq x$.

4)  The following conditions describe coordination with equalities

$\qquad $ 4.1) $s_* (w\equiv w') = (sw \equiv sw')$ for $s\colon
W(X) \to W(Y)$, $w, w' \in W(X)$ are of the same sort.

$\qquad $ 4.2) $s^x_w a \wedge (w \equiv w') < s^x_{w'} a $ for an
arbitrary $a \in H(X), x \in X, w, w'$ of the same sort with $x$
 in $W(X)$, and $s^x_w\colon W(X) \to W(X)$ is defined by the rule:
$ s^x_w (x) = w, sy = y, y \in X,\; \;  y \neq x$.

 So, the definition of Halmos category is given.

\subheading{2.3 The example $\Hal_\Theta(G)$}

 Fix an algebra $G$ in the variety $\Theta$.  Define the Halmos
 category $\Hal_\Theta (G)$ for the given $G$.  Take a finite set
 $X$ and consider the space $\Hom(W(X), G)$.  We have defined the
 action of quantifiers $\exists x $ for all $ x\in X$ in the
 Boolean algebra $\Bool (W(X), G)$.  The equality $w\equiv w'$ in
$ \Bool
 (W(X),G)$ is defined as a diagonal, coinciding with the set of
 all $\mu: W(X) \to G$ for which $w^\mu = {w'}^{\mu}$ holds true.
 It is easy to check that in this case the algebra $\Bool (W(X), G)$
 turns to be a quantifier $X$-algebra with equalities.  We set
 $$
 \Hal_\Theta (G) (X) = \Bool (W(X), G).
 $$
 Let now $s\colon W(X) \to W(Y)$ be given in $\Theta^0$.  We have:
$$
\tilde s\colon \Hom (W(Y), G) \to \Hom (W(X), G)$$
defined by $\tilde s (\nu) = \nu s $ for any arbitrary $\nu \colon W(Y)
\to G$.

Now, if $A$ is a subset in $\Hom (W(X), G)$, then $\nu \in s_*A =
s A$ if and only if $\tilde s (\nu) = \nu s \in A$.

We have a mapping:
$$
s_* \colon \Bool (W(X), G) \to \Bool (W(Y), G)
$$ which is a Boolean homomorphism.  One can also check that $s_*$
satisfies the conditions 3--4.  Thus,
the Halmos category $\Hal_\Theta (G)$ is defined.

Note that to each $s_*$ it corresponds a conjugate mapping
$$s^*\colon \Bool (W(Y), G) \to \Bool (W(X), G),$$ where the set
$s^*B$ is the $\tilde s$-image of the set $B$ for every $B \subset
\Hom (W(Y), G)$.

Here, $s^*$ is not a boolean homomorphism, but it preserves
sums and zero.

It may be seen {\cite{Pl1}} that such conjugate mapping can be defined in any
Halmos category. Note the obvious relation between the categories
$Hal_\Theta(G)$ and $Bool_\Theta(G)$.

\subheading{2.4 Multisorted Halmos algebras}

Fix some infinite set $X^0 = (X^0_i, i \in \G)$ and let $\G^0$ be
the set of all finite subsets $X=(X_i, i \in \G)$ in $X^0$.  In this
section multisorted algebra means $\G^0$-sorted.  Every such algebra is of
the form $H=(H(X), X \in \G^0)$.

 A few words about the signature of the algebras to be constructed.
First,  the signature includes $L_X$ for every $X$
 together with equalities $(w\equiv w',  w, w'$ of one sort in
 $W(X))$ as nullary operations.  This is the signature in $H(X)$.
 Second, we consider symbols of operations of the type $s\colon W(X) \to
 W(Y)$.  To each such symbol corresponds an unary operation
 $s\colon H(X) \to H(Y)$.  Denote this
 signature of all $L_X$, all equalities, and all  $s\colon W(X) \to
 W(Y)$ by $L_\Theta$.  This is the signature of FOL in
 $\Theta$ in the multisorted variant.

 Consider further the variety of $\Gamma^0$-sorted
 $L_\Theta$-algebras, denoted by $\Hal_\Theta$.   The identities
 of this variety exactly copy  the definition of Halmos category.
 We call algebras from $\Hal_\Theta$ multisorted Halmos algebras.

Every such algebra can be considered as a small Halmos category.

\subheading{2.5 Algebras of formulas}

First consider a multisorted set of atomic formulas $M = (M(X), X
\in \G^0)$, with $M(X) = M_X$ defined as above.  All $w\equiv w'$
are viewed as symbols of nullary operations-equalities.  The set
of symbols of relations $\Phi $ is fixed.

Denote by $H_{\Phi\Theta} = (H_{\Phi \Theta} (X), X \in \G^0)$ an
absolutely free $L_\Theta$-algebra over the set $M$.  This is the
algebra of formulas of pure FOL in the given $\Theta$.

Now denote by $\tilde H_{\Phi\Theta}$ the result of factorization
of the algebra $H_{\Phi \Theta}$  by the identities of the variety
$\Hal_\Theta$.  It is a free Halmos algebra over the set of atomic
formulas $M$.

Let us introduce the following defining relations:
$$
s_* \vp (w_1,\ldots, w_n) = \vp(sw_1,\ldots, sw_n)\tag {*}
$$
for all $s\colon W(X) \to W(Y)$ and all formulas of the type
$\vp(w_1, \ldots, w_n)$ in $M(X)$.

In the sequel the principal role
will play the Halmos algebra $\Hal_\Theta(\Phi) =
\Hal_{\Phi\Theta}$, defined as a factor algebra of the
free algebra $\tilde H_{\Phi\Theta}$ by the relations of the (*)
 type.  Elements of this algebra are defined as squeezed formulas.

 Consider now values of formulas.  First of all consider a mapping
$$
\Val_f = (\Val_f^X, X \in \G^0)\colon \ M \to \Hal_\Theta (G)
$$
for the model $(G, \Phi, f)$. Here the mappings $\Val^X_f\colon
M_X \to \Bool (W(X), G)=\Hal_\Theta(G)(X)$ have been defined.

This mapping is uniquely extended up to the homomorphisms
$$
\eqalign{
&\Val_f\colon H_{\Phi\Theta}\to\Hal_\Theta (G),\cr
&\Val_f\colon \tilde H_{\Phi\Theta}\to\Hal_\Theta(G).\cr}
$$
Note that the relations $(*)$ hold in every algebra
$\Hal_\Theta(G)$ and  this gives a canonical
homomorphism of Halmos algebras $$Val_f :\Hal_\Theta (\Phi) \to
\Hal_\Theta (G).$$  It determines the value of the formulas
$\Val_f(u)$ (pure and squeezed) in the given model $(G, \Phi, f)$.

It is easy to see that the kernel $\Ker(Val_f)$ is precisely the
elementary theory of the model $(G,\Phi,f)$ in the logic
of the variety $\Theta$.

In fact, elementary theory of the algebra $G$ or the model $(G,\Phi,f)$
is considered also on the logic of the variety $\Theta^G$. This logic
is more reach in respect to the given $G$.

We call two pure formulas $u$ and $v$ of the given type $X$
{\it semantically equivalent}, if $\Val_f (u) = \Val_f (v)$ for every
model $(G, \Phi, f)$.

 The following main theorem takes place:

 \proclaim{Theorem 1}  Two formulas $u$ and $v$ are semantically
 equivalent if and only if the corresponding squeezed formulas $\ol u $ and
 $\ol v$ coincide in the algebra $\Hal_\Theta (\Phi)$.
\endproclaim

This theorem explains the role of algebra $\Hal_\Theta(\Phi)$ as a
main structure of the multisorted algebraic logic for FOL in the
given $\Theta$.  The same algebra plays the essential part in the
algebraic geometry in the FOL  in $\Theta$.  Besides that, the
role of the algebras $\Hal_\Theta (G)$ is underlined by the
following theorem:
\proclaim{Theorem 2}  All algebras $\Hal_\Theta (G)$ over
different $G \in \Theta$ generate the variety of Halmos algebras
$\Hal_\Theta$.
\endproclaim

Define the notion of the logical kernel of a homomorphism.

Let the homomorphism $\mu\colon W(X) \to G$ be given.  One can view
its kernel $\Ker \mu$ as a system of all formulas $w \equiv w'$
with $w, w'$ of the same sort in $W(X)$, for which $\mu \in \Val
(w \equiv w')$.

{\it The logical kernel} $\Log\Ker \mu$ naturally generates the standard
$\Ker \mu$.
We set: the formula $u \in \Hal_{\Phi\Theta}(X)$ belongs to
$\Log\Ker (\mu)$ if the point  $\mu$ lies in $\Val_f (u)$, i.e.,
$\mu$ is a solution of the ``equation" $u$ in the given model $(G,
\Phi, f)$.  It is easy to understand, that for every point $\mu$
its logical kernel is an ultrafilter of the Boolean algebra
$\Hal_{\Phi\Theta} (X)$.  It is also clear, that the kernel $\Ker
\mu$ is the set of all equalities in the logical kernel.

\heading
3. Elementary (algebraic) sets
\endheading

\subheading{3.1 Preliminary remarks}

In the sequel we call the sets below algebraic sets
although it would be more sensible to call them elementary
sets.

Algebraic sets are the sets, determined by FOL formulas.
We work with squeezed formulas, i.e., formulas of the algebra
$\Hal_\Theta (\Phi) = \Hal_{\Phi\Theta}$.

For the given place $X$ consider sets of formulas $T$ in
$\Hal_{\Phi\Theta} (X)$ and the sets of points $A$ in the space
$\Hom (W(X), G)$.
We establish  a Galois correspondence, determined by
the given model $(G, \Phi, f)$.  It looks like
$$
\eqalign{
&T^f = A = \bigcap_{u\in T} \Val_f(u) = \{ \mu \vert T\subset
\Log\Ker (\mu)\}\cr
&A^f = T = \{ u \vert A\subset \Val_f (u) \} = \bigcap_{\mu
\in A} \Log\Ker (\mu).\cr}
$$
As above, we call a set $A$ represented as $A=T^f$ an {\it algebraic
set} or {\it elementary set} for the given model $(G, \Phi, f)$.

The set $T$, represented as $A^f=T$, is always a filter of Boolean
algebra $\Hal_{\Phi\Theta} (X)$, since by the definition it is
an intersection of
ultrafilters.  We call it $f$-closed filter. If $A$ is an algebraic
set then $T=A^f$ can be considered also as the elementary theory
of the given $A$. One
can consider here the Boolean algebra $\Hal_{\Phi\Theta} (X)/T$ for this
$T$.  If $T^f=A$ and $A^f=T$, then the algebra $\Hal_{\Phi\Theta}
(X)/T$ is considered as an invariant of the algebraic set $A$.
 This invariant is a {\it coordinate algebra} of the set $A$.  It can be
 represented as an algebra of regular functions determined on
$A$ (see [Pl2]).

 Every algebraic set, defined in Section 1, is also an algebraic
 set according to this new definition.  The opposite is not true,
 because in the new variant additional operations of the type
 $s\colon W(X) \to W(Y)$ are involved in the formulas.

 Consider now the relations between Galois correspondence  and
 morphisms of Halmos categories.

 For every $ s\colon W(X) \to W(Y)$ and every $A$ of the type $X$
 we considered a set $B=s_* A$ of the type $Y$.  If $B$ is of the
 type $Y$, then $A = s^* B $ is of the type $X$.  Define the
 operations $s_*$ and $s^*$ on the sets of formulas.

 If $T$ is a set of formulas in $\Hal_{\Phi\Theta}(Y)$, then
 $s_*T$ is a set of formulas in $\Hal_{\Phi\Theta} (X)$ defined by
 the rule:
 $$
 u \in s_*T \Leftrightarrow s  u \in T.
 $$
 If $T$ is a set of formulas in $\Hal_{\Phi\Theta}(X)$, then $s^*
 T$ is contained in $\Hal_{\Phi\Theta} (Y)$ and it is defined by
 $$
 u\in s^*T \; \; \hbox{\rm if } \; \; u = sv, \; \;  \; v \in T.
 $$
 The following theorem takes place:
 \proclaim{Theorem 3}

 1. If $T$ lies in $\Hal_{\Phi\Theta} (X)$, then
 $$
 (s^* T)^f = s_* T^f = sT^f.
 $$

 2.  If $B \subset \Hom (W(Y), G)$, then
 $$
 (s^*B)^f = s_* B^f.
 $$

 3.  If $A \subset \Hom(W(X), G)$, then $s^*A^f\subset(s_*A)^f$.
\endproclaim

 It follows from these rules that

{\it
 1.  If $A = T^f$ is an algebraic set, then $sA$ is also an
 algebraic set.

 2.  If $T=B^f$ is $f$-closed, then $sT=s_*T$ is $f$-closed.
}

\subheading{3.2. Categories $K_{\Phi\Theta} (f)$ and $C_{\Phi\Theta}(f)$}

 Fix a model $(G,\Phi, f)$ and define a category of algebraic
 sets $K_{\Phi\Theta} (f)$ for this model.  Objects of this
 category have the form $(X, A)$, where $A = T^f$ for some $T$.
 $X$ is the place for both $A$ and $T$.

 Let us now define  morphisms
 $
 (X,A)\to (Y,B).
 $
 Proceed from $s\colon W(Y) \to W(X)$.  We say that $s$ is
 admissible for
$A$ and $B$ if $\tilde s(\nu) = \nu s \in B$ for any
$\nu \in A$.  It is clear that $s$ is admissible for $A$
and $B$ if $A\subset sB$. A mapping $[s]:A\to B$ corresponds to each $s$
admissible for $A$ and $B$.

We consider {\it weak} and {\it exact} categories $K_{\Phi\Theta}(f)$.
In the first
one the morphisms look like $s:(X,A)\to (Y,B)$, while in the second one
like $[s]:(X,A)\to (Y,B)$.

If $s_1$ is admissible for $A$ and $B$
and $s_2$ for $B$ and $C$, then $A \subset s_1 B$,
 $B \subset s_2 C$, $ s_1 B \subset s_1s_2 C$, and $s_1s_2$ is
admissible for $A$ and $C$.

Define now a category $C_{\Phi\Theta} (f)$.  Its objects are
Boolean algebras of the type
\newline
 $\Hal_{\Phi\Theta} (X) / T$, where $T
= A^f$ for some $A$.

Consider morphisms
$$
\Hal_{\Phi
\Theta}(Y) / T_2 \overset \ol s\to{\longrightarrow}
\Hal_{\Phi\Theta} (X) /T_1.
$$
We proceed here from $s\colon W(Y) \to W(X)$ and pass to the new
$s : \Hal_{\Phi\Theta} (Y) \to \Hal_{\Phi\Theta} (X)$.
Assume that $s u
\in T_1$ for every $u \in T_2$.
The homomorphism $s$ is admissible for $T_2$ and $T_1$
in this sense.  Define homomorphisms  $\ol s$ for such $s$.  This
defines morphisms in $C_{\Phi \Theta}(f)$.

The next two propositions determine the correspondence between the
categories $K_{\Phi\Theta}(f)$ and $C_{\Phi\Theta} (f)$.

\proclaim{Proposition 1}  Homomorphism $s\colon W(Y)\to W(X)$ is
admissible for varieties $(X, A)$ and $(Y, B)$ if and only if this
$s$ is admissible for $T_2 = B^f$ and $T_1 = A^f$.
\endproclaim
\proclaim{Proposition 2}   If $s_1, s_2 \colon W(Y) \to W(X)$ are
admissible for $A$ and $B$, then $[s_1]=[s_2]$ implies
$\ol{s_1} = \ol{s_2}$.
\endproclaim

It follows from these two propositions
that the transition
$$
(X,A) \to\Hal_{\Phi\Theta} (X)/A^f
$$
determines a contravariant functor
$$
K_{\Phi\Theta} (f) \to C_{\Phi\Theta} (f)
$$
for weak and exact categories $K_{\Phi\Theta} (f)$.
Duality takes place under some additional conditions.
In particular, this is the case when we proceed from
the variety $\Theta^G$ for the given $G$ in $\Theta$.

\subheading{3.3 Categories $K_{\Phi\Theta} $ and $C_{\Phi\Theta}$}

In the categories $K_{\Phi\Theta}$ and $C_{\Phi\Theta}$ the model
$(G, \Phi, f)$ is not fixed.  Objects of $K_{\Phi\Theta}$ have the
form $(X,A; G,f)$.  Here $f$ is a realization of the set
$\Phi$, fixed for the category $K_{\Phi\Theta},$ in the algebra $G$,
and $A=T^f$ for some $T\subset \Hal_{\Phi\Theta}(X)$.

Define morphisms
$$
(X,A; G_1, f_1) \to (Y,B; G_2, f_2).
$$
They act on all components of the objects.  Proceed from the
commutative diagram
$$
\CD
W(Y) @>s>> W(X) \\
@V\nu' VV @VV\nu V\\
G_2 @<\delta<< G_1\\
\endCD
$$
Consider a pair $(s, \delta)$ and write $(s, \delta) (\nu) = \nu'
= \delta\nu s$.

Let now $A = T_1^{f_1} $ be of the type $X$ and $B=T_2^{f_2}$ of
the type $Y$.  We say that the pair $(s, \delta)$ is {\it admissible}
for $A$ and $B$ if $(s, \delta) (\nu) \in B$ for every $\nu
\in A$.

We need some further auxiliary remarks.
For every $\delta: G_1 \to G_2$ and every $X$ we have a mapping
$$
\tilde\delta: \Hom(W(X), G_1) \to \Hom (W(X), G_2)
$$
defined by the rule
$$
\tilde \delta(\nu) = \delta\nu, \nu \in \Hom (W(X), G_1).
$$

Determine $\delta_* A\subset \Hom (W(X), G_1)$ for every $A\subset
\Hom (W(X), G_2)$, setting
$$
\nu \in \delta_* A \; \; \hbox{\rm if } \; \; \delta\nu=\tilde
\delta(\nu) \in A.
$$
We change notations to $\delta_* A=\delta A$.

If $A \subset \Hom (W(X), G_1)$, then $\delta^* A \subset \Hom
(W(X), G_2)$ and $\nu \in \delta^* A$ if $\nu = \delta\nu_1, \nu_1
\in A$.

Now we can  say that the pair $(s, \delta)$ is admissible
for $A$ and $B$ if $\delta^*A \subset s B$, or,  the same,
$A\subset \delta sB = s\delta B$.

We have morphisms
$$
(s, \delta) : (X, A; G_1, f_1)\to(Y, B; G_2, f_2)$$
and
$$
([s], \delta)\colon (X, A; G_1, f_1) \to (Y, B; G_2, f_2)
$$
for the admissible $(s, \delta)$.   Here $[s]: A \to B$ is a
mapping, induced by the pair $(s, \delta)$.  We get weak and
exact categories $K_{\pt}$.  It can be proven that the pair $(s,
\delta)$ is admissible for $A$ and $B$ if and only if
the homomorphism $s\colon \Hal_{\pt} (Y)\to Hal_{\pt} (X) $
is admissible in respect
to $T_2 = B^{f_2}$ and $T_1 = (\delta^* A)^{f_2}$.  This leads to
a natural definition of the category $C_{\pt}$ with contravariant
functor $K_{\pt} \to C_{\pt}$.

Let us define the categories $K_{\pt}(G)$ and $C_{\pt} (G)$.
Here $G$ is a fixed algebra in $\Theta$, while the realizations
$f$ of the set $\Phi$ in $G$ change.

The objects in $K_{\pt}(G)$ have the form
$$
(X, A; f).
$$
The morphisms
$$
(X, A, f_1) \to (Y, B, f_2)
$$
are defined according to the general definition of the morphisms
in $K_{\pt}$ with identical $\delta = \vare \colon G\to G$.

Objects in $C_{\pt} (G)$ have the form
$$
(\Hal_{\pt} (X) /T, f), \; \; \hbox{\rm where} \; \; T=A^f
$$
for some $A$ of the type $X$.

The transition
$$
(X, A; f) \to (\Hal_{\pt} (X) /A^f, f)
$$
determines the functor $K_{\pt} (G) \to C_{\pt} (G)$.
Here $K_{\pt} (G$) is a subcategory in $K_{\pt}$ and every
$K_{\pt}(f) $ is a subcategory in $K_{\pt} (G)$.  The same for
$C$.

\heading
4. Galois theory in logic.
\endheading

Galois theory we are talking about is tied in with the
considered algebraic geometry and, besides, it will be used
in the next lecture.

\subheading{4.1 Automorphisms}

Let $\delta\colon G\to G$ be an automorphism of the algebra $G$,
$\delta \in Aut G$.  Then for every $X$ we have a substitution
$$
\tilde\delta: \Hom (W(X), G) \to \Hom (W(X), G).
$$
This substitution induces an automorphism of the Boolean algebra
$\Bool (W(X), G)$.  For each element $A$ from this Boolean algebra
 we have $\delta_* A
= \delta A$.

Actually, it is an automorphism of the quantifier $X$-algebra with
equalities.  This leads to the automorphism $\sigma=\delta_*$ of
the Halmos algebra $\Hal_\Theta (G)$. We have a representation
$$
Aut (G)\to Aut (\Hal_\Theta (G)).
$$
\proclaim{Theorem 4}  The given representation is an isomorphism
of the groups of automorphisms $Aut (G)$ and $ Aut (\Hal_\Theta (G))$.

Thus, the group of automorphisms $Aut G$ can be considered as a
group of automorphisms of the algebra $\Hal_\Theta (G)$.
\endproclaim

Let, further, $\delta\colon G_2 \to G_1$ be an isomorphism of
algebras in $\Theta$.  As above, we may take the isomorphism of
Halmos algebras
$$
\delta_*\colon \Hal_\Theta (G_1) \to \Hal_\Theta (G_2).
$$
This isomorphism is well correlated with the homomorphisms of the
type $\Val_f$.

Let the model $(G_2, \Phi, f)$ be given.  The isomorphism
$\delta:G_2\to G_1$ uniquely determines the model $(G_1, \Phi,
f^\delta)$,  isomorphic to the initial model.

We have the commutative diagram

$$
\CD
\; \; \; \; \; \;\;\;\; \Hal_\Theta(\Phi) \\
@/SW/\Val_{f^\delta} //\; \; \quad  @[2]/SE//\Val_f / \\
\!\!\!\!\!\!\!\! \!\!\!\!\!\!\!\!\!\!\!\!\!\!\!\!\! \Hal_\Theta (G_1)@>\delta_* >> \Hal_\Theta (G_2) \\
\endCD
$$
The arrow $\Val_{f^\delta}$ is uniquely determined by the other
two arrows due to the diagram commutativity.
\subheading{4.2 \   The main theorems of Galois theory}

Consider the algebra $\Hal_\Theta (G)$ together with its group of
automorphisms $\Aut (G)$.

Let us define the standard Galois correspondence.

Let $H$ be a subset of $\Aut (G)$.  Then $H' = R$ is a subalgebra
of $\Hal_\Theta(G)$ (with equalities), consisting of all $A$ with
$\delta A = A$ for every $\delta \in H$.

Let $R$ be a subset in $\Hal_\Theta(G)$.  Then $R'=H$ is a
subgroup in $\Aut (G)$ consisting of all $\delta$ with $\delta A =
A $ for every $A \in R$.

Define closures $R^{''}$ and $H^{''}$.
\proclaim{Theorem 5}  If algebra $G$ is finite, then every
subalgebra in $\Hal_\Theta (G)$ and every subgroup in $\Aut (G)$
are closed.  There is one-to-one correspondence between them.
\endproclaim

Let now two finite algebras, $G_1$ and $G_2$, be given.  Take a
subalgebra $R_1$ in $\Hal_\Theta(G_1)$ and a subalgebra $R_2$ in
$\Hal_\Theta(G_2)$.

\proclaim{Theorem 6}  If there is an isomorphism $\gamma\colon
R_1\to R_2$, then there exists  an isomorphism $\delta: G_2 \to G_1$,
 such that $\delta_*: \Hal_\Theta (G_1) \to \Hal_\Theta
(G_2)$ induces the given $\gamma: R_1\to R_2$.

Here,

1. $R_2 = \delta_* (R_1)$,

2. $R'_1 =\delta R'_2 \delta^{-1}$,
\endproclaim

 $R'_1 $ and $R'_2$ determine $R_1$ and $R_2$ with $R_2 =
\delta_* (R_1)$.

Note that the considered Galois theory goes back to the papers
by M.I.Krasner {\cite{KR}} and it would be reasonable to call it
the Galois-Krasner theory.

\subheading{4.3  Automorphisms and algebraic sets}

Denote a group of automorphisms of the model $(G,\Phi,f)$ by $\Aut
(f)$. It is a subgroup of $\Aut (G)$.

\proclaim{Theorem 7}  Every algebraic set $A=T^f$ is invariant
under the action of the group $\Aut (f)$.   If $A= T^f$ and $\delta
\in \Aut (f)$, then $\delta A = A$.
 \endproclaim

This theorem to some extent
determines the group structure of algebraic sets.

Consider $\Val_f\colon \Hal_\Theta (\Phi)\to\Hal_\Theta (G)$ for
the model $(G, \Phi, f)$ and let again $R_f$ be the image of this
homomorphism.  $R_f$ is a subalgebra in $\Hal_\Theta (G)$
consisting
 of algebraic sets of the form $A = \Val_f(u)$.  Here $T$ consists
 of one element $u$.  We call $A$ a simple algebraic set.

\proclaim{Theorem 8}
$R_f'=\Aut(f)$ holds for every model $(G,\Phi,f).$
If the algebra $G$ is finite, then $\Aut(f)'=R_f.$\ep

In particular, it follows that if the algebra $G$ is finite, then every
algebraic set $A=T^f$
is a simple algebraic set, $A=\Val_f(u)$ for some
$u\in\Hal_{\Phi\Th}=\Hal_\Th(\Phi).$

Note in the addition to Theorem 7 that the set $A\subset \Hom(W,G)$  is
elementary
(algebraic) set over a model with the algebra $G,$ if and only if $\dl
A=A,$\ $\dl\ne 1$ holds
for some $\dl\in \Aut(G).$
We can select a set $\Phi$ and its  realization $f$ in $G$ by $\dl$ in such
a way that
$\dl\in\Aut(f)$ and, simultaneously, we have $\dl A=A$ for every
$\dl\in\Aut(f).$

\subheading{5. Lattices of elementary sets}

For every model $(G,\Phi,f)$ and every finite set $X$ denote: $\Alv_f(X)$
is the set of all
elementary sets for the model $(G,\Phi,f)$ in the space $\Hom(W,G);$ \
$\Cl_f(X)$ is the set of
all closed sets $T$ in $\Hal_{\Phi\Th}(X).$

We consider both these sets as lattices as well.
The intersection operation is defined here as set theory intersection,
while the union is
defined by the rules:
$$A\bigtriangledown B=(A\cup B)^{ff};\ A,B\in \Alv_f(X)$$
$$T_1\bigtriangledown T_2=(T_1\cup T_2)^{ff};\ T_1,T_2\in\Cl_f(X).$$
The transitions $A\to A^f=T$ and $T\to T^f=A$ give the duality of lattices.

Further, we can specify various levels of logic, all of which are the
parts of the algebra
$\Hal_\Th(\Phi).$
Levels of algebraic geometry are connected with them.
In particular, equational geometry is connected with equational logic.
The attitude to the pointed lattices changes.
Let us popint out some interesting levels.

$L_0$ -- equational logic.
Its formulas have the form $w\equiv w'.$
They are equalities.

$L_1$ -- pseudoequational logic.
Its formulas are pseudoequalities of the form $w_1\equiv w_1'\vee\dots\vee
w_n=w_n'.$

$L_2$ -- universal logic over equalities.
  The formulas have the form $w_1\equiv
w_1'\vee\dots\vee w_n\equiv w_n'\vee v_1\not\equiv v_1'\vee\dots v_m\not\equiv
v_m'.$

$L_3$ -- we mean here positive logic whose formulas are built without
negations.

$L_4$ -- universal logic (without quantifiers).

$L_5$ -- all that is built from atoms without quantifiers and negations
(only $\vee$ and
$\wedge).$

$L_6$ -- the whole algebra $\Hal_\Th(\Phi).$

Finally, we say that the logic $L$ is $\vee$-closed, if $u,v\in L$ implies
$u\vee v\in L.$
With each of these logics $L$ a definite level $\ell$ of geometry is
associated.
A notion of logical kernel of a homomorphism relates to the level $\ell.$
If the homomorphism $\mu: W(X)\to G$ in $\Th$ is given, then we set
$\ell-\log\ker(\mu)$ is
the set of formulas $u$ of the level $\ell$ in $\Hal_{\Th\Phi}(X)$ for which
$\mu\in\Val_f(u).$
We assume that the algebra $G$ is included in the model $(G,\Phi,f).$

According to these considerations, we localize Galois connections,
all sets of formulas $T$ are on a definite level $\ell$ in a
definite logic $L.$ Let us rewrite  these connections

For every $A\subset \Hom(W,G)$ we have
$$T=A^{f,\ell}=\capl_{\mu\in A}(\ell-\log\ker(\mu))=\{u\bigm| u\in
\ell\qquad \text{and}\qquad A\subset
\Val_f(u)\}.$$
If, further $T$ is of the level $\ell,$ then $$T^f=A=\cap_{u\in
T}\Val_f(u)=\{\mu: W\to
G\bigm| T\subset\ell-\log\ker(\mu)\}.$$

Consider corresponding Galois closures.
For the given $A\subset \Hom(W(X),G)$ we have: $A^{ff,\ell}=(A^{f,\ell})^f.$
For the given $T$ of the level $\ell$ we have $T^{ff,\ell}=(T^f)^{f,\ell}.$
Consider also $\ell,f$-closed $A$ and $T$.
We consider the lattices $\Alv_f(X)$ and $\Cl_f(X)$  on the level $\ell$ as
well.
We denote them $\Alv_{f,\ell}(X)$ and $\Cl_{f,\ell}(X)$, respectively.

\subheading{6. Zariski topology}

We speak of topology in the affine space $\Hom(W,G),W=W(X),$
for the given model
$(G,\Phi,f)$ on the given level of logic $\ell.$
We suppose here to consider only positive formulas.

For the given level $\ell$ a Zariski topology in $\Hom(W,G)$ over the model
$(G,\Phi,f)$ is a
topology generated by all $\ell$-algebraic sets $A$ in $\Hom(W,G)$ as
closed sets.

\proclaim{Theorem 9}
If the level
$\ell$ is a $\vee$-closed logic, then all closed sets in the Zariski
topology are
exactly all $\ell$-algebraic sets.\ep

In this case the lattice $\Alv_{f,\ell}(X)$ is a sublattice of the lattice
$\Bool(W,G).$
It is distributive, as is the lattice $\Cl_{f,\ell}(X).$

Interesting case of Zariski topology is provided by empty set
$\Phi$. In the special situation when $\Theta=\Var$-$R$ and the
field $R$ is the field of real numbers, the case when $\Phi$ is
an one-element  relation $\varphi$ which is the order relation is
of special interest.

It is natural to treat all the considerations of this lecture in the
classical situation of
the variety $\Var$-$P$.

\subheading{7. Geometrical properties}

We work at the beginning in the largest FOL-- $L_6.$

\definition{Definition 1}
The model $(G,\Phi,f)$ is called {\it geometrically noetherian}
if for every finite
set $X$ and every
set of formulas $T$ in $\Hal_{\Phi\Th}(X)$ in $T$ there is some finite part
$T_0$ with
$T_0^f=T^f.$
This means that  $T_0^{ff}=T^{ff}.$\enddefinition

\proclaim{Theorem 10}
The model $(G,\Phi,f)$ is geometrically noetherian if and only if the
minimality condition
holds in the lattice $\Alv_f(X).$
Correspondingly, in the lattice $\Cl_f(X)$ we have the maximality condition.\ep

Now let $T_0=\{u_1,\dots, u_n\}$ and $u$ is $u_1\wedge\dots\wedge u_n$.
Then $T_0^f=\Val_f(u)=\{u\}^f.$
Thus, we may claim that if the model $(G,\Phi,f)$ is geometrically
noetherian, then every
algebraic set over this  model is  a simple algebraic set.

However, the corresponding element $u=u_1\wedge\dots\wedge u_n$ does not necessarily
belong to the
initial set $T.$
We call a model $(G,\Phi,f)$ weak geometrically noetherian if
every
algebraic set over it is a simple algebraic set. Weak noetherian
model is not necessarily geometrically noetherian.

However, we may claim that every finite model is geometrically noetherian
(compare Theorem 8
in Galois theory).
We may also claim that any finite cartesian product of geometrically
noetherian models is also
a geometrically noetherian model.
A submodel of a geometrically noetherian model is also geometrically
noetherian.
 The similar properties not true in general in the transition to cartesian
powers.

The notion of the model to be geometrically noetherian may be also defined
on different
special levels $\ell.$
It is easy to understand that if the model $(G,\Phi,f)$ is noetherian on
the absolute level
$\ell_6,$ then it is noetherian on each other level $\ell.$

We may assume that by lowering the level $\ell,$ one may obtain good
properties.
For example, for which level $\ell$ the property of the model $(G,\Phi,f)$
to be noetherian
implies that each of its cartesian powers is noetherian as well.
It is valid for $L_0$ and for the logic of all atomic formulas.
It is not clear when else this is true.
It is most likely never (see also Theorem 11 below).

Let us pass to the notion of geometrical equivalence of two models.
Let two models $(G_1,\Phi,f_1)$ and $(G_2,\Phi,f_2)$ with the same $\Phi$
be given.

\definition {Definition 2}
The models $(G_1,\Phi,f_1)$ and $(G_2,\Phi,f_2)$ are geometrically
equivalent if
$T^{f_1f_1}=T^{f_2f_2}$ holds for every finite $X$ and every $T$ in
$\Hal_{\Phi\Th}(X).$

If the models $(G_1,\Phi,f_1)$ and $(G_2,\Phi,f_2)$ are geometrically
equivalent, then
\roster\item"1." The lattices $\Cl_{f_1}(X)$ and $\Cl_{f_2}(X)$ coincide,
while the lattices
$\Alv_{f_1}(X)$ and $\Alv_{f_2}(X)$ are isomorphic.
\item"2." The categories $C_{\Phi\Th}(f_1)$ and $C_{\Phi\Th}(f_2)$
coincide, while the
categories $K_{\Phi\Th}(f_1)$ and $K_{\Phi\Th}(f_2)$ are isomorphic.
\item"3." These models are elementary equivalent.\endroster
It follows from the first claim that if the models are geometrically
equivalent and one of
them is geometrically noetherian, then the second one is also geometrically
noetherian.
\enddefinition

\proclaim{Theorem 11}
Let the model $(G,\Phi,f)$ be geometrically noetherian.
Then each of its ultrapower is also geometrically noetherian, and all these
ultrapowers are
geometrically equivalent to the initial $(G,\Phi,f).$\ep

All the described notions are naturally connected with the logic of
generalized (infinitary)
formulas of the kind $(\wedge_{u\in T})\to v,$\ or, what is the same
$T\to v.$
For the geometrically noetherian models it is sufficient to proceed from
the usual finite
formulas.

Note that the notion of geometrical equivalence of models can be
also considered on
different logical levels.
Here, if the models are equivalent on the absolute level $\ell_6, $ then
they are equivalent
on every other level.

Elementary equivalence of the models does not generally imply their
geometrical equivalence.
However, we may claim the following
\proclaim{Theorem 12}
If two models are elementary equivalent and one of them is geometrically
noetherian, then the
second one is aso geometrically noetherian and these models are
geometrically equivalent.\ep

In concern with the notion of geometrically noetherian model let us
return to the notion of the logical kernel of a homomorphism of
the form $\mu: W(X) \to G$. It is easy to see that if $\Log_f\Ker(\mu)$
is such a kernel and $\Val_f(\Log_f\Ker(\mu))$ is its image in the algebra
$R_f(X)$ then this image is a principal ultrafilter if the model
$(G,\Phi, f)$ is geometrically noetherian. This ultrafilter is
generated by the algebraic set ${\{\mu\}}^{ff}$.

\newpage

\baselineskip 12pt
 \topmatter \rightheadtext{Lecture
7} \leftheadtext{B. Plotkin}
\title Lecture 7\\
\quad\\
Databases and Knowledge bases\\
Geometrical aspect\endtitle
\endtopmatter

\bigskip

\bigskip
\centerline{\smc Contents}
\bigskip
\roster\item "1." {\smc Introduction}
\newline
\item"2." {\smc Categories of elementary knowledge}
\newline
\item"3." {\smc Databases and Knowledge bases}
\newline
\item"4." {\smc Equivalence of Databases and Knowledge bases}
\newline
\item"5." {\smc Conclusion}
\newline
\endroster

\newpage

\subheading{1. Introduction}

We use here the material from lecture 6, mainly, in knowledge theory.

Elementary knowledge is considered to be a
 first order knowledge, i.e., knowledge
that can be represented by the means of the First Order Logic
(FOL). The corresponding applied field (field of knowledge) is
based on some variety of algebras $\Theta$, which is arbitrary but
fixed and can be multisorted. This variety $\Theta$ is considered
as a knowledge type, like in database theory, databases of a data
type $\Theta$ are considered. We also fix a set of symbols of
relations $\Phi$. Finally, the subject of knowledge is a model
$(G, \Phi, f)$, where $G$ is an algebra in $\Theta$ and $f$ is a
realization of the set $\Phi$ in $G$. It is a model in the
ordinary mathematical meaning.  Similar to the previous lecture,
we write $f$ instead of $(G, \Phi, f)$ for short.  Given $\Phi, $
we denote the corresponding applied field by $\Phi \Theta$.

FOL is also oriented on the variety $\Theta$.

We assume that
every knowledge under consideration is represented by three
components:
\pmf
1) The description of the knowledge.
\newline
It is a syntactical part of knowledge, written out in the language
of the given logic.
The description reflects, what  we do want to know.
\pmf
2) The subject of the knowledge which is an object in the given
applied field, i.e., an object for which we determine knowledge.
\pmf
3) Content of the knowledge (its semantics).

It is assumed also that there is a finite type of the description
of knowledge, which is denoted by $X$. This $X$ determines
the space where the content of knowledge is calculated.

We consider knowledge as a triple of the form $(Description,
Subject, Content)$. In the knowledge base a query is the
description of the knowledge. The reply to query is content of
knowledge. Subject of the knowledge is fixed for the whole
knowledge base.

The first two components of knowledge are relatively independent, while the
third one is uniquely determined by the previous two components.
In the theory under consideration, this third component has a
geometrical nature.  In some sense it is an algebraic variety in
an affine space.
If $T$ is a description  of knowledge and $(G, \Phi, f)$ is
a subject, then $T^f$ denotes the content of knowledge.
We would like to fill the content with its own structure,
algebraic or geometric,
and we consider some elements of such structure.

On the one hand, the language is constructed according
to primary notions of algebraic geometry.  On the other hand,
 it uses algebraic logic.
 We want to underline that there are three aspects in our approach
 to knowledge representation:  logical (for knowledge
 description), algebraic (for the subject of knowledge) and
 geometric (in the content of knowledge).

The core point of the Lecture is elementary knowledge,
i.e., First Order Knowledge.  The main goal is
to construct a model to represent some non-elementary
knowledge about elementary knowledge using universal algebra approach.
For the solution of this problem we join the methods of
algebraic logic and universal algebraic geometry in logic, both defined
over an arbitrary variety of algebras $\Theta$.
Let me stress also that we make emphasis on the geometric
nature of knowledge.

 We consider categories of elementary knowledge.  Language of
 categories in the knowledge theory is a good way to organize and
 systematize primary elementary knowledge.  Morphisms in the knowledge
 category  give links between knowledge.  In particular, one
 can speak of isomorphic knowledge.  The categorical approach also
 allows us to use ideas of monada and comonada [ML].  It turns out that
 this leads to some general views on enrichment and
 computation of knowledge.  Enrichment of a structure can be
 associated with a suitable monada over a category, while the
 corresponding computation is organized by comonada.

 Let us make one more remark.  In every well described
 field of knowledge one can study the category of elementary
 knowledge, belonging to this field.  Consideration of such
 categories might be of special interest.

\subheading{2. Categories of elementary knowledge}

\subheading{2.1 The category Know$_{\Phi \Theta} (f)$}

Fix a model (subject of knowledge) $(G, \Phi, f)$.
Let us define a category of knowledge for this model and denote it
by Know$_{\Phi \Theta}(f)$.
It is a category for the given subject of knowledge.  The objects
of the category \knowf \  have the form $(X, T, A)$. Their meaning is
knowledge.
We do not fix the subject of knowledge, although it occurs here
implicitly, since it is fixed once for all.  The set $X$ is
multisorted in general.  It marks the ``place" where knowledge is
situated.
The set $X$ points also the ``place of the knowledge'' - the space of
the knowledge $Hom (W(X), G)$ while the subject of the knowledge
$(G, \Phi, f)$ is given. The set
$T$ is the description of the knowledge in the algebra
$Hal_{\empty\Phi \Theta}(X)$, and $A=T^f$, $A\subset Hom (W(X),G),$
 is the content of knowledge,
depending on $T$ and $f$.
The set $T^{ff} = A^f$ is the full description of the knowledge
$(X, T, A)$ which is a Boolean filter in Hal$_{\Phi
\Theta}(X)$.

Now about morphisms $(X, T_1, A) \to (Y, T_2, B)$.
Take $s\colon W(Y) \to W(X)$.  We have also $s\colon
\Hal_{\Phi \Theta}(Y)\to \Hal_{\Phi \Theta}(X)$.
This is a homomorphism of Boolean algebras.  The first $s$ gives
also
$$
\tilde s \colon Hom (W(X), G) \to Hom (W(Y), G).$$
As above, the first $s$ is admissible for $A$ and $B$
if $\tilde s (\nu) =\nu s \in B$ for every point $\nu\colon W(X)\to G$
in $A$.

As we know, $s$ is admissible for $A$ and $B$ if and only if for every $u \in
B^f$ we have $su \in A^f$.  This is for the second $s$, for which
we have also a homomorphism $\ol s:\Hal\em (Y)/B^f\to \Hal\em
(X)/A^f$.  It is easy to prove that $s$ is admissible for $A$ and
$B$ if and only if $s u \in A^f$ holds for every $u \in T_2$.
We consider such $s$ as a morphism
$$
s\colon (X, T_1, A) \to (Y, T_2, B)
$$
in the weak category \knowf.  We have $\tilde s (\nu) = \nu s \in
B$ if $\nu \in A$, and $s$ induces a mapping $[s]\colon A \to B$.
Simultaneously arises a mapping $s\colon T_2 \to A^f$ and a we have a
homomorphism
$$
\ol s\colon \Hal\em (Y)/B^f \to \Hal\em(X)/A^f.
$$
We have already mentioned  that
$\ol s_1 = \ol s_2$ follows from $[s_1]=[s_2]$.
Thus, we can consider the morphisms of the form
$$
[s]\colon (X, T_1, A) \to (Y, T_2, B)
$$
in the exact category \knowf.

The canonical functors

$$
{Know_{\Phi \Theta}(f)}\to K\em(f)
$$
 for weak and exact
categories are given
by the transition $(X, T, A) \to (X, A)$.
In this transition we ``forget" to fix the description of knowledge.

\subheading{2.2 The category \know}

Let us define the category of elementary knowledge for the whole
applied field $\Phi\Theta$; the subject of the knowledge
$(G, \Phi, f)$ is not fixed.  As earlier, we proceed from the
category $\Phi \Theta$ whose morphisms are homomorphisms in
$\Theta$.
They ignore the relations from $\Phi$.

Objects of the category \know \  are knowledge, represented by
$$
(X, T, A; (G, \Phi, f)),
$$
and we write $(X, T, A; G,  f)$, because $\Phi$ is fixed
for the category.
Here $X$ denotes the place of
knowledge; $A = T^f$,  $G$ and $f$ may change.

Consider morphisms:
$$
(X, T_1, A; G_1, f_1) \to (Y, T_2, B; G_2, f_2).
$$
We apply the same approach as before with some extensions.

Start from $s: W(Y) \to W(X)$ and $\delta: G_1 \to G_2$.
These $s$ and $\delta$ should correlate.
Let us explain the correlation condition.  Take a set $A_1 =
\{ \delta \nu, \nu \in A\}=\delta^* A$ and
 take further $T^\delta_1 = A_1^{f_2}$.
 Correlation of $s$ and $\delta$ means that $su\in T^\delta_1$ holds for
 any $u \in T_2$.  The same for every $u \in B^{f_2}$.
 The last also says that there is a homomorphism
 $$
 \ol s\colon \Hal\em (Y) /B^{f_2} \to \Hal\em (X)/A_1^{f_2}.
 $$
 The first of the two mappings $(s, \delta) \colon A \to B$ and
 $s\colon T_2 \to T^\delta_1$ transforms the content of knowledge,
 while the second one acts on the description.  Here $T_2$ and
 $T^\delta_1$ describe knowledge, associated with the same subject
 $(G_2, \Phi, f_2)$.

 With fixed $\delta$ there is also an exact mapping $([s], \delta): A \to
 B$.  We come to  weak and exact categories \know.  The morphisms of
 the first one are $(s, \delta)$ and in the second one they are $([s],
 \delta)$ for $(X, T_1, A; G_2, f_1) \to (Y, T_2, B;
 G_2, f_2)$.

The canonical functor \know$ \to K\em$ is defined by the
 transition
 $$
 (X, T, A; G, f) \to (X, A; G, f).
 $$
 As above, we remove description of knowledge from the information
 about it.

\subheading{2.3 Categories $K\em(G)$ and \know$(G)$}

Algebra $G \in \Theta$ is fixed in the categories $K\em (G)$ and
\know$(G)$.  A set of symbols of relations $\Phi $ is fixed
as usual, but realizations $f$ of $\Phi$ in $G$ may change.
Thus, $K\em(G)$ is a subcategory in $K\em$ and \know$(G)$ is a
subcategory in \know.   Here the corresponding $\delta: G \to G$ are
identical homomorphisms.  Objects of the category $K\em(G)$  are
now recorded as $(X, A, f)$, and those of the category \know$(G)$
as $(X, T, A, f)$.  There is a canonical functor \know$(G) \to
K\em(G)$.  As for morphisms
$$
\eqalign{
&(X, A, f_1) \to (Y, B, f_2) \; \; \hbox{\rm and} \cr
&(X, T_1, A, f_1) \to (Y, T_2, B, f_2),\cr}
$$
we can note that $A=A_1, A_1^{f_2} = T^\delta_1$ and
$A^{f_2}=T_1^{f_1f_2}$.
Hence, the corresponding admissible $s: W(Y) \to W(X)$ transfers
each $u \in T_2$ into $s u \in T_1^{f_1 f_2}$ and it induces a
homomorphism
$$
\ol s: \Hal\em (Y) /B^{f_2} \to \Hal\em (X) /A^{f_2}.
$$
Every $s$ gives a mapping $[s]: A \to B$.

\heading 3. Databases and Knowledge bases
\endheading

\subheading{3.1 Databases}

The proposed model of a database differs from that of [Pl 1].  We
want to compare databases and knowledge bases.  Geometrical aspect
in databases reflects the fact
 that the reply to the query can be considered as
an algebraic variety. It is a simple algebraic variety.  In
knowledge bases we deal with arbitrary algebraic varieties.  But
this is not the only difference.

Database is represented as a category.

Let us fix an algebra $G \in \Theta$ and consider an (admissible)
set $F$ of realizations $f$ of the set of symbols of relations
$\Phi$ in $G$.  These $f$ are instances of a database.
For every instance $f \in F$ we have a model $(G, \Phi, f)$,
all of them forming a system of models, denoted by $(G, \Phi,
F)$.  We call this system  a {\it multimodel}.

Consider a $DB$ as a category whose  objects have the form
$$
\Val_f\colon \Hal\em \to R_f, \ \ \ f \in F.
$$
Here, as above, $R_f $ is a  subalgebra in the algebra $\Hal_\Theta (G)$,
coinciding with the image of the homomorphism $\Val_f$.
If $u$ is a query to a database, then the reply to this query in
the instance $f$ is $u\ast f = Val_f (u)$.

The morphisms are
homomorphisms $\gamma: R_{f_1} \to R_{f_2}$  with the commutative
diagrams

$$
\CD
\Hal_{\Phi \Theta} @>\Val_{f_1}>> R_{f_1} \\
@. @/SE/ \Val_{f_2}// @VV\gamma V\\
@. R_{f_2}\\
\endCD
$$
We call these $\gamma$ correct homomorphisms.  The diagram
associates replies to the same query in different instances and for
different $\gamma$.

Note that all algebras $R_f$ are simple algebras and hence all
$\gamma : R_{f_1} \to R_{f_2}$ are injective.

Denote the database by $DB(G, \Phi, F)$.

\subheading{3.2 Knowledge bases}

We fix again a multimodel $(G, \Phi, F)$ and consider a knowledge
base $KB(G, \Phi, F)$.  This knowledge base is a category, whose
objects are:
$$
\Val_f\colon \Hal\em \to K\em (f), f \in F.
$$
The mapping $\Val_f$ transforms formulas of the algebra $\Hal\em$
into the objects of the category $K\em(f)$, which is a subcategory
in $K\em(G)$.  Denote by $R\em (f)$ a full subcategory in
$K\em(f)$, whose objects form a subalgebra $R_f$ in $\Hal_\Theta
(G)$.

In each object for every description of knowledge $T$ the content
of knowledge $A=T^f$ is calculated. This $A$ is considered as an
object of the category $K\em(f)$ with all its internal and
external ties in this category.

Morphisms of the category $KB(G, \Phi, F)$ are represented as
follows:

$$
\CD
\Hal_{\Phi \Theta} @>\Val_{f_1}>> K\em(f_1) \\
@. @/SE/ \Val_{f_2}// @VV\gamma V\\
@. K\em(f_2)\\
\endCD
$$
This diagram needs explanation.  Here, $\gamma$ is a functor of
categories and commutativity of the diagram is supposed on the
level of the objects of categories who are elements of the algebra
$\Hal_\Theta (G)$.  This $\gamma$ induces the diagram

$$
\CD
\Hal_{\Phi \Theta} @>\Val_{f_1}>> R_{f_1} \\
@. @/SE/ \Val_{f_2}// @VV\gamma V\\
@. R_{f_2}\\
\endCD
$$

Thus, there is a canonical functor
$$
KB(G, \Phi, f) \to DB(G, \Phi, f).
$$
This functor shows what is in common for
 databases and knowledge bases. A homomorphism
of algebras $R_{f_1} \to R_{f_2}$ corresponds to the functor
$\gamma : K_{\Phi\Theta}(f_1) \to  K_{\Phi\Theta}(f_2)$.
Furthermore, we assume that if $T$ is a set of formulas in
$\Hal\em (X)$, then
$$
\gamma (T^{f_1}) = T^{f_2}.
$$
This is a strengthened commutativity of the diagram.  It connects
knowledge content for the same description in different instances
and for different $\gamma$.

In the knowledge category \know$(G)$ we can distinguish a
subcategory for the given set of instances $F$.

Note also the ties between knowledge base and knowledge category.

We consider the following commutative diagram:

$$
\CD
\Hal_\Theta (\Phi) @[2]>\Val_f>> K_{\Phi\Theta}(f) \\
@[2]/NW/// @.@.  @/NE// /\\
@.\!\! Know_{\Phi\Theta}(f)\\
\endCD
$$

for every $f\in F.$

The right arrow in the object $(X,T,A)$ ``forgets'' the component
$T$, while the left one ``forgets'' the component $A$. Such
diagrams associate the category of knowledge with the knowledge
base.

\heading 4. Equivalence of databases and knowledge bases
\endheading

\subheading{4.1. Equivalence of databases}

Let two databases with different $(G_1, \Phi_1, F_1)$ and
$(G_2, \Phi_2, F_2)$ be given.  We are interested in
informational equivalence of these databases; another approach see
in [PT].

Consider pairs $(\a, \gamma)$ where $\a: F_1 \to F_2$ is a mapping of
sets and $\gamma$ is a function, defining a homomorphism $\gamma_f
\colon R_f \to R_{f^\a}$ for every $f \in F_1$.

The pair $(\a, \gamma)$ is called an {\it equivalence} of the
corresponding databases if $\a$ is a bijection and every
$\gamma_f, f \in F_1$, is an isomorphism of algebras.
Databases are equivalent, if there exists an equivalence
$(\a, \gamma)$ between them.

Let us motivate this definition.  Take first a kernel $\Ker(\Val_
f)$ for every $\Val_f\colon$ $\Hal\em \to R_f$ and pass to a
factor algebra $Q_f = \Hal\em / \Ker(\Val_f)$.
Let here $\delta_f \colon \Hal\em \to Q_f $ be a natural
homomorphism.
Represent a homomorphism $\Val_f$ as $\Val_f =
\Val^\circ_f\delta_f$ where $\Val_f^\circ\colon Q_f \to R_f$ is an
isomorphism.  We have a diagram of isomorphisms for $f \in F_1$

$$
\CD
Q_f @>\beta^\circ_f>> Q_{f^\a} \\
@V\Val^\circ_f VV @VV\Val^\circ_{f^\a} V\\
R_f @>\gamma_f >> R_{f^\a}\\
\endCD
$$

Here $\b^\circ_f$ =
$(\Val^\circ_{f^\a})^{-1} \gamma_f \Val_f^\circ.$

Along with the natural homomorphism $\delta_f$ we fix also some
function of choice $\delta^{-1}_f\colon Q_f \to \Hal\em$ which chooses a
definite $u \in \Hal\em$ with $\delta_f (u) = q$ for every $q \in
Q_f$.

Consider  special functions $\beta$ and $\beta'$.  The
function $\beta$ gives a mapping $\beta_f\colon$
$\Hal_{\Phi_1\Theta} \to \Hal_{\Phi_2 \Theta}$ for every
$f \in F_1$.
First we define a homomorphism $\b^1_f \colon \Hal_{\Phi_1\Theta}
\to Q_{f^\a}$ by $\beta^1_f = \beta^\circ_f \delta_f$.
Now, $\beta_f = \delta^{-1}_{f^\a} \beta^1_f = \delta^{-1}_{f^\a}
\beta^\circ_f \delta_f$.
Here $\b_f$ is a multisorted mapping of algebras, which is not,
obviously, a homomorphism.

The function $\b'$ chooses a mapping $\b'_f\colon
\Hal_{\Phi_2\Theta} \to \Hal_{\Phi_1 \Theta}$ for
every
$f^\a \in F_2$.  It is constructed similarly:
$$
\b'_f=\delta^{-1}_f(\b^\circ_f)^{-1} \delta_{f^\a}.
$$
Let now $u \in \Hal_{\Phi_1 \Theta}$ and $f \in F_1$.  Then
$$
\eqalign{
(&u*f)^{\gamma_f} = \gamma_f (\Val_f(u)) = \gamma_f (\Val^\circ_f
\delta_f (u))=\cr
&=\, \Val^\circ_{f^\a} (\Val^\circ_{f^\a})^{-1} \gamma_f \Val_f^\circ
\delta_f (u) = \Val^\circ_{f^\a} \b^\circ_f \delta_f(u)=\cr
&=\, \Val^\circ_{f^\a} \delta_{f^\a} \delta^{-1}_{f^\a}
\b^\circ_f\delta_f (u) = \Val_{f^\a} \b_f (u) = u^{\b_f} *
f^\a.\cr}
$$
Analogously, if $u \in \Hal_{\Phi_2 \Theta}$ and $f_1 = f^\a
\in F_2$, then
$$
(u*f^\a)^{\gamma_f^{-1}} = u^{\b'_f} * f.
$$
Thus, the reply to the query in the first database can be obtained
via the second one and vice versa.

Consider separately a case when $\Phi_1 = \Phi_2 = \Phi$
and every isomorphism $\gamma_f$ is correct.  The last means that
$$
\CD
\Hal_{\Phi \Theta} @>\Val_{f}>> R_{f} \\
 @. @/SE/\Val_{f^{\a}}// @VV\gamma_f V\\
 @. R_{f^{\a}}\\
\endCD
$$
holds for every $f \in F_1$.  Also, for every query $u \in
\Hal\em$ and every $f \in F_1$ we have $(f*u)^{\gamma_f} =
u*f^\a$.

Take a morphism in the first database:
$$
\CD
\Hal_{\Phi \Theta} @>\Val_{f_1}>> R_{f_1} \\
 @. @/SE/\Val_{f_2}// @VV\gamma V\\
 @. R_{f_2}\\
\endCD
$$
and construct the morphism in the second database:
$$
\CD
\Hal_{\Phi \Theta} @>\Val_{f_1^\a}>> R^\a_{f_1} \\
 @. @/SE/\Val_{f_2^{\a}}// @VV\gamma^\a V\\
 @. R_{f_2^{\a}}\\
\endCD
$$
Here, $\gamma^\a= \gamma_{f_2} \gamma \gamma^{-1}_{f_1}$
can be found from the diagram
$$
\CD
R_{f_1} @>\gamma_{f_1}>> R_{f_1^{\a}} \\
@V\gamma VV @VV\gamma^{\a} V\\
R_{f_2} @>\gamma_{f_2} >> R_{f_2^{\a}}\\
\endCD
$$

We need to check that $\gamma^\a$ is correct: $\gamma^\a
\Val_{f_1^\a} = \Val_{f_2^\a}$.  We have: $\gamma^\a \Val_{f_1^ \a} =
\gamma_{f_2} \gamma \gamma^{-1}_{f_1} \gamma_{f_1}\Val_{f_1} =
\gamma_{f_2} \gamma \Val_{f_1} = \gamma_{f_2} \Val_{f_2} =
\Val_{f_2^\a}$. Hence, the transition $\gamma \to \gamma^\a$ gives
the functor which is an isomorphism of databases, and in this
specific situation equivalence of databases turns out to be their
isomorphism.

Let us return to the general case.  The relations with the special
functions of choice $\beta$ and $\beta'$ for the sets
$F_1$ and $F_2$ can be now represented as commutative diagrams:
$$
\CD
\Hal_{\Phi_1\Theta} @>\Val_f>> R_f \\
@V\beta_f VV @VV\gamma_f V\\
\Hal_{\Phi_2 \Theta} @>\Val_{f^\a}  >> R_{f^\a}\\
\endCD
\qquad
\CD
\Hal_{\Phi_1\Theta} @>\Val_f>> R_f \\
@A\beta'_f AA @AA\gamma_f^{-1} A\\
\Hal_{\Phi_2 \Theta} @>\Val_{f^\a}  >> R_{f^\a}\\
\endCD
$$

These diagrams along with further remarks mean that equivalence of
databases in general can be treated as some semi-isomorphism or
skew isomorphism [Pl1]. The remarks are as follows.

Let $\gamma\colon R_{f_1} \to R_{f_2}$ be a morphism in the
first database.  Take $\gamma^\a = \gamma_{f_2} \gamma
\gamma^{-1}_{f_1}$.  Check that for every $u \in \Hal_{\Phi_1
\Theta}$ it holds
$$
\gamma^\a \Val_{f_1^\a} (\beta_{f_1}(u)) = \Val_{f_2^\a}
(\b_{f_2}(u)).$$
Indeed,
$$
\eqalign{
\gamma^\a &\Val_{f_1^\a} (\beta_{f_1} (u)) = \gamma^\a (\beta_{f_1}
(u) * f_1^\a) =\cr
&= \gamma^\a\gamma_{f_1} (u*f_1) = \gamma_{f_2} \gamma
\gamma^{-1}_{f_1} \gamma_{f_1} (u*f_1)=\cr
&=\gamma_{f_2} \gamma(u*f_1) = \gamma_{f_2} (u*f_2) =
\beta_{f_2} (u) * f^\a_2=\cr
&=\Val_{f_2^\a} (\beta_{f_2} (u)).\cr}
$$
Thus, $\gamma^\a$ is not anymore a morphism in the second database,
but some ``skew" morphism.  The same for the transition from the
second database to the first one with the mapping $\beta'$.

\subheading{4.2 Equivalence of knowledge bases}

Again we regard multimodels $(G_1, \Phi_1, F_1)$ and $(G_2,
\Phi_2, F_2)$ and the related know\-ledge\-bases, denoted by $KB_1$ and
$KB_2$.

Consider pairs $(\a, \gamma)$ where $\a: F_1\to F_2$ is a
bijection of sets and $\gamma$ is a function, determining an
isomorphism of weak categories:
$$
\gamma_f\colon K_{\Phi_1 \Theta} (f) \to K_{\Phi_2
\Theta} (f^\a).
$$
We assume that the isomorphism $\gamma_f$ induces isomorphism of
algebras $\gamma_f: R_f\to R_{f^\a}$.  It preserves the type $X$
and is coordinated with the inclusion of sets of points in affine
spaces.

Consider every such pair $(\a, \gamma)$ as an equivalence of
knowledge bases:  knowledge bases are equivalent if there is some
equivalence $(\a, \gamma)$.

As it was done for databases, let us pass to
motivations and distinguish the case when
$\Phi_1=\Phi_2=\Phi$ and isomorphism $\gamma_f $ is
correct.  The last means that
$$
\CD
\Hal\em @>\Val_f >> K\em (f) \\
@. @/SE/ \Val_{f^\a} //  @VV\gamma_f V\\
@. K\em (f^\a)
\endCD
$$
for every $f \in F_1$.

Take an arbitrary set of formulas $T$ in the definite
$\Hal\em(X)$.
Then:
$$
\eqalign{
&T^{f^\a} = \mathop{\cap}\limits_{u \in T} \Val_{f^\a} (u) = \mathop{\cap}\limits_{u
\in T} \gamma_f \Val_f (u) =\cr
&=\gamma_f(\mathop{\cap}\limits_{u \in T} \Val_f(u)) = \gamma_fT^f.\cr}
$$
We have used here the correlation of $\gamma_f $ with inclusions
which brings correlation with intersections.  Thus, all $T^{f^\a}$
and $T^f$ are well correlated by the isomorphism $\gamma_f$.

Check now that the pair $(\alpha, \gamma)$ gives an isomorphism of
categories $KB_1$ and $KB_2$.  If $\Val_f: \Hal\em \to K\em(f)$ is
an object in the first category, then the corresponding object
 of the second category is $\Val_{f^\a}: \Hal\em \to K\em(f^\a)$,
and vice versa.  Consider a morphism
$$
\CD
\Hal\em @>\Val_{f_1} >> K\em (f_1)\\
@. @ /SE/ \Val_{f_2} // @VV\gamma V\\
@. K\em(f_2)\\
\endCD
$$
The corresponding diagram is
$$
\CD
\Hal\em @>\Val^\a_{f_1} >> K\em (f^\a_1) \\
@. @ /SE/ \Val_{f_2}^\a // @VV\gamma^\a=\gamma_{f_2}
\gamma\gamma^{-1}_{f_1} V\\
@. K\em (f_2^\a)\\
\endCD
$$

This second diagram is actually a morphism  in the second
category.

Thus, in this special situation equivalence of knowledge bases is
reduced to their isomorphism. The same fact has been established
for the databases.

Now let us make some remarks on the general case.  For the given
set of formulas $T$ in $\Hal_{\Phi_1 \Theta} (X) $ and the
mapping $\beta_f: \Hal_{\Phi_1 \Theta} \to \Hal_{\Phi_2
\Theta}$ consider a set $T^{\beta_f}$ in $\Hal_{\Phi_2\Theta}
(X)$, defined by
$$
T^{\beta_f} = \{ \beta_f(u)|u \in T\}.
$$
Check, that $T^{\beta_f f^\a} = \gamma_f T^f$ for every $f \in
F_1$.
We have
$$
\eqalign{
T^{\beta_f f^\a} &=\mathop{\cap}\limits_{u \in T} \Val_{f^\a} (\beta_f(u))
= \mathop{\cap}\limits_{u \in T} (\beta_f(u) * f^\a)=\cr
&=\mathop{\cap}\limits_{u \in T} \gamma_f(u*f) = \gamma_f (\mathop{\cap}\limits_{u
\in T} \Val_f (u)) =  \gamma_f T^f.\cr}
$$
Similarly, if $T$ is a set of formulas in
$\Hal_{\Phi_2\Theta} (X)$, then
$$
T^{\beta'_f f}= \gamma^{-1 }_f T^{f^\a}.
$$
The pointed relations mean that equivalence of knowledge bases
$KB_1$ and $KB_2$ provides the corresponding informational
equivalence on the level of the transition from the description of
knowledge to its content.

Let now $\gamma \colon K_{\Phi_1 \Theta} (f_1) \to K_{\Phi_1 \Theta}
(f_2)$ be a morphism in $KB_1$.  Take $\gamma^\a = \gamma_{f_2}
\gamma \gamma^{-1}_{f_1}$.
One can check that
$\gamma^\a T^{\beta_{f^1} f^\a_1} = T^{\beta_{f_2} f^\a_2}$ for
every $T$ in $\Hal_{\Phi_1 \Theta} (X)$.
In particular, it means that $\gamma^\a: K_{\Phi_2 \Theta}
(f^\a_1) \to K_{\Phi_2\Theta} (f^\a_2)$ is not a morphism in
$KB_2$. Also here we have some ``skew" property which needs additional
motivation [Pl1].

\subheading{4.3. Main results}

 Two
{\it multimodels} $(G_1, \Phi, F_1)$ and $(G_2, \Phi, F_2)$ with
the same $\Phi$ are called {\it geometrically equivalent} if there
is a bijection $\a: F_1 \to F_2$ such that for every $f \in F_1$
models $(G_1, \Phi, f)$ and $(G_2, \Phi, f^\a)$ are geometrically
equivalent. \proclaim{Theorem 1}  If multimodels $(G_1, \Phi,
F_1)$ and $(G_2, \Phi, F_2)$ are geometrically equivalent, then
the knowledge bases $KB_1$ and $KB_2$ are isomorphic and, hence,
equivalent.
\endproclaim

Let us say that two {\it multimodels} $(G_1, \Phi,  F_1)$ and $(G_2,
\Phi, F_2)$ are {\it isomorphic} if there is a bijection $\a: F_1 \to F_2$
such that for every $f \in F_1 $ the models $(G_1, \Phi, f)$ and
$(G_2, \Phi, f^\a)$ are isomorphic.  If multimodels are isomorphic,
then they are geometrically equivalent and the corresponding
bases $KB_1 $ and $KB_2$ are isomorphic.  If here the algebras $G_1$
and $G_2$ are finite, then the opposite is true as well.

\proclaim{Theorem 2} If algebras $G_1$ and $G_2$ in $(G_1, \Phi,
F_1)$ and $(G_2, \Phi, F_2)$ are finite, then the knowledge bases
$KB_1 $ and $KB_2$ are isomorphic if and only if the multimodels
are isomorphic.
\endproclaim

This theorem gives the algorithm of
verification of isomorphism of two finite knowledge bases.

Consider now the question of equivalence of two finite knowledge
bases.  Here we need additional definitions.

\definition{Definition 1} Let $(G_1, \Phi_1, f_1)$ and $(G_2, \Phi_2,
f_2)$ be two models ($\Phi_1$ and $\Phi_2$ may be different).
We call them automorphic equivalent, if
\item{1)} Algebras $G_1$ and $G_2$ are isomorphic
\item{2)} Groups $\Aut(f_1) $ and $Aut(f_2)$ are conjugated by
some isomorphism of algebras $G_1$ and $G_2$.
\enddefinition

In other words, there exists an isomorphism $\delta:
G_2 \to G_1$ such that
$$
\Aut (f_2) = \delta^{-1} \Aut (f_1) \delta.
$$
\definition{Definition 2}  Two multimodels $(G_1, \Phi_1, F_1)$ and $(G_2,
\Phi_2, F_2)$ are automorphic equivalent, if
for some bijection $\a: F_1 \to F_2$ the models
$(G_1, \Phi_1, f)$ and $(
G_2, \Phi_2, f^\alpha)$ are automorphic equivalent for every $f\in F_1$.
\enddefinition

\proclaim{Proposition 1} If the multimodels $(G_1, \Phi_1, F_1)$
and $(G_2, \Phi_2, F_2)$ are automorphic equivalent then the
corresponding knowledge bases $KB_1$ and $KB_2$ are equivalent.
\endproclaim

Let now the multimodels $(G_1, \Phi_1, F_1)$ and $(G_2, \Phi_2,
F_2)$ with finite $G_1, G_2 \in \Theta $ be given; $KB_1$ and
$KB_2$ are the corresponding knowledge bases.

\proclaim{Theorem 3}  Knowledge bases $KB_1$ and $KB_2$ are
equivalent if and only if their multimodels are automorphic
equivalent.
\endproclaim

The proof of the theorem is based on Galois theory from the
previous lecture.

It also gives the algorithm of verification of two finite
knowledge bases (compare [PT]).

Let us note that it is natural to consider semi-isomorphisms along
with the isomorphisms of knowledge bases with the same $\Phi$.
Semi-isomorphisms are described by the diagrams of the form
$$
\CD
\Hal\em @>\Val_f>> K\em(f)\\
@V\sigma VV @VV\gamma_f V\\
\Hal\em @>\Val_{f^\a} >> K\em(f^\a) \\
\endCD
$$
where $\sigma$ is an automorphism of the algebra $\Hal\em$ and
$\gamma_f$ is an isomorphism of the categories.
According to such a diagram, we have
$$
\gamma_f\Val_f(u) = (u*f)^{\gamma_f} = \Val_{f^\a} (\sigma u) =
\sigma u *f^\a
$$
for every $u \in \Hal\em$.  Thus, we may replace the mappings
$\beta_f$ and $\beta_f'$ in the general situation by universal
mappings $\sigma$ and $\sigma^{-1}$.  In this case all the
pictures with ``skew" become more visible. We see that
semi-isomorphism implies equivalence of knowledge bases, as well
as isomorphism does.

We can also consider a category of knowledge bases for fixed
$\Phi$ and $\Theta$.  These are categories with usual (special)
morphisms, and categories with semimorphisms.  We can consider
monads and comonads in them.

\heading
5. Conclusion
\endheading

The main problem in computation of knowledge
is to find the content of knowledge $A=T^f$ by
the given description of knowledge $T$.  Since $A =
\mathop{\cap}\limits_{u\in T}\Val_f(u)$, we need to compute the
sets $\Val_f (u) $ for various $u \in \Hal_{\Phi\Theta}$.

Pass from the algebra $\Hal_{\Phi\Theta}$ to the corresponding
free Halmos algebra $\tilde H =\tilde H\em$ (see Lecture 6).  We
treat every $u $ as an element in $\tilde H$.  We consider an
algebra $\tilde H$ to be the constructive one.  In other words,
every element $u$ is well represented by atoms.  Thus, the
computation of the sets $\Val_f (u)$ with arbitrary $u$ is reduced
to that for atoms. The last is somehow made for constructive
models $(G, \Phi, f)$. All above can be realized via some comonada
in the category of knowledge bases.  One may also consider monads
in these categories for enrichment of knowledge.

Let us note the following general problem: what is the inference of
knowledge in the categories of knowledge under consideration?
The problem is to well formalize  the corresponding idea.


\newpage
\centerline{BIBLIOGRAPHY}
\bigskip

\item{[B]} H.Bauer, About Hilbert and Ruckkert Nullstellensatz, Manuscript.

\item{[BelP]} A.Belov, B.Plotkin, Abstract and natural Zariski topology, Manuscript.

\item{[BMR1]} G.Baumslag, A.Myasnikov, V.Remeslennikov, Algebraic geometry
over groups, J. Algebra, {\bf 219} (1999),  16 -- 79.

\item{[BMR2]} G.Baumslag, A.Myasnikov, V.Remeslennikov, Algebraic geometry
over groups, in book ``Algorithmic problems in groups and
semigroups'', Birkhauser, 1999, p. 35 -- 51.

\item{[BMRo]} G.Baumslag, A.Myasnikov, V.Roman'kov, Two theorems about
equationally noetherian groups, J.~Algebra, {\bf 194} (1997), 654
-- 664.

\item{[Be1]} A. Berzins, Geometrical equivalence of algebras,
International Journal of Algebra and Computations, {\bf 11:4}
(2001), 447 -- 456.

\item{[Be2]} A. Berzins, Variety Grp-F is semiperfect, Preprint,
Riga, 1998.

\item{BPP} A. Berzins, B. Plotkin, E. Plotkin, Algebraic
geometry in varieties of algebras with the given algebra of
constants,
 Journal of Math. Sciences, {\bf 102:3}, (2000), 4039 -- 4070.

\item{[BlG]} V.Bludov, D.Gusev, On geometric equivalence of
groups, (2002), to appear

\item{[Br]} R.Bryant, The verbal topology of a group, J.~Algebra, {\bf
48} (1977), 340 -- 346.

\item{[Co]} P.Cohn, Universal algebra,  Harper \& Row Publ., New York,
1965.

\item{[D]} Y.Y. Diers,  Categories of algebraic sets, Applied categorical
structures, {\bf 4} (1996) 329 -- 341.

\item{[DF]} J. Dyer, E. Formanek,  The automorphism group of a 
free group is complete. J. London Math. Soc. (2), {\bf{11:2}} 
(1975),  181 -- 190.

\item{[For]} E.Formanek, A question of B.Plotkin about semigroup
of endomorphisms of a free group, Proc. Amer. Math.Soc., {\bf 130}
(2002), 935 -- 937.

\item{[F]} L.Fuchs, Infinite abelian groups, Academic Press, New York,
London, 1970.

\item{[Gro]} A.Grothendieck, The cohomology theory of abstract algebraic
varieties, Proc. Int. Congress Mathematik, Edinburgh 1958,
Cambridge Univ. Press., 1960, 103 -- 118.

\item{[GoSh]} R.Gobel, S. Shelah. Radicals and Plotkin's problem
concerning geometrically equivalent groups. Proc. Amer. Math.Soc.,
{\bf 130} (2002), 673 -- 674.

\item{[Go]} V.A. Gorbunov, Algebraic theory of quasivarieties,
Doctoral Thesis, Novosibirsk, 1996, Plenum Publ. Co., 1998.

\item{[Gu]} V.Guba, Equivalence of infinite systems of equations in
free groups and semigroups to finite systems, Mat. Zametki, {\bf
40:3} (1986), 321 -- 324.

\item{[GK]} R.I. Grigorchuk, P.F. Kurchanov, Some questions in group
theory related to geometry, in book ``Combinatorial group theory
and applications to geometry'', Springer, 1998,  169 -- 232.

\item{[Hal]} P.Halmos, Algebraic logic, New York, 1962.

\item{[Har]} R.Hartshorne, Algebraic geometry, Springer-Verlag, New-York,
Heidelberg, 1977.

\item{[HMT]} L. Henkin, J.Monk, A,Tarski, Cylindric algebras, North Holland,
1985

\item{[Hi]} P.Higgins, Groups with multioperators, Proc. London Math. Soc.,
6, (1956), 366 -- 373.

\item{[HZ]} E.Hrushovski, B.Zilber, Zariski geometries, J. Amer. Math. Soc,
{\bf 1}, (1996), 1 -- 56.

\item {[KM]}
O.Kharlampovich, A. Myasnikov  Irreducible affine varieties over a
free group. I: irreducibility of quadratic equations and
Nullstellensatz J. of Algebra, {\bf 200: 2}, (1998) 472 -- 516.

\item {[KM1]}
O.Kharlampovich, A. Myasnikov  Irreducible affine varieties over a
free group. II:  J. of Algebra, {\bf 200: 2}, (1998) 517 -- 570.

\item{[KP]} E.Koublanova, B. Plotkin, Algebraic geometry in
group representations, Pre\-print. Jerusalem, 1999.

\item{[KR]} M.I.Krasner, Generalization et analogues de la theorie
de Galois, Comptes Rendus de congress de la Victorie de l'Ass.
Frac.pour I'Avancem. Sci., 1945, 54 -- 58.

\item{[Ku]} A.Kurosh, Lectures on general algebra, Moscow, Nauka, 1973.

\item{[Lav]} F. Lawvere, Continiously variable sets; algebraic geometry =
geometric logic, Logic Colloquium, Bristol  1973, North Holland,
Amsterdam, 1975, p. 135 -- 156.

\item{[L]} S. Lang, Algebra, Adison-Wesley Publ., Reading, MA, 1965.

\item{[Li]}  A.Lichtman, Two-generator division algebras and
Plotkin's question, to appear

\item{[L1]} R.C.Lyndon. Equations in groups, Bol.Soc. Brasil. Mat.,
{\bf 11} (1980), 79 -- 102.

\item{[L2]} R.C.Lyndon. Groups with parametric exponents, Trans.
Amer. Math. Soc., {\bf 9:6} (1960),  518 -- 533.

\item{[LS]} R.C.Lyndon, P.E. Shupp, Combinatorial group theory,
 Springer, 1977.

\item{[Ma]} G.Makanin, Decidability of the universal and positive theories
of a free group, Izv. Math. USSR., {\bf 25 :1} (1985), 75 -- 88.

\item{[Mal1]} A.I. Malcev, Algebraic systems, North Holland, 1973.

\item{[Mal2]} A.I. Malcev, Some remarks on quasivarieties of algebraic
structures, Algebra and Logic,  {\bf 5:3} (1966) 3 -- 9.

\item{[Man1]} J.I.Manin, Quantum groups and noncommutative
geometry, Publ. du C.R.M., University de Montreal, 1988.

\item{[Man2]} J.I.Manin, Topics in noncommutative geometry,
Princeton Univ. Press, New Jersey, 1991.

\item{MPP} A. Mashevitzky, B.Plotkin. E.Plotkin   Automorphisms
of categories of free algebras of varieties, Electronic Research
Announcements of AMS, {\bf 8} (2002), 1 -- 10.

\item{[Me]} Ju. Merzljakov, Positive formulas on free groups,
Algebra and Logica, {\bf 5:4} (1966), 25 -- 42.

\item{[ML]} S. MacLane, Categories for the working mathematicians,
Springer, 1971.

\item{MR1]} A.Myasnikov, V.Remeslennikov, Algebraic geometry
over groups I,  J. of Algebra, {\bf 219:1} (1999) 16 -- 79.

\item{MR2]} A.Myasnikov, V.Remeslennikov, Algebraic geometry
over groups II, Logical foundations J. of Algebra, {\bf 234:1}
(2000) 225 -- 276.

\item{[OMe]} O. O'Meara, Lectures on linear groups, Providence, 1974.

\item{[Ne1]} B.Neumann, Adjunction of elements to groups, J. London Math.
Soc., {\bf 18} (1943), 4 -- 11.

\item{[Ne2]} B.Neumann, A note on algebraically closed groups, J. London Math.
Soc., {\bf 27} (1952), 247 -- 249.

\item{[NP]} D. Nikolova, B. Plotkin, Some notes on universal algebraic
geometry, in book ``Algebra. Proc. International Conf. on Algebra
on the Occasion of the 90th Birthday of A. G. Kurosh, Moscow,
Russia, 1998'' Walter De Gruyter Publ., Berlin, 1999,  237 -- 261.

\item{[Pl1]} B. Plotkin, Universal algebra, algebraic logic and databases
       Kluwer Acad. Publ., 1994.

\item{[Pl2]} B. Plotkin,   Algebra, categories and databases , in "Handbook of
       algebra", v.2, Elsevier Publ., (2000), 81 -- 148.

\item{[Pl3]} B. Plotkin,  Algebraic logic, varieties of
algebras and algebraic varieties, in Proc. Int. Alg. Conf., St.  Petersburg,
1995,  St.Petersburg,  1999, p. 189 - 271.

\item{[Pl4]} B. Plotkin, Varieties of algebras and algebraic varieties
       Israel. Math. Journal, {\bf 96:2} (1996),  511 -- 522.

\item{[Pl5]} B. Plotkin, Varieties of algebras and algebraic varieties.
Categories of  algebraic varieties. Siberian Advanced Mathematics,
       Allerton Press, {\bf 7:2} (1997), p.64 -- 97.

\item{[Pl6]} B. Plotkin, Some notions of algebraic geometry in
universal algebra,  Algebra and Analysis, {\bf 9:4} (1997), 224 --
248, St.Peterburg Math. J., {\bf 9:4}, (1998)  859 -- 879.

\item{[Pl7]} B. Plotkin, Groups of automorphisms of algebraic
systems, Walter Nordhoff, 1972.

\item{[Pl8]} B. Plotkin, Zero divisors in group-based algebras.
Algebras without zero divisors, Bul. Acad. Sci. Mold. Mat., {\bf
2}, (1999), 67 -- 84.

\item{[Pl9]} B. Plotkin, Algebraic geometry in first order logic,
Preprint, Jerusalem, 2002.

\item{[Pl10]} B. Plotkin, Infinitary quasi-identities and infinitary
quasivarieties, Proc. Latvian Acad. Sci., Section B, {\bf
56}(2002), to appear.

\item{[Pl11]} B. Plotkin, Algebras with the same algebraic
geometry, Proceedings of the International Conference on Mathematical
 Logic, Algebra and Set Theory,  dedicated to 100 anniversary 
of P.S.Novikov, Proceedings of MIAN,  (2002), to appear

\item{[PPT]} B.Plotkin, E.Plotkin, A.Tsurkov, "Geometrical equivalence of
groups", Communications in Algebra, {\bf 27:8} (1999), 4015 --
4025.

\item{[Ra]} A.Razborov, On systems of equations in free groups, in
``Combinatorial and Geometric Group Theory, Edinburgh, 1993'', pp.
269 -- 283, Cambridge Univ. Press, Cambridge, 1995.

\item{[R]} E.Rips, Cyclic splittings of the finitely
presented groups and the canonical JSJ decomposition, in
Proceedings of Int. Math. Congress in Zurich, v1. Birkhauser,
(1995), 595 -- 600.

\item {[RS]} E.Rips, Z.Sela, Cyclic splittings of the finitely
presented groups and the canonical JSJ decomposition, Ann. of
Math., {\bf 146:1} (1997), 53 -- 109.

\item{[Ro1]} A.L.Rosenberg, Noncommutative schemes, Compositio Mathematica,
{\bf 112} (1998), 93 -- 125.

\item{[Ro2]} A.L.Rosenberg, Noncommutative algebraic geometry and
representations of quantized algebras, Kluwer Acad. Publ., 1995

\item{[Se]} Z.Sela, Diophantine geometry over groups I, IHES, {\bf 93},
(2001), 31 -- 105.

\item{[S]} W.R. Scott. Algebraically closed groups, Proc. American
Math. Soc., {\bf 2} (1951),  118 -- 121.

\item{[Sha]} I.Shapharevich, The Foundations of Algebraic Geometry, Nauka,
Moscow, 1972.

\item{[TPP]} T.Plotkin, B.Plotkin,
 Geometrical aspect of databases and knowledge bases, Algebra
 Universalis, {\bf 46} (2001) 131 -- 161.

\item{[TP]} T.Plotkin, Relational databases equivalence problem,
in book Advances of databases and information systems, Springer,
1996,  391 -- 404.

\item{[Shi]} Jan Shi-Jain, Linear groups over a ring, Chinese {\bf 7:2}
(1965),  163 -- 179.

\item{[SZ]} O.Zariski, P.Samuel, Commutative algebra, Van-Nostrand, Toronto,
1958.

\bye